
  
  %
  %

  %
  %
  \let\footnotea=\footnote
  \def\anote#1#2{\footnotea{\hbox{$^{#1}$}}{\eightpoint#2}}  
  \catcode`@=12 

 \def\defrefnote#1{\definexref{#1}{{\the\footnotenumber}}{refnotes}}

  %
  %


\ifx\eplain\undefined
  \let\next\relax
\else
  \expandafter\let\expandafter\next\csname endinput\endcsname
\fi
\next
\begingroup
  \catcode44 12 
  \catcode45 12 
  \catcode46 12 
  \catcode58 12 
  \catcode64 11 
  \expandafter\let\expandafter\x\csname ver@ifpdf.sty\endcsname
  \ifcase 0%
    \ifx\x\relax 
    \else
      \ifx\x\empty 
      \else
        1%
      \fi
    \fi
  \else
    \catcode35 6 
    \catcode123 1 
    \catcode125 2 
    \expandafter\ifx\csname PackageInfo\endcsname\relax
      \def\x#1#2{%
        \immediate\write-1{Package #1 Info: #2.}%
      }%
    \else
      \def\x#1#2{\PackageInfo{#1}{#2, stopped}}%
    \fi
    \x{ifpdf}{The package is already loaded}%
    \endgroup
  \fi
\endgroup
\begingroup
  \catcode35 6 
  \catcode40 12 
  \catcode41 12 
  \catcode44 12 
  \catcode45 12 
  \catcode46 12 
  \catcode47 12 
  \catcode58 12 
  \catcode64 11 
  \catcode123 1 
  \catcode125 2 
  \expandafter\ifx\csname ProvidesPackage\endcsname\relax
    \def\x#1#2#3[#4]{\endgroup
      \xdef#1{#4}%
    }%
  \else
    \def\x#1#2[#3]{\endgroup
      #2[{#3}]%
      \ifx#1\@undefined
        \xdef#1{#3}%
      \fi
      \ifx#1\relax
        \xdef#1{#3}%
      \fi
    }%
  \fi
\expandafter\x\csname ver@ifpdf.sty\endcsname
\ProvidesPackage{ifpdf}%
  [2009/04/10 v2.0 Provides the ifpdf switch (HO)]
\begingroup
  \catcode123 1 
  \catcode125 2 
  \def\x{\endgroup
    \expandafter\edef\csname ifpdf@AtEnd\endcsname{%
      \catcode35 \the\catcode35\relax
      \catcode64 \the\catcode64\relax
      \catcode123 \the\catcode123\relax
      \catcode125 \the\catcode125\relax
    }%
  }%
\x
\catcode35 6 
\catcode64 11 
\catcode123 1 
\catcode125 2 
\def\TMP@EnsureCode#1#2{%
  \edef\ifpdf@AtEnd{%
    \ifpdf@AtEnd
    \catcode#1 \the\catcode#1\relax
  }%
  \catcode#1 #2\relax
}
\TMP@EnsureCode{10}{12}
\TMP@EnsureCode{39}{12}
\TMP@EnsureCode{40}{12}
\TMP@EnsureCode{41}{12}
\TMP@EnsureCode{44}{12}
\TMP@EnsureCode{45}{12}
\TMP@EnsureCode{46}{12}
\TMP@EnsureCode{47}{12}
\TMP@EnsureCode{58}{12}
\TMP@EnsureCode{60}{12}
\TMP@EnsureCode{61}{12}
\TMP@EnsureCode{94}{7}
\TMP@EnsureCode{96}{12}
\begingroup
  \expandafter\ifx\csname ifpdf\endcsname\relax
  \else
    \edef\i/{\expandafter\string\csname ifpdf\endcsname}%
    \expandafter\ifx\csname PackageError\endcsname\relax
      \def\x#1#2{%
        \edef\z{#2}%
        \expandafter\errhelp\expandafter{\z}%
        \errmessage{Package ifpdf Error: #1}%
      }%
      \def\y{^^J}%
      \newlinechar=10 %
    \else
      \def\x#1#2{%
        \PackageError{ifpdf}{#1}{#2}%
      }%
      \def\y{\MessageBreak}%
    \fi
    \x{Name clash, \i/ is already defined}{%
      Incompatible versions of \i/ can cause problems,\y
      therefore package loading is aborted.%
    }%
    \endgroup
    \ifpdf@AtEnd
  \fi
\endgroup
\begingroup
  \expandafter\ifx\csname pdfoutput\endcsname\relax
  \else
    \def\skip#1\relax\endgroup{\csname fi\endcsname\endgroup}%
    \skip
  \fi
  \expandafter\ifx\csname directlua\endcsname\relax
    \def\skip#1\endgroup{\csname fi\endcsname\endgroup}%
    \skip
  \fi
  \expandafter\ifx\csname RequirePackage\endcsname\relax
    \input ifluatex.sty\relax
  \else
    \RequirePackage{ifluatex}[2009/04/10]%
  \fi
  \ifluatex
    \ifnum\luatexversion<36 %
    \else
      \directlua{tex.enableprimitives('ifpdf', {'pdfoutput'})}%
      \global\let\pdfoutput\ifpdfpdfoutput
    \fi
  \fi
  \relax
\endgroup
\newif\ifpdf
\begingroup\expandafter\expandafter\expandafter\endgroup
\expandafter\ifx\csname pdfoutput\endcsname\relax
\else
  \ifnum\pdfoutput<1 %
  \else
    \pdftrue
  \fi
\fi
\begingroup
  \expandafter\ifx\csname pdfoutput\endcsname\relax
  \else
    \escapechar=-1 %
    \edef\m{\meaning\pdfoutput}%
    \edef\p{%
      \string p\string d\string f%
      \string o\string u\string t\string p\string u\string t%
    }%
    \ifx\m\p
    \else
      \expandafter\ifx\csname PackageWarningNoLine\endcsname\relax
        \def\PackageWarningNoLine#1#2{%
          \immediate\write16{%
            Package `#1' Warning: #2.%
          }%
        }%
      \fi
      \PackageWarningNoLine{ifpdf}{%
        Someone has redefined \string\\pdfoutput%
      }%
    \fi
  \fi
\endgroup
\begingroup
  \expandafter\ifx\csname PackageInfo\endcsname\relax
    \def\x#1#2{%
    }%
  \else
    \let\x\PackageInfo
    \expandafter\let\csname on@line\endcsname\empty
  \fi
  \x{ifpdf}{pdfTeX in pdf mode \ifpdf\else not \fi detected}%
\endgroup
\ifpdf@AtEnd
\def\makeactive#1{\catcode`#1 = \active \ignorespaces}%
\chardef\letter = 11
\chardef\other = 12
\def\makeatletter{%
  \edef\resetatcatcode{\catcode`\noexpand\@\the\catcode`\@\relax}%
  \catcode`\@11\relax
}%
\def\makeatother{%
  \edef\resetatcatcode{\catcode`\noexpand\@\the\catcode`\@\relax}%
  \catcode`\@12\relax
}%
\edef\leftdisplays{\the\catcode`@}%
\catcode`@ = \letter
\let\@eplainoldatcode = \leftdisplays
\toksdef\toks@ii = 2
\def\uncatcodespecials{%
   \def\do##1{\catcode`##1 = \other}%
   \dospecials
}%
{%
   \makeactive\^^M %
   \long\gdef\letreturn#1{\let^^M = #1}%
}%
\let\@eattoken = \relax  
\def\eattoken{\let\@eattoken = }%
\def\gobble#1{}%
\def\gobbletwo#1#2{}%
\def\gobblethree#1#2#3{}%
\def\identity#1{#1}%
\def\@emptymarkA{\@emptymarkB} 
\def\ifempty#1{\@@ifempty #1\@emptymarkA\@emptymarkB}%
\def\@@ifempty#1#2\@emptymarkB{\ifx #1\@emptymarkA}%
\def\@gobbleminus#1{\ifx-#1\else#1\fi}%
\def\ifinteger#1{\ifcat_\ifnum9<1\@gobbleminus#1 _\else A\fi}%
\def\isinteger{TT\fi\ifinteger}%
\def\@gobblemeaning#1:->{}%
\def\sanitize{\expandafter\@gobblemeaning\meaning}%
\def\ifundefined#1{\expandafter\ifx\csname#1\endcsname\relax}%
\def\csn#1{\csname#1\endcsname}%
\def\ece#1#2{\expandafter#1\csname#2\endcsname}%
\def\expandonce{\expandafter\noexpand}%
\let\@plainwlog = \wlog
\let\wlog = \gobble
\newlinechar = `^^J
\def\loggingall{\tracingcommands\tw@\tracingstats\tw@
   \tracingpages\@ne\tracingoutput\@ne\tracinglostchars\@ne
   \tracingmacros\tw@\tracingparagraphs\@ne\tracingrestores\@ne
   \showboxbreadth\maxdimen\showboxdepth\maxdimen
}%
\def\tracingoff{\tracingonline\z@\tracingcommands\z@\tracingstats\z@
  \tracingpages\z@\tracingoutput\z@\tracinglostchars\z@
  \tracingmacros\z@\tracingparagraphs\z@\tracingrestores\z@
  \showboxbreadth5 \showboxdepth3
}%
\begingroup
  \catcode`\{ = 12 \catcode`\} = 12
  \catcode`\[ = 1 \catcode`\] = 2
  \gdef\lbracechar[{]%
  \gdef\rbracechar[}]%
  \catcode`\% = \other
  \gdef\percentchar[
\def\vpenalty{\ifhmode\par\fi \penalty}%
\def\hpenalty{\ifvmode\leavevmode\fi \penalty}%
\def\iterate{%
  \let\loop@next\relax
  \body
  \let\loop@next\iterate
  \fi
  \loop@next
}%
\def\edefappend#1#2{%
  \toks@ = \expandafter{#1}%
  \edef#1{\the\toks@ #2}%
}%
\def\allowhyphens{\nobreak\hskip\z@skip}%
\long\def\hookprepend{\@hookassign{\the\toks@ii \the\toks@}}%
\long\def\hookappend{\@hookassign{\the\toks@ \the\toks@ii}}%
\let\hookaction = \hookappend 
\long\def\@hookassign#1#2#3{%
  \expandafter\ifx\csname @#2hook\endcsname \relax
    \toks@ = {}%
  \else
    \expandafter\let\expandafter\temp \csname @#2hook\endcsname
    \toks@ = \expandafter{\temp}%
  \fi
  \toks2 = {#3}
  \ece\edef{@#2hook}{#1}%
}%
\long\def\hookactiononce#1#2{%
  \edefappend#2{\global\let\noexpand#2\relax}
  \hookaction{#1}#2%
}%
\def\hookrun#1{%
  \expandafter\ifx\csname @#1hook\endcsname \relax \else
    \def\temp{\csname @#1hook\endcsname}%
    \expandafter\temp
  \fi
}%
\def\setpropertyglobal#1#2#3{\ece\xdef{#1@p#2}{#3}}%
\def\getproperty#1#2{%
  \expandafter\ifx\csname#1@p#2\endcsname\relax
  \else \csname#1@p#2\endcsname
  \fi
}%
\ifx\@undefinedmessage\@undefined
  \def\@undefinedmessage
    {No .aux file; I won't warn you about undefined labels.}%
\fi
\edef\cite{\the\catcode`@}%
\catcode`@ = 11
\let\@oldatcatcode = \cite
\chardef\@letter = 11
\chardef\@other = 12
\def\@innerdef#1#2{\edef#1{\expandafter\noexpand\csname #2\endcsname}}%
\@innerdef\@innernewcount{newcount}%
\@innerdef\@innernewdimen{newdimen}%
\@innerdef\@innernewif{newif}%
\@innerdef\@innernewwrite{newwrite}%
\def\@gobble#1{}%
\ifx\inputlineno\@undefined
   \let\@linenumber = \empty 
\else
   \def\@linenumber{\the\inputlineno:\space}%
\fi
\def\@futurenonspacelet#1{\def\cs{#1}%
   \afterassignment\@stepone\let\@nexttoken=
}%
\begingroup 
\def\\{\global\let\@stoken= }%
\\ 
\endgroup
\def\@stepone{\expandafter\futurelet\cs\@steptwo}%
\def\@steptwo{\expandafter\ifx\cs\@stoken\let\@@next=\@stepthree
   \else\let\@@next=\@nexttoken\fi \@@next}%
\def\@stepthree{\afterassignment\@stepone\let\@@next= }%
\def\@getoptionalarg#1{%
   \let\@optionaltemp = #1%
   \let\@optionalnext = \relax
   \@futurenonspacelet\@optionalnext\@bracketcheck
}%
\def\@bracketcheck{%
   \ifx [\@optionalnext
      \expandafter\@@getoptionalarg
   \else
      \let\@optionalarg = \empty
      \expandafter\@optionaltemp
   \fi
}%
\def\@@getoptionalarg[#1]{%
   \def\@optionalarg{#1}%
   \@optionaltemp
}%
\def\@nnil{\@nil}%
\def\@fornoop#1\@@#2#3{}%
\def\@for#1:=#2\do#3{%
   \edef\@fortmp{#2}%
   \ifx\@fortmp\empty \else
      \expandafter\@forloop#2,\@nil,\@nil\@@#1{#3}%
   \fi
}%
\def\@forloop#1,#2,#3\@@#4#5{\def#4{#1}\ifx #4\@nnil \else
       #5\def#4{#2}\ifx #4\@nnil \else#5\@iforloop #3\@@#4{#5}\fi\fi
}%
\def\@iforloop#1,#2\@@#3#4{\def#3{#1}\ifx #3\@nnil
       \let\@nextwhile=\@fornoop \else
      #4\relax\let\@nextwhile=\@iforloop\fi\@nextwhile#2\@@#3{#4}%
}%
\@innernewif\if@fileexists
\def\@testfileexistence{\@getoptionalarg\@finishtestfileexistence}%
\def\@finishtestfileexistence#1{%
   \begingroup
      \def\extension{#1}%
      \immediate\openin0 =
         \ifx\@optionalarg\empty\jobname\else\@optionalarg\fi
         \ifx\extension\empty \else .#1\fi
         \space
      \ifeof 0
         \global\@fileexistsfalse
      \else
         \global\@fileexiststrue
      \fi
      \immediate\closein0
   \endgroup
}%
\toks0 = {%
\def\bibliographystyle#1{%
   \@readauxfile
   \@writeaux{\string\bibstyle{#1}}%
}%
\let\bibstyle = \@gobble
\let\bblfilebasename = \jobname
\def\bibliography#1{%
   \@readauxfile
   \@writeaux{\string\bibdata{#1}}%
   \@testfileexistence[\bblfilebasename]{bbl}%
   \if@fileexists
      \nobreak
      \@readbblfile
   \fi
}%
\let\bibdata = \@gobble
\def\nocite#1{%
   \@readauxfile
   \@writeaux{\string\citation{#1}}%
}%
\@innernewif\if@notfirstcitation
\def\cite{\@getoptionalarg\@cite}%
\def\@cite#1{%
   \let\@citenotetext = \@optionalarg
   \printcitestart
   \nocite{#1}%
   \@notfirstcitationfalse
   \@for \@citation :=#1\do
   {%
      \expandafter\@onecitation\@citation\@@
   }%
   \ifx\empty\@citenotetext\else
      \printcitenote{\@citenotetext}%
   \fi
   \printcitefinish
}%
\def\@onecitation#1\@@{%
   \if@notfirstcitation
      \printbetweencitations
   \fi
   \expandafter \ifx \csname\@citelabel{#1}\endcsname \relax
      \if@citewarning
         \message{\@linenumber Undefined citation `#1'.}%
      \fi
      \expandafter\gdef\csname\@citelabel{#1}\endcsname{%
         {\tt
            \escapechar = -1
            \nobreak\hskip0pt
            \expandafter\string\csname#1\endcsname
            \nobreak\hskip0pt
         }%
      }%
   \fi
   \printcitepreitem{#1}%
   \csname\@citelabel{#1}\endcsname
   \printcitepostitem
   \@notfirstcitationtrue
}%
\def\@citelabel#1{b@#1}%
\def\@citedef#1#2{{\expandafter\gdef\csname\@citelabel{#1}\endcsname{#2}}}%
\def\@readbblfile{%
   \ifx\@itemnum\@undefined
      \@innernewcount\@itemnum
   \fi
   \begingroup
      \ifx\begin\@undefined
         \def\begin##1##2{%
            \setbox0 = \hbox{\biblabelcontents{##2}}%
            \biblabelwidth = \wd0
         }%
         \let\end = \@gobble 
      \fi
      \@itemnum = 0
      \def\bibitem{\@getoptionalarg\@bibitem}%
      \def\@bibitem{%
         \ifx\@optionalarg\empty
            \expandafter\@numberedbibitem
         \else
            \expandafter\@alphabibitem
         \fi
      }%
      \def\@alphabibitem##1{%
         \expandafter \xdef\csname\@citelabel{##1}\endcsname {\@optionalarg}%
         \ifx\biblabelprecontents\@undefined
            \let\biblabelprecontents = \relax
         \fi
         \ifx\biblabelpostcontents\@undefined
            \let\biblabelpostcontents = \hss
         \fi
         \@finishbibitem{##1}%
      }%
      \def\@numberedbibitem##1{%
         \advance\@itemnum by 1
         \expandafter \xdef\csname\@citelabel{##1}\endcsname{\number\@itemnum}%
         \ifx\biblabelprecontents\@undefined
            \let\biblabelprecontents = \hss
         \fi
         \ifx\biblabelpostcontents\@undefined
            \let\biblabelpostcontents = \relax
         \fi
         \@finishbibitem{##1}%
      }%
      \def\@finishbibitem##1{%
         \bblitemhook{##1}%
         \biblabelprint{\csname\@citelabel{##1}\endcsname}%
         \@writeaux{\string\@citedef{##1}{\csname\@citelabel{##1}\endcsname}}%
         \ignorespaces
      }%
      \ifx\undefined\em \let\em=\bblem \fi
      \ifx\undefined\emph \let\emph=\bblemph \fi
      \ifx\undefined\mbox \let\mbox=\bblmbox \fi
      \ifx\undefined\newblock \let\newblock=\bblnewblock \fi
      \ifx\undefined\sc \let\sc=\bblsc \fi
      \ifx\undefined\textbf \let\textbf=\bbltextbf \fi
      \frenchspacing
      \clubpenalty = 4000 \widowpenalty = 4000
      \tolerance = 10000 \hfuzz = .5pt
      \everypar = {\hangindent = \biblabelwidth
                      \advance\hangindent by \biblabelextraspace}%
      \bblrm
      \parskip = 1.5ex plus .5ex minus .5ex
      \biblabelextraspace = .5em
      \bblhook
      \input \bblfilebasename.bbl
   \endgroup
}%
\@innernewdimen\biblabelwidth
\@innernewdimen\biblabelextraspace
\def\biblabelprint#1{%
   \noindent
   \hbox to \biblabelwidth{%
      \biblabelprecontents
      \biblabelcontents{#1}%
      \biblabelpostcontents
   }%
   \kern\biblabelextraspace
}%
\def\biblabelcontents#1{{\bblrm [#1]}}%
\def\bblrm{\rm}%
\def\bblem{\it}%
\def\bblemph#1{{\bblem #1\/}}
\def\textbf#1{{\bf #1}}
\def\bblmbox{\leavevmode\hbox}
\def\bblsc{\ifx\@scfont\@undefined
              \font\@scfont = cmcsc10
           \fi
           \@scfont
}%
\def\bblnewblock{\hskip .11em plus .33em minus .07em }%
\let\bblhook = \empty
\let\bblitemhook = \@gobble
\def\printcitestart{[}
\def\printcitefinish{]}
\def\printbetweencitations{, }
\let\printcitepreitem\@gobble 
\let\printcitepostitem\empty
\def\printcitenote#1{, #1}
\let\citation = \@gobble
\@innernewcount\@numparams
\ifx\newcommand\undefined
\def\newcommand#1{%
   \def\@commandname{#1}%
   \@getoptionalarg\@continuenewcommand
}%
\fi
\ifx\renewcommand\undefined
\let\renewcommand = \newcommand
\fi
\ifx\providecommand\undefined
\def\providecommand#1{%
   \def\@commandname{#1}%
   \expandafter\ifx\@commandname \@undefined
     \let\cs=\@continuenewcommand  
   \else
     \let\cs=\@gobble              
   \fi
   \@getoptionalarg\cs
}%
\fi
\def\@continuenewcommand{%
   \@numparams = \ifx\@optionalarg\empty 0\else\@optionalarg \fi \relax
   \@newcommand
}%
\def\@newcommand#1{%
   \def\@startdef{\expandafter\edef\@commandname}%
   \ifnum\@numparams=0
      \let\@paramdef = \empty
   \else
      \ifnum\@numparams>9
         \errmessage{\the\@numparams\space is too many parameters}%
      \else
         \ifnum\@numparams<0
            \errmessage{\the\@numparams\space is too few parameters}%
         \else
            \edef\@paramdef{%
               \ifcase\@numparams
                  \empty  No arguments.
               \or ####1%
               \or ####1####2%
               \or ####1####2####3%
               \or ####1####2####3####4%
               \or ####1####2####3####4####5%
               \or ####1####2####3####4####5####6%
               \or ####1####2####3####4####5####6####7%
               \or ####1####2####3####4####5####6####7####8%
               \or ####1####2####3####4####5####6####7####8####9%
               \fi
            }%
         \fi
      \fi
   \fi
   \expandafter\@startdef\@paramdef{#1}%
}%
}%
\ifx\nobibtex\@undefined \the\toks0 \fi
\def\@readauxfile{%
   \if@auxfiledone \else 
      \global\@auxfiledonetrue
      \@testfileexistence{aux}%
      \if@fileexists
         \begingroup
            \endlinechar = -1
            \catcode`@ = 11
            \input \jobname.aux
         \endgroup
      \else
         \message{\@undefinedmessage}%
         \global\@citewarningfalse
      \fi
      \immediate\openout\@auxfile = \jobname.aux
   \fi
}%
\newif\if@auxfiledone
\ifx\noauxfile\@undefined \else \@auxfiledonetrue\fi
\@innernewwrite\@auxfile
\def\@writeaux#1{\ifx\noauxfile\@undefined \write\@auxfile{#1}\fi}%
\ifx\@undefinedmessage\@undefined
   \def\@undefinedmessage{No .aux file; I won't give you warnings about
                          undefined citations.}%
\fi
\@innernewif\if@citewarning
\ifx\noauxfile\@undefined \@citewarningtrue\fi
\catcode`@ = \@oldatcatcode
\let\auxfile = \@auxfile
\let\for = \@for
\let\futurenonspacelet = \@futurenonspacelet
\def\iffileexists{\if@fileexists}%
\let\innerdef = \@innerdef
\let\innernewcount = \@innernewcount
\let\innernewdimen = \@innernewdimen
\let\innernewif = \@innernewif
\let\innernewwrite = \@innernewwrite
\let\linenumber = \@linenumber
\let\readauxfile = \@readauxfile
\let\spacesub = \@spacesub
\let\testfileexistence = \@testfileexistence
\let\writeaux = \@writeaux
\def\innerinnerdef#1{\expandafter\innerdef\csname inner#1\endcsname{#1}}%
\innerinnerdef{newbox}%
\innerinnerdef{newfam}%
\innerinnerdef{newhelp}%
\innerinnerdef{newinsert}%
\innerinnerdef{newlanguage}%
\innerinnerdef{newmuskip}%
\innerinnerdef{newread}%
\innerinnerdef{newskip}%
\innerinnerdef{newtoks}%
\def\immediatewriteaux#1{%
  \ifx\noauxfile\@undefined
    \immediate\write\@auxfile{#1}%
  \fi
}%
\def\bblitemhook#1{\gdef\@hlbblitemlabel{#1}}%
\def\biblabelprint#1{%
   \noindent
   \hbox to \biblabelwidth{%
      \hldest@impl{bib}{\@hlbblitemlabel}%
      \biblabelprecontents
      \biblabelcontents{#1}%
      \biblabelpostcontents
   }%
   \kern\biblabelextraspace
}%
\def\eplainprintcitepreitem#1{\hlstart@impl{cite}{#1}}%
\def\eplainprintcitepostitem{\hlend@impl{cite}}%
\def\printcitepreitem#1{%
  \testfileexistence[\bblfilebasename]{bbl}%
  \iffileexists
    \global\let\printcitepreitem\eplainprintcitepreitem
    \global\let\printcitepostitem\eplainprintcitepostitem
  \else
    \global\let\printcitepreitem\gobble
    \global\let\printcitepostitem\relax
  \fi
  \printcitepreitem{#1}%
}%
\def\@Nnil{\@Nil}%
\def\@Fornoop#1\@@#2#3{}%
\def\For#1:=#2\do#3{%
   \edef\@Fortmp{#2}%
   \ifx\@Fortmp\empty \else
      \expandafter\@Forloop#2,\@Nil,\@Nil\@@#1{#3}%
   \fi
}%
\def\@Forloop#1,#2,#3\@@#4#5{\@Fordef#1\@@#4\ifx #4\@Nnil \else
       #5\@Fordef#2\@@#4\ifx #4\@Nnil \else#5\@iForloop #3\@@#4{#5}\fi\fi
}%
\def\@iForloop#1,#2\@@#3#4{\@Fordef#1\@@#3\ifx #3\@Nnil
       \let\@Nextwhile=\@Fornoop \else
      #4\relax\let\@Nextwhile=\@iForloop\fi\@Nextwhile#2\@@#3{#4}%
}%
\def\@Forspc{ }%
\def\@Fordef{\futurelet\@Fortmp\@@Fordef}
\def\@@Fordef{%
  \expandafter\ifx\@Forspc\@Fortmp 
    \expandafter\@Fortrim
  \else
    \expandafter\@@@Fordef
  \fi
}%
\expandafter\def\expandafter\@Fortrim\@Forspc#1\@@{\@Fordef#1\@@}%
\def\@@@Fordef#1\@@#2{\def#2{#1}}%
\def\tmpfileextension{.tmp}%
\let\tmpfilebasename = \jobname
\ifx\eTeXversion\undefined
  \innernewwrite\eplain@tmpfile
  \def\scantokens#1{%
    \toks@={#1}%
    \immediate\openout\eplain@tmpfile=\tmpfilebasename\tmpfileextension
    \immediate\write\eplain@tmpfile{\the\toks@}%
    \immediate\closeout\eplain@tmpfile
    \input \tmpfilebasename\tmpfileextension\relax
  }%
\fi
\begingroup
   \makeactive\^^M \makeactive\ 
\gdef\obeywhitespace{%
\makeactive\^^M\def^^M{\par\futurelet\next\@finishobeyedreturn}%
\makeactive\ \let =\ %
\aftergroup\@removebox%
\futurelet\next\@finishobeywhitespace%
}%
\gdef\@finishobeywhitespace{{%
\ifx\next %
\aftergroup\@obeywhitespaceloop%
\else\ifx\next^^M%
\aftergroup\gobble%
\fi\fi}}%
\gdef\@finishobeyedreturn{%
\ifx\next^^M\vskip\blanklineskipamount\fi%
\indent%
}%
\endgroup
\def\@obeywhitespaceloop#1{\futurelet\next\@finishobeywhitespace}%
\def\@removebox{%
  \ifhmode
    \setbox0 = \lastbox
    \ifdim\wd0=\parindent
      \setbox2 = \hbox{\unhcopy0}
      \ifdim\wd2=0pt
        \ignorespaces
      \else
        \box0 
      \fi
    \else
       \box0 
    \fi
  \fi
}%
\newskip\blanklineskipamount
\blanklineskipamount = 0pt
\def\frac#1/#2{\leavevmode
   \kern.1em \raise .5ex \hbox{\the\scriptfont0 #1}%
   \kern-.1em $/$%
   \kern-.15em \lower .25ex \hbox{\the\scriptfont0 #2}%
}%
\newdimen\hruledefaultheight  \hruledefaultheight = 0.4pt
\newdimen\hruledefaultdepth   \hruledefaultdepth = 0.0pt
\newdimen\vruledefaultwidth   \vruledefaultwidth = 0.4pt
\def\ehrule{\hrule height\hruledefaultheight depth\hruledefaultdepth}%
\def\evrule{\vrule width\vruledefaultwidth}%
\ifx\sc\undefined
    \def\sc{%
      \expandafter\ifx\the\scriptfont\fam\nullfont
        \font\temp = cmr7 \temp
      \else
        \the\scriptfont\fam
      \fi
      \def\uppercasesc{\char\uccode`}%
    }%
\fi
\ifx\uppercasesc\undefined
  \let\uppercasesc = \relax
\fi
\def\TeX{T\kern-.1667em\lower.5ex\hbox{E}\kern-.125emX\spacefactor1000 }%
\ifx\AmS\undefined
    \def\AmS{{\the\textfont2 A}\kern-.1667em\lower.5ex\hbox
        {\the\textfont2 M}\kern-.125em{\the\textfont2 S}}
\fi
\ifx\AMS\undefined \let\AMS=\AmS \fi
\ifx\AmSLaTeX\undefined
    \def\AmSLaTeX{\AmS-\LaTeX}
\fi
\ifx\AMSLaTeX\undefined \let\AMSLaTeX=\AmSLaTeX \fi
\ifx\AmSTeX\undefined
    \def\AmSTeX{$\cal A$\kern-.1667em\lower.5ex\hbox{$\cal M$}%
            \kern-.125em$\cal S$-\TeX}%
\fi
\ifx\AMSTEX\undefined \let\AMSTEX=\AmSTeX \fi
\ifx\AMSTeX\undefined \let\AMSTeX=\AmSTeX \fi
\ifx\BibTeX\undefined
    \def\BibTeX{B{\sc \uppercasesc i\kern-.025em \uppercasesc b}\kern-.08em
                \TeX}%
\fi
\ifx\BIBTeX\undefined \let\BIBTeX=\BibTeX \fi
\ifx\BIBTEX\undefined \let\BIBTEX=\BibTeX \fi
\ifx\LAMSTeX\undefined
    \def\LAMSTeX{L\raise.42ex\hbox{\kern-.3em\the\scriptfont2 A}%
                 \kern-.2em\lower.376ex\hbox{\the\textfont2 M}%
                 \kern-.125em {\the\textfont2 S}-\TeX}%
\fi
\ifx\LamSTeX\undefined \let\LamSTeX=\LAMSTeX \fi
\ifx\LAmSTeX\undefined \let\LAmSTeX=\LAMSTeX \fi
\ifx\LaTeX\undefined
    \def\LaTeX{L\kern-.36em\raise.3ex\hbox{\sc \uppercasesc a}\kern-.15em\TeX}%
\fi
\ifx\LATEX\undefined \let\LATEX=\LaTeX \fi
\ifx\LaTeXe\undefined
    \def\LaTeXe{\LaTeX{}\kern.05em2$_{\textstyle\varepsilon}$}
\fi
\ifx\MF\undefined
    \ifx\manfnt\undefined
            \font\manfnt=logo10
    \fi
    \ifx\manfntsl\undefined
            \font\manfntsl=logosl10
    \fi
    \def\MF{{\ifdim\fontdimen1\font>0pt \let\manfnt = \manfntsl \fi
      {\manfnt META}\-{\manfnt FONT}}\spacefactor1000 }%
\fi
\ifx\METAFONT\undefined \let\METAFONT=\MF \fi
\ifx\SLITEX\undefined
    \def\SLITEX{S\kern-.065em L\kern-.18em\raise.32ex\hbox{i}\kern-.03em\TeX}%
\fi
\ifx\SLiTeX\undefined \let\SLiTeX=\SLITEX \fi
\ifx\SliTeX\undefined \let\SliTeX=\SLITEX \fi
\ifx\SLITeX\undefined \let\SLITeX=\SLITEX \fi
\edef\path{\the\catcode`@}%
\catcode`@ = 11
\let\@oldatcatcode = \path
\newcount \c@tcode
\newcount \c@unter
\newif \ifspecialpathdelimiters
\begingroup
\catcode `\ = 10
\gdef \passivesp@ce { }%
\catcode `\ = 13\relax%
\gdef\activesp@ce{ }%
\endgroup
\def \discretionaries 
    {\begingroup
        \c@tcodes = 13
        \discr@tionaries
    }%
\def \discr@tionaries #1
    {\def \discr@ti@naries ##1#1
         {\endgroup
          \def \discr@ti@n@ries ####1
              {\if   \noexpand ####1\noexpand #1%
                     \let \n@xt = \relax
               \else
                     \catcode `####1 = 13
                     \def ####1{\discretionary
                                  {\char `####1}{}{\char `####1}}%
                     \let \n@xt = \discr@ti@n@ries
               \fi
               \n@xt
              }%
          \def \discr@ti@n@ri@s {\discr@ti@n@ries ##1#1}%
         }%
     \discr@ti@naries
    }%
\let\pathafterhook = \relax
\def \path
    {\ifspecialpathdelimiters
        \begingroup
        \c@tcodes = 12
        \def \endp@th {\endgroup \endgroup \pathafterhook}%
     \else
        \def \endp@th {\endgroup \pathafterhook}%
     \fi
     \p@th
    }%
\def \p@th #1
    {\begingroup
        \tt
        \c@tcode = \catcode `#1
        \discr@ti@n@ri@s
        \catcode `\ = \active
        \expandafter \edef \activesp@ce {\passivesp@ce \hbox {}}%
        \catcode `#1 = \c@tcode
        \def \p@@th ##1#1
            {\leavevmode \hbox {}##1%
             \endp@th
            }%
     \p@@th
    }%
\def \c@tcodes {\afterassignment \c@tc@des \c@tcode}%
\def \c@tc@des
    {\c@unter = 0
     \loop
            \ifnum \catcode \c@unter = \c@tcode
            \else
                \catcode \c@unter = \c@tcode
            \fi
     \ifnum \c@unter < 255
            \advance \c@unter by 1
     \repeat
     \catcode `\ = 10
    }%
\catcode `\@ = \@oldatcatcode
\discretionaries |~!@$
\ifx\eTeX\undefined
  \def\eTeX{\hbox{\mathsurround=0pt $\varepsilon$-\kern-.125em\TeX}}%
\fi
\ifx\ExTeX\undefined
  \def\ExTeX{\hbox{\mathsurround=0pt
    $\textstyle\varepsilon_{\kern-0.15em\cal{X}}$\kern-.2em\TeX}}%
\fi
\def\eplain@Xe@reflect#1{%
  \ifx\reflectbox\undefined
    \errmessage{A graphics package must be loaded for \string\XeTeX}%
  \else
    \ifdim \fontdimen1\font>0pt
      \raise 1.75ex \hbox{\kern.1em\rotatebox{180}{#1}}\kern-.1em
    \else
      \reflectbox{#1}%
    \fi
  \fi
}%
\def\eplain@Xe#1{\leavevmode
  \smash{\hbox{X%
    \setbox0=\hbox{\TeX}\setbox2=\hbox{E}%
    \lower\dp0\hbox{\raise\dp2\hbox{\kern-.125em\eplain@Xe@reflect{E}}}%
    \kern-.1667em #1}}}%
\ifx\XeTeX\undefined
  \def\XeTeX{\eplain@Xe\TeX}%
\fi
\ifx\XeLaTeX\undefined
  \def\XeLaTeX{\eplain@Xe{\thinspace\LaTeX}}%
\fi
\def\blackbox{\vrule height .8ex width .6ex depth -.2ex \relax}
\def\makeblankbox#1#2{%
  \ifvoid0
    \errhelp = \@makeblankboxhelp
    \errmessage{Box 0 is void}%
  \fi
  \hbox{\lower\dp0
    \vbox{\hidehrule{#1}{#2}%
      \kern -#1
      \hbox to \wd0{\hidevrule{#1}{#2}%
        \raise\ht0\vbox to #1{}
        \lower\dp0\vtop to #1{}
        \hfil\hidevrule{#2}{#1}%
      }%
      \kern-#1\hidehrule{#2}{#1}%
    }%
  }%
}%
\newhelp\@makeblankboxhelp{Assigning to the dimensions of a void^^J%
  box has no effect.  Do `\string\setbox0=\string\null' before you^^J%
  define its dimensions.}%
\def\hidehrule#1#2{\kern-#1\hrule height#1 depth#2 \kern-#2}%
\def\hidevrule#1#2{%
  \kern-#1%
  \dimen@=#1\advance\dimen@ by #2%
  \vrule width\dimen@
  \kern-#2%
}%
\newdimen\boxitspace \boxitspace = 3pt
\long\def\boxit#1{%
  \vbox{%
    \ehrule
    \hbox{%
      \evrule
      \kern\boxitspace
      \vbox{\kern\boxitspace \parindent = 0pt #1\kern\boxitspace}%
      \kern\boxitspace
      \evrule
    }%
    \ehrule
  }%
}%
\def\numbername#1{\ifcase#1%
   zero%
   \or one%
   \or two%
   \or three%
   \or four%
   \or five%
   \or six%
   \or seven%
   \or eight%
   \or nine%
   \or ten%
   \or #1%
   \fi
}%
\let\@plainnewif = \newif
\let\@plainnewdimen = \newdimen
\let\newif = \innernewif
\let\newdimen = \innernewdimen
\edef\@eplainoldandcode{\the\catcode`& }%
\catcode`& = 11
\toks0 = {%
\edef\thinlines{\the\catcode`@ }%
\catcode`@ = 11
\let\@oldatcatcode = \thinlines
\def\smash@@{\relax 
  \ifmmode\def\next{\mathpalette\mathsm@sh}\else\let\next\makesm@sh
  \fi\next}
\def\makesm@sh#1{\setbox\z@\hbox{#1}\finsm@sh}
\def\mathsm@sh#1#2{\setbox\z@\hbox{$\m@th#1{#2}$}\finsm@sh}
\def\finsm@sh{\ht\z@\z@ \dp\z@\z@ \box\z@}
\edef\@oldandcatcode{\the\catcode`& }%
\catcode`& = 11
\def\&whilenoop#1{}%
\def\&whiledim#1\do #2{\ifdim #1\relax#2\&iwhiledim{#1\relax#2}\fi}%
\def\&iwhiledim#1{\ifdim #1\let\&nextwhile=\&iwhiledim 
        \else\let\&nextwhile=\&whilenoop\fi\&nextwhile{#1}}%
\newif\if&negarg
\newdimen\&wholewidth
\newdimen\&halfwidth
\font\tenln=line10
\def\thinlines{\let\&linefnt\tenln \let\&circlefnt\tencirc
  \&wholewidth\fontdimen8\tenln \&halfwidth .5\&wholewidth}%
\def\thicklines{\let\&linefnt\tenlnw \let\&circlefnt\tencircw
  \&wholewidth\fontdimen8\tenlnw \&halfwidth .5\&wholewidth}%
\def\drawline(#1,#2)#3{\&xarg #1\relax \&yarg #2\relax \&linelen=#3\relax
  \ifnum\&xarg =0 \&vline \else \ifnum\&yarg =0 \&hline \else \&sline\fi\fi}%
\def\&sline{\leavevmode
  \ifnum\&xarg< 0 \&negargtrue \&xarg -\&xarg \&yyarg -\&yarg
  \else \&negargfalse \&yyarg \&yarg \fi
  \ifnum \&yyarg >0 \&tempcnta\&yyarg \else \&tempcnta -\&yyarg \fi
  \ifnum\&tempcnta>6 \&badlinearg \&yyarg0 \fi
  \ifnum\&xarg>6 \&badlinearg \&xarg1 \fi
  \setbox\&linechar\hbox{\&linefnt\&getlinechar(\&xarg,\&yyarg)}%
  \ifnum \&yyarg >0 \let\&upordown\raise \&clnht\z@
  \else\let\&upordown\lower \&clnht \ht\&linechar\fi
  \&clnwd=\wd\&linechar
  \&whiledim \&clnwd <\&linelen \do {%
    \&upordown\&clnht\copy\&linechar
    \advance\&clnht \ht\&linechar
    \advance\&clnwd \wd\&linechar
  }%
  \advance\&clnht -\ht\&linechar
  \advance\&clnwd -\wd\&linechar
  \&tempdima\&linelen\advance\&tempdima -\&clnwd
  \&tempdimb\&tempdima\advance\&tempdimb -\wd\&linechar
  \hskip\&tempdimb \multiply\&tempdima \@m
  \&tempcnta \&tempdima \&tempdima \wd\&linechar \divide\&tempcnta \&tempdima
  \&tempdima \ht\&linechar \multiply\&tempdima \&tempcnta
  \divide\&tempdima \@m
  \advance\&clnht \&tempdima
  \ifdim \&linelen <\wd\&linechar \hskip \wd\&linechar
  \else\&upordown\&clnht\copy\&linechar\fi}%
\def\&hline{\vrule height \&halfwidth depth \&halfwidth width \&linelen}%
\def\&getlinechar(#1,#2){\&tempcnta#1\relax\multiply\&tempcnta 8
  \advance\&tempcnta -9 \ifnum #2>0 \advance\&tempcnta #2\relax\else
  \advance\&tempcnta -#2\relax\advance\&tempcnta 64 \fi
  \char\&tempcnta}%
\def\drawvector(#1,#2)#3{\&xarg #1\relax \&yarg #2\relax
  \&tempcnta \ifnum\&xarg<0 -\&xarg\else\&xarg\fi
  \ifnum\&tempcnta<5\relax \&linelen=#3\relax
    \ifnum\&xarg =0 \&vvector \else \ifnum\&yarg =0 \&hvector
    \else \&svector\fi\fi\else\&badlinearg\fi}%
\def\&hvector{\ifnum\&xarg<0 \rlap{\&linefnt\&getlarrow(1,0)}\fi \&hline
  \ifnum\&xarg>0 \llap{\&linefnt\&getrarrow(1,0)}\fi}%
\def\&vvector{\ifnum \&yarg <0 \&downvector \else \&upvector \fi}%
\def\&svector{\&sline
  \&tempcnta\&yarg \ifnum\&tempcnta <0 \&tempcnta=-\&tempcnta\fi
  \ifnum\&tempcnta <5 
    \if&negarg\ifnum\&yarg>0                   
      \llap{\lower\ht\&linechar\hbox to\&linelen{\&linefnt
        \&getlarrow(\&xarg,\&yyarg)\hss}}\else 
      \llap{\hbox to\&linelen{\&linefnt\&getlarrow(\&xarg,\&yyarg)\hss}}\fi
    \else\ifnum\&yarg>0                        
      \&tempdima\&linelen \multiply\&tempdima\&yarg
      \divide\&tempdima\&xarg \advance\&tempdima-\ht\&linechar
      \raise\&tempdima\llap{\&linefnt\&getrarrow(\&xarg,\&yyarg)}\else
      \&tempdima\&linelen \multiply\&tempdima-\&yarg 
      \divide\&tempdima\&xarg
      \lower\&tempdima\llap{\&linefnt\&getrarrow(\&xarg,\&yyarg)}\fi\fi
  \else\&badlinearg\fi}%
\def\&getlarrow(#1,#2){\ifnum #2 =\z@ \&tempcnta='33\else
\&tempcnta=#1\relax\multiply\&tempcnta \sixt@@n \advance\&tempcnta
-9 \&tempcntb=#2\relax\multiply\&tempcntb \tw@
\ifnum \&tempcntb >0 \advance\&tempcnta \&tempcntb\relax
\else\advance\&tempcnta -\&tempcntb\advance\&tempcnta 64
\fi\fi\char\&tempcnta}%
\def\&getrarrow(#1,#2){\&tempcntb=#2\relax
\ifnum\&tempcntb < 0 \&tempcntb=-\&tempcntb\relax\fi
\ifcase \&tempcntb\relax \&tempcnta='55 \or 
\ifnum #1<3 \&tempcnta=#1\relax\multiply\&tempcnta
24 \advance\&tempcnta -6 \else \ifnum #1=3 \&tempcnta=49
\else\&tempcnta=58 \fi\fi\or 
\ifnum #1<3 \&tempcnta=#1\relax\multiply\&tempcnta
24 \advance\&tempcnta -3 \else \&tempcnta=51\fi\or 
\&tempcnta=#1\relax\multiply\&tempcnta
\sixt@@n \advance\&tempcnta -\tw@ \else
\&tempcnta=#1\relax\multiply\&tempcnta
\sixt@@n \advance\&tempcnta 7 \fi\ifnum #2<0 \advance\&tempcnta 64 \fi
\char\&tempcnta}%
\def\&vline{\ifnum \&yarg <0 \&downline \else \&upline\fi}%
\def\&upline{\hbox to \z@{\hskip -\&halfwidth \vrule width \&wholewidth
   height \&linelen depth \z@\hss}}%
\def\&downline{\hbox to \z@{\hskip -\&halfwidth \vrule width \&wholewidth
   height \z@ depth \&linelen \hss}}%
\def\&upvector{\&upline\setbox\&tempboxa\hbox{\&linefnt\char'66}\raise 
     \&linelen \hbox to\z@{\lower \ht\&tempboxa\box\&tempboxa\hss}}%
\def\&downvector{\&downline\lower \&linelen
      \hbox to \z@{\&linefnt\char'77\hss}}%
\def\&badlinearg{\errmessage{Bad \string\arrow\space argument.}}%
\thinlines
\countdef\&xarg     0
\countdef\&yarg     2
\countdef\&yyarg    4
\countdef\&tempcnta 6
\countdef\&tempcntb 8
\dimendef\&linelen  0
\dimendef\&clnwd    2
\dimendef\&clnht    4
\dimendef\&tempdima 6
\dimendef\&tempdimb 8
\chardef\@arrbox    0
\chardef\&linechar  2
\chardef\&tempboxa  2           
\let\lft^%
\let\rt_
\newif\if@pslope 
\def\@findslope(#1,#2){\ifnum#1>0
  \ifnum#2>0 \@pslopetrue \else\@pslopefalse\fi \else
  \ifnum#2>0 \@pslopefalse \else\@pslopetrue\fi\fi}%
\def\generalsmap(#1,#2){\getm@rphposn(#1,#2)\plnmorph\futurelet\next\addm@rph}%
\def\sline(#1,#2){\setbox\@arrbox=\hbox{\drawline(#1,#2){\sarrowlength}}%
  \@findslope(#1,#2)\d@@blearrfalse\generalsmap(#1,#2)}%
\def\arrow(#1,#2){\setbox\@arrbox=\hbox{\drawvector(#1,#2){\sarrowlength}}%
  \@findslope(#1,#2)\d@@blearrfalse\generalsmap(#1,#2)}%
\newif\ifd@@blearr
\def\bisline(#1,#2){\@findslope(#1,#2)%
  \if@pslope \let\@upordown\raise \else \let\@upordown\lower\fi
  \getch@nnel(#1,#2)\setbox\@arrbox=\hbox{\@upordown\@vchannel
    \rlap{\drawline(#1,#2){\sarrowlength}}%
      \hskip\@hchannel\hbox{\drawline(#1,#2){\sarrowlength}}}%
  \d@@blearrtrue\generalsmap(#1,#2)}%
\def\biarrow(#1,#2){\@findslope(#1,#2)%
  \if@pslope \let\@upordown\raise \else \let\@upordown\lower\fi
  \getch@nnel(#1,#2)\setbox\@arrbox=\hbox{\@upordown\@vchannel
    \rlap{\drawvector(#1,#2){\sarrowlength}}%
      \hskip\@hchannel\hbox{\drawvector(#1,#2){\sarrowlength}}}%
  \d@@blearrtrue\generalsmap(#1,#2)}%
\def\adjarrow(#1,#2){\@findslope(#1,#2)%
  \if@pslope \let\@upordown\raise \else \let\@upordown\lower\fi
  \getch@nnel(#1,#2)\setbox\@arrbox=\hbox{\@upordown\@vchannel
    \rlap{\drawvector(#1,#2){\sarrowlength}}%
      \hskip\@hchannel\hbox{\drawvector(-#1,-#2){\sarrowlength}}}%
  \d@@blearrtrue\generalsmap(#1,#2)}%
\newif\ifrtm@rph
\def\@shiftmorph#1{\hbox{\setbox0=\hbox{$\scriptstyle#1$}%
  \setbox1=\hbox{\hskip\@hm@rphshift\raise\@vm@rphshift\copy0}%
  \wd1=\wd0 \ht1=\ht0 \dp1=\dp0 \box1}}%
\def\@hm@rphshift{\ifrtm@rph
  \ifdim\hmorphposnrt=\z@\hmorphposn\else\hmorphposnrt\fi \else
  \ifdim\hmorphposnlft=\z@\hmorphposn\else\hmorphposnlft\fi \fi}%
\def\@vm@rphshift{\ifrtm@rph
  \ifdim\vmorphposnrt=\z@\vmorphposn\else\vmorphposnrt\fi \else
  \ifdim\vmorphposnlft=\z@\vmorphposn\else\vmorphposnlft\fi \fi}%
\def\addm@rph{\ifx\next\lft\let\temp=\lftmorph\else
  \ifx\next\rt\let\temp=\rtmorph\else\let\temp\relax\fi\fi \temp}%
\def\plnmorph{\dimen1\wd\@arrbox \ifdim\dimen1<\z@ \dimen1-\dimen1\fi
  \vcenter{\box\@arrbox}}%
\def\lftmorph\lft#1{\rtm@rphfalse \setbox0=\@shiftmorph{#1}%
  \if@pslope \let\@upordown\raise \else \let\@upordown\lower\fi
  \llap{\@upordown\@vmorphdflt\hbox to\dimen1{\hss 
    \llap{\box0}\hss}\hskip\@hmorphdflt}\futurelet\next\addm@rph}%
\def\rtmorph\rt#1{\rtm@rphtrue \setbox0=\@shiftmorph{#1}%
  \if@pslope \let\@upordown\lower \else \let\@upordown\raise\fi
  \llap{\@upordown\@vmorphdflt\hbox to\dimen1{\hss
    \rlap{\box0}\hss}\hskip-\@hmorphdflt}\futurelet\next\addm@rph}%
\def\getm@rphposn(#1,#2){\ifd@@blearr \dimen@\morphdist \advance\dimen@ by
  .5\channelwidth \@getshift(#1,#2){\@hmorphdflt}{\@vmorphdflt}{\dimen@}\else
  \@getshift(#1,#2){\@hmorphdflt}{\@vmorphdflt}{\morphdist}\fi}%
\def\getch@nnel(#1,#2){\ifdim\hchannel=\z@ \ifdim\vchannel=\z@
    \@getshift(#1,#2){\@hchannel}{\@vchannel}{\channelwidth}%
    \else \@hchannel\hchannel \@vchannel\vchannel \fi
  \else \@hchannel\hchannel \@vchannel\vchannel \fi}%
\def\@getshift(#1,#2)#3#4#5{\dimen@ #5\relax
  \&xarg #1\relax \&yarg #2\relax
  \ifnum\&xarg<0 \&xarg -\&xarg \fi
  \ifnum\&yarg<0 \&yarg -\&yarg \fi
  \ifnum\&xarg<\&yarg \&negargtrue \&yyarg\&xarg \&xarg\&yarg \&yarg\&yyarg\fi
  \ifcase\&xarg \or  
    \ifcase\&yarg    
      \dimen@i \z@ \dimen@ii \dimen@ \or 
      \dimen@i .7071\dimen@ \dimen@ii .7071\dimen@ \fi \or
    \ifcase\&yarg    
      \or 
      \dimen@i .4472\dimen@ \dimen@ii .8944\dimen@ \fi \or
    \ifcase\&yarg    
      \or 
      \dimen@i .3162\dimen@ \dimen@ii .9486\dimen@ \or
      \dimen@i .5547\dimen@ \dimen@ii .8321\dimen@ \fi \or
    \ifcase\&yarg    
      \or 
      \dimen@i .2425\dimen@ \dimen@ii .9701\dimen@ \or\or
      \dimen@i .6\dimen@ \dimen@ii .8\dimen@ \fi \or
    \ifcase\&yarg    
      \or 
      \dimen@i .1961\dimen@ \dimen@ii .9801\dimen@ \or
      \dimen@i .3714\dimen@ \dimen@ii .9284\dimen@ \or
      \dimen@i .5144\dimen@ \dimen@ii .8575\dimen@ \or
      \dimen@i .6247\dimen@ \dimen@ii .7801\dimen@ \fi \or
    \ifcase\&yarg    
      \or 
      \dimen@i .1645\dimen@ \dimen@ii .9864\dimen@ \or\or\or\or
      \dimen@i .6402\dimen@ \dimen@ii .7682\dimen@ \fi \fi
  \if&negarg \&tempdima\dimen@i \dimen@i\dimen@ii \dimen@ii\&tempdima\fi
  #3\dimen@i\relax #4\dimen@ii\relax }%
\catcode`\&=4  
}%
\catcode`& = 4
\toks2 = {%
\catcode`\&=4  
\def\generalhmap{\futurelet\next\@generalhmap}%
\def\@generalhmap{\ifx\next^ \let\temp\generalhm@rph\else
  \ifx\next_ \let\temp\generalhm@rph\else \let\temp\m@kehmap\fi\fi \temp}%
\def\generalhm@rph#1#2{\ifx#1^
    \toks@=\expandafter{\the\toks@#1{\rtm@rphtrue\@shiftmorph{#2}}}\else
    \toks@=\expandafter{\the\toks@#1{\rtm@rphfalse\@shiftmorph{#2}}}\fi
  \generalhmap}%
\def\m@kehmap{\mathrel{\smash@@{\the\toks@}}}%
\def\mapright{\toks@={\mathop{\vcenter{\smash@@{\drawrightarrow}}}\limits}%
  \generalhmap}%
\def\mapleft{\toks@={\mathop{\vcenter{\smash@@{\drawleftarrow}}}\limits}%
  \generalhmap}%
\def\bimapright{\toks@={\mathop{\vcenter{\smash@@{\drawbirightarrow}}}\limits}%
  \generalhmap}%
\def\bimapleft{\toks@={\mathop{\vcenter{\smash@@{\drawbileftarrow}}}\limits}%
  \generalhmap}%
\def\adjmapright{\toks@={\mathop{\vcenter{\smash@@{\drawadjrightarrow}}}\limits}%
  \generalhmap}%
\def\adjmapleft{\toks@={\mathop{\vcenter{\smash@@{\drawadjleftarrow}}}\limits}%
  \generalhmap}%
\def\hline{\toks@={\mathop{\vcenter{\smash@@{\drawhline}}}\limits}%
  \generalhmap}%
\def\bihline{\toks@={\mathop{\vcenter{\smash@@{\drawbihline}}}\limits}%
  \generalhmap}%
\def\drawrightarrow{\hbox{\drawvector(1,0){\harrowlength}}}%
\def\drawleftarrow{\hbox{\drawvector(-1,0){\harrowlength}}}%
\def\drawbirightarrow{\hbox{\raise.5\channelwidth
  \hbox{\drawvector(1,0){\harrowlength}}\lower.5\channelwidth
  \llap{\drawvector(1,0){\harrowlength}}}}%
\def\drawbileftarrow{\hbox{\raise.5\channelwidth
  \hbox{\drawvector(-1,0){\harrowlength}}\lower.5\channelwidth
  \llap{\drawvector(-1,0){\harrowlength}}}}%
\def\drawadjrightarrow{\hbox{\raise.5\channelwidth
  \hbox{\drawvector(-1,0){\harrowlength}}\lower.5\channelwidth
  \llap{\drawvector(1,0){\harrowlength}}}}%
\def\drawadjleftarrow{\hbox{\raise.5\channelwidth
  \hbox{\drawvector(1,0){\harrowlength}}\lower.5\channelwidth
  \llap{\drawvector(-1,0){\harrowlength}}}}%
\def\drawhline{\hbox{\drawline(1,0){\harrowlength}}}%
\def\drawbihline{\hbox{\raise.5\channelwidth
  \hbox{\drawline(1,0){\harrowlength}}\lower.5\channelwidth
  \llap{\drawline(1,0){\harrowlength}}}}%
\def\generalvmap{\futurelet\next\@generalvmap}%
\def\@generalvmap{\ifx\next\lft \let\temp\generalvm@rph\else
  \ifx\next\rt \let\temp\generalvm@rph\else \let\temp\m@kevmap\fi\fi \temp}%
\toksdef\toks@@=1
\def\generalvm@rph#1#2{\ifx#1\rt 
    \toks@=\expandafter{\the\toks@
      \rlap{$\vcenter{\rtm@rphtrue\@shiftmorph{#2}}$}}\else 
    \toks@@={\llap{$\vcenter{\rtm@rphfalse\@shiftmorph{#2}}$}}%
    \toks@=\expandafter\expandafter\expandafter{\expandafter\the\expandafter
      \toks@@ \the\toks@}\fi \generalvmap}%
\def\m@kevmap{\the\toks@}%
\def\mapdown{\toks@={\vcenter{\drawdownarrow}}\generalvmap}%
\def\mapup{\toks@={\vcenter{\drawuparrow}}\generalvmap}%
\def\bimapdown{\toks@={\vcenter{\drawbidownarrow}}\generalvmap}%
\def\bimapup{\toks@={\vcenter{\drawbiuparrow}}\generalvmap}%
\def\adjmapdown{\toks@={\vcenter{\drawadjdownarrow}}\generalvmap}%
\def\adjmapup{\toks@={\vcenter{\drawadjuparrow}}\generalvmap}%
\def\vline{\toks@={\vcenter{\drawvline}}\generalvmap}%
\def\bivline{\toks@={\vcenter{\drawbivline}}\generalvmap}%
\def\drawdownarrow{\hbox to5pt{\hss\drawvector(0,-1){\varrowlength}\hss}}%
\def\drawuparrow{\hbox to5pt{\hss\drawvector(0,1){\varrowlength}\hss}}%
\def\drawbidownarrow{\hbox to5pt{\hss\hbox{\drawvector(0,-1){\varrowlength}}%
  \hskip\channelwidth\hbox{\drawvector(0,-1){\varrowlength}}\hss}}%
\def\drawbiuparrow{\hbox to5pt{\hss\hbox{\drawvector(0,1){\varrowlength}}%
  \hskip\channelwidth\hbox{\drawvector(0,1){\varrowlength}}\hss}}%
\def\drawadjdownarrow{\hbox to5pt{\hss\hbox{\drawvector(0,-1){\varrowlength}}%
  \hskip\channelwidth\lower\varrowlength
  \hbox{\drawvector(0,1){\varrowlength}}\hss}}%
\def\drawadjuparrow{\hbox to5pt{\hss\hbox{\drawvector(0,1){\varrowlength}}%
  \hskip\channelwidth\raise\varrowlength
  \hbox{\drawvector(0,-1){\varrowlength}}\hss}}%
\def\drawvline{\hbox to5pt{\hss\drawline(0,1){\varrowlength}\hss}}%
\def\drawbivline{\hbox to5pt{\hss\hbox{\drawline(0,1){\varrowlength}}%
  \hskip\channelwidth\hbox{\drawline(0,1){\varrowlength}}\hss}}%
\def\commdiag#1{\null\,
  \vcenter{\commdiagbaselines
  \m@th\ialign{\hfil$##$\hfil&&\hfil$\mkern4mu ##$\hfil\crcr
      \mathstrut\crcr\noalign{\kern-\baselineskip}
      #1\crcr\mathstrut\crcr\noalign{\kern-\baselineskip}}}\,}%
\def\commdiagbaselines{\baselineskip15pt \lineskip3pt \lineskiplimit3pt }%
\def\gridcommdiag#1{\null\,
  \vcenter{\offinterlineskip
  \m@th\ialign{&\vbox to\vgrid{\vss
    \hbox to\hgrid{\hss\smash@@{$##$}\hss}}\crcr
      \mathstrut\crcr\noalign{\kern-\vgrid}
      #1\crcr\mathstrut\crcr\noalign{\kern-.5\vgrid}}}\,}%
\newdimen\harrowlength \harrowlength=60pt
\newdimen\varrowlength \varrowlength=.618\harrowlength
\newdimen\sarrowlength \sarrowlength=\harrowlength
\newdimen\hmorphposn \hmorphposn=\z@
\newdimen\vmorphposn \vmorphposn=\z@
\newdimen\morphdist  \morphdist=4pt
\dimendef\@hmorphdflt 0       
\dimendef\@vmorphdflt 2       
\newdimen\hmorphposnrt  \hmorphposnrt=\z@
\newdimen\hmorphposnlft \hmorphposnlft=\z@
\newdimen\vmorphposnrt  \vmorphposnrt=\z@
\newdimen\vmorphposnlft \vmorphposnlft=\z@

\newdimen\hgrid \hgrid=15pt
\newdimen\vgrid \vgrid=15pt
\newdimen\hchannel  \hchannel=0pt
\newdimen\vchannel  \vchannel=0pt
\newdimen\channelwidth \channelwidth=3pt
\dimendef\@hchannel 0         
\dimendef\@vchannel 2         
\catcode`& = \@oldandcatcode
\catcode`@ = \@oldatcatcode
}%
\let\newif = \@plainnewif
\let\newdimen = \@plainnewdimen
\ifx\noarrow\@undefined \the\toks0 \the\toks2 \fi
\catcode`& = \@eplainoldandcode
\def\environment#1{%
   \ifx\@groupname\@undefined\else
      \errhelp = \@unnamedendgrouphelp
      \errmessage{`\@groupname' was not closed by \string\endenvironment}%
   \fi
   \edef\@groupname{#1}%
   \begingroup
      \let\@groupname = \@undefined
}%
\def\endenvironment#1{%
   \endgroup
   \edef\@thearg{#1}%
   \ifx\@groupname\@thearg
   \else
      \ifx\@groupname\@undefined
         \errhelp = \@isolatedendenvironmenthelp
         \errmessage{Isolated \string\endenvironment\space for `#1'}%
      \else
         \errhelp = \@mismatchedenvironmenthelp
         \errmessage{Environment `#1' ended, but `\@groupname' started}%
         \endgroup 
      \fi
   \fi
   \let\@groupname = \@undefined
}%
\newhelp\@unnamedendgrouphelp{Most likely, you just forgot an^^J%
   \string\endenvironment.  Maybe you should try inserting another^^J%
   \string\endgroup to recover.}%
\newhelp\@isolatedendenvironmenthelp{You ended an environment X, but^^J%
   no \string\environment{X} to start it is anywhere in sight.^^J%
   You might also be at an \string\endenvironment\space that would match^^J%
   a \string\begingroup, i.e., you forgot an \string\endgroup.}%
\newhelp\@mismatchedenvironmenthelp{You started an environment named X, but^^J%
   you ended one named Y.  Maybe you made a typo in one^^J%
   or the other of the names?}%
\newif\ifenvironment
\def\checkenv{\ifenvironment \errhelp = \@interwovenenvhelp
   \errmessage{Interwoven environments}%
   \egroup \fi
}%
\newhelp\@interwovenenvhelp{Perhaps you forgot to end the previous^^J%
   environment? I'm finishing off the current group,^^J%
   hoping that will fix it.}%
\newtoks\previouseverydisplay
\let\@leftleftfill\relax 
\newdimen\leftdisplayindent \leftdisplayindent=\parindent
\newif\if@leftdisplays
\def\leftdisplays{%
  \if@leftdisplays\else
    \previouseverydisplay = \everydisplay
    \everydisplay = {\the\previouseverydisplay \leftdisplaysetup}%
    \let\@save@maybedisableeqno = \@maybedisableeqno
    \let\@saveeqno = \eqno
    \let\@saveleqno = \leqno
    \let\@saveeqalignno = \eqalignno
    \let\@saveleqalignno = \leqalignno
    \let\@maybedisableeqno = \relax
    \def\eqno{\hfill\textstyle\enspace}%
    \def\leqno{%
      \hfill
      \hbox to0pt\bgroup
        \kern-\displaywidth
        \kern-\leftdisplayindent    
        $\aftergroup\@leftleqnoend  
    }%
    \@redefinealignmentdisplays
    \@leftdisplaystrue
  \fi
}%
\newbox\@lignbox
\newdimen\disprevdepth
\def\centereddisplays{%
  \if@leftdisplays
    \everydisplay = \previouseverydisplay
    \let\@maybedisableeqno = \@save@maybedisableeqno
    \let\eqno = \@saveeqno
    \let\leqno = \@saveleqno
    \let\eqalignno = \@saveeqalignno
    \let\leqalignno = \@saveleqalignno
    \@leftdisplaysfalse
  \fi
}%
\def\leftdisplaysetup{%
   \dimen@ = \leftdisplayindent
   \advance\dimen@ by \leftskip
   \advance\displayindent by \dimen@
   \advance\displaywidth by -\dimen@
   \halign\bgroup##\cr \noalign\bgroup
      \disprevdepth = \prevdepth
      \setbox\z@ = \hbox to\displaywidth\bgroup
      $\displaystyle
      \aftergroup\@lefteqend 
}
\def\@lefteqend{
   \hfil\egroup
   \@putdisplay}
\def\@leftleqnoend{\hss \egroup $}
\def\@putdisplay{%
   \ifvoid\@lignbox 
     \moveright\displayindent\box\z@ 
   \else 
     \prevdepth = \dp\@lignbox 
     \unvbox\@lignbox
   \fi
   \egroup\egroup 
   $
}
\def\@redefinealignmentdisplays{%
  \def\displaylines##1{
    \global\setbox\@lignbox\vbox{%
      \prevdepth = \disprevdepth
      \displ@y
      \tabskip\displayindent
      \halign{\hbox to\displaywidth
        {$\@lign\displaystyle####\hfil$\hfil}\crcr
              ##1\crcr}}}%
  \def\eqalignno##1{%
    \def\eqno{&}%
    \def\leqno{&}%
    \global\setbox\@lignbox\vbox{%
      \prevdepth = \disprevdepth
      \displ@y
      \advance\displaywidth by \displayindent
      \tabskip\displayindent
      \halign to\displaywidth{%
         \hfil $\@lign\displaystyle{####}$\@leftleftfill\tabskip\z@skip
        &$\@lign\displaystyle{{}####}$\hfil\tabskip\centering
        &\llap{$\@lign####$}\tabskip\z@skip\crcr
        ##1\crcr}}}%
  \def\leqalignno##1{%
    \def\eqno{&}%
    \def\leqno{&}%
    \global\setbox\@lignbox\vbox{%
      \prevdepth = \disprevdepth
      \displ@y
      \advance\displaywidth by \displayindent
      \tabskip\displayindent
      \halign to\displaywidth{%
         \hfil $\@lign\displaystyle{####}$\@leftleftfill\tabskip\z@skip
        &$\@lign\displaystyle{{}####}$\hfil\tabskip\centering
        &\kern-\displaywidth 
         \rlap{\kern\displayindent \kern-\leftdisplayindent$\@lign####$}%
         \tabskip\displaywidth\crcr
        ##1\crcr}}}%
}%
\let\@primitivenoalign = \noalign
\newtoks\@everynoalign
\def\@lefteqalignonoalign#1{%
  \@primitivenoalign{%
    \advance\leftskip by -\parindent
    \advance\leftskip by -\leftdisplayindent
    \parskip = 0pt
    \parindent = 0pt
    \the\@everynoalign
    #1%
  }%
}%
\def\monthname{%
   \ifcase\month
      \or Jan\or Feb\or Mar\or Apr\or May\or Jun%
      \or Jul\or Aug\or Sep\or Oct\or Nov\or Dec%
   \fi
}%
\def\fullmonthname{%
   \ifcase\month
      \or January\or February\or March\or April\or May\or June%
      \or July\or August\or September\or October\or November\or December%
   \fi
}%
\def\timestring{\begingroup
   \count0 = \time
   \divide\count0 by 60
   \count2 = \count0   
   \count4 = \time
   \multiply\count0 by 60
   \advance\count4 by -\count0   
   \ifnum\count4<10
      \toks1 = {0}%
   \else
      \toks1 = {}%
   \fi
   \ifnum\count2<12
      \toks0 = {a.m.}%
   \else
      \toks0 = {p.m.}%
      \advance\count2 by -12
   \fi
   \ifnum\count2=0
      \count2 = 12
   \fi
   \number\count2:\the\toks1 \number\count4 \thinspace \the\toks0
\endgroup}%
\def\today{\the\day\ \fullmonthname\ \the\year}%
\newskip\abovelistskipamount      \abovelistskipamount = .5\baselineskip
  \newcount\abovelistpenalty      \abovelistpenalty    = 10000
  \def\abovelistskip{\vpenalty\abovelistpenalty \vskip\abovelistskipamount}%
\newskip\interitemskipamount      \interitemskipamount = 0pt
  \newcount\belowlistpenalty      \belowlistpenalty    = -50
\newskip\belowlistskipamount      \belowlistskipamount = .5\baselineskip
  \newcount\interitempenalty      \interitempenalty    = 0
  \def\interitemskip{\vpenalty\interitempenalty \vskip\interitemskipamount}%
\newdimen\listleftindent    \listleftindent = 0pt
\newdimen\listrightindent   \listrightindent = 0pt        
\let\listmarkerspace = \enspace
\newtoks\everylist
\newdimen\@listindent
\def\beginlist{%
  \abovelistskip
  \@listindent = \parindent
  \advance\@listindent by \listleftindent
  \advance\leftskip by \@listindent
  \advance\rightskip by \listrightindent
  \itemnumber = 1
  \the\everylist
}%
\def\li{\@getoptionalarg\@finli}%
\def\@finli{%
  \let\@lioptarg\@optionalarg
  \ifx\@lioptarg\empty \else
    \begingroup
      \@@hldestoff
      \expandafter\writeitemxref\expandafter{\@lioptarg}%
    \endgroup
  \fi
  \ifnum\itemnumber=1 \else \interitemskip \fi
  \begingroup
    \ifx\@lioptarg\empty \else
      \toks@ = \expandafter{\@lioptarg}%
      \let\li@nohldest@marker\marker
      \edef\marker{\noexpand\hldest@impl{li}{\the\toks@}\noexpand\li@nohldest@marker}%
    \fi
    \printitem
  \endgroup
  \advance\itemnumber by 1
  \advance\itemletter by 1
  \advance\itemromannumeral by 1
  \ignorespaces
}%
\def\writeitemxref#1{\definexref{#1}\marker{item}}%
\def\printitem{%
  \par
  \nobreak
  \vskip-\parskip
  \noindent
  \printmarker\marker
}%
\def\printmarker#1{\llap{\marker \enspace}}%
\newcount\numberedlistdepth
\newcount\itemnumber
\newcount\itemletter
\newcount\itemromannumeral
\def\numberedmarker{%
  \ifcase\numberedlistdepth
      (impossible)%
  \or \printitemnumber
  \or \printitemletter
  \or \printitemromannumeral
  \else *%
  \fi
}%
\def\printitemnumber{\number\itemnumber}%
\def\printitemletter{\char\the\itemletter}%
\def\printitemromannumeral{\romannumeral\itemromannumeral}%
\def\numberedprintmarker#1{\llap{#1) \listmarkerspace}}%
\def\numberedlist{\environment{@numbered-list}%
  \advance\numberedlistdepth by 1
  \itemletter = `a
  \itemromannumeral = 1
  \beginlist
  \let\marker = \numberedmarker
  \let\printmarker = \numberedprintmarker
}%

\newcount\unorderedlistdepth
\def\unorderedmarker{%
  \ifcase\unorderedlistdepth
      (impossible)%
  \or \blackbox
  \or ---%
  \else *%
  \fi
}%
\def\unorderedprintmarker#1{\llap{#1\listmarkerspace}}%
\def\unorderedlist{\environment{@unordered-list}%
  \advance\unorderedlistdepth by 1
  \beginlist
  \let\marker = \unorderedmarker
  \let\printmarker = \unorderedprintmarker
}%
\def\listing#1{%
   \par \begingroup
   \@setuplisting
   \setuplistinghook
   \input #1
   \endgroup
}%
\let\setuplistinghook = \relax
\def\linenumberedlisting{%
  \ifx\lineno\undefined \innernewcount\lineno \fi
  \lineno = 0
  \everypar = {\advance\lineno by 1 \printlistinglineno}%
}%
\def\printlistinglineno{\llap{[\the\lineno]\quad}}%
\def\nolastlinelisting{\aftergroup\@removeboxes}%
\def\@removeboxes{%
  \setbox0 = \lastbox
  \ifvoid0
    \ignorespaces 
  \else
    \expandafter\@removeboxes
  \fi
}%
{%
  \makeactive\^^L
  \let^^L = \relax
  \gdef\@setuplisting{%
     \uncatcodespecials
     \obeywhitespace
     \makeactive\`
     \makeactive\^^I
     \makeactive\^^L
     \def^^L{\vfill\break}%
     \parskip = 0pt
     \listingfont
  }%
}%
\def\listingfont{\tt}%
{%
   \makeactive\`
   \gdef`{\relax\lq}
}%
{%
   \makeactive\^^I
   \gdef^^I{\hskip8\fontdimen2}%
}%
\def\verbatimescapechar#1{%
  \gdef\@makeverbatimescapechar{%
    \@makeverbatimdoubleescape #1%
    \catcode`#1 = 0
  }%
}%
\def\@makeverbatimdoubleescape#1{%
  \catcode`#1 = \other
  \begingroup
    \lccode`\* = `#1%
    \lowercase{\endgroup \ece\def*{*}}%
}%
\verbatimescapechar\|  
\def\verbatim{\begingroup
  \uncatcodespecials
  \makeactive\` 
  \@makeverbatimescapechar
  \tt\obeywhitespace}

\def\definecontentsfile#1{%
  \ece\innernewwrite{#1file}%
  \ece\innernewif{if@#1fileopened}%
  \ece\let{#1filebasename} = \jobname
  \ece\def{open#1file}{\opencontentsfile{#1}}%
  \ece\def{write#1entry}{\writecontentsentry{#1}}%
  \ece\def{writenumbered#1entry}{\writenumberedcontentsentry{#1}}%
  \ece\def{writenumbered#1line}{\writenumberedcontentsline{#1}}%
  \ece\innernewif{ifrewrite#1file} \csname rewrite#1filetrue\endcsname
  \ece\def{read#1file}{\readcontentsfile{#1}}%
}%
\definecontentsfile{toc}%
\def\opencontentsfile#1{%
  \csname if@#1fileopened\endcsname \else
     \ece{\immediate\openout}{#1file} = \csname #1filebasename\endcsname.#1
     \ece\global{@#1fileopenedtrue}%
  \fi
}%
\def\writecontentsentry#1#2#3{\writenumberedcontentsentry{#1}{#2}{#3}{}}%
\def\writenumberedcontentsentry#1#2#3#4{%
  \csname ifrewrite#1file\endcsname
    \writenumberedcontents@cmdname{#1}{#2}%
    \def\temp{#3}
    \toks2 = \expandafter{#4}%
    \edef\cs{\the\toks2}%
    \edef\@wr{%
      \write\csname #1file\endcsname{%
        \the\toks0 
        {\sanitize\temp}
        \ifx\empty\cs\else {\sanitize\cs}\fi 
        {\noexpand\folio}
      }%
    }%
    \@wr
  \fi
  \ignorespaces
}%
\def\writenumberedcontentsline#1#2#3#4{%
  \csname ifrewrite#1file\endcsname
    \writenumberedcontents@cmdname{#1}{#2}%
    \def\temp{#4}
    \toks2 = \expandafter{#3}%
    \edef\cs{\the\toks2}%
    \edef\@wr{%
      \write\csname #1file\endcsname{%
        \the\toks0 
        \ifx\empty\cs\else {\sanitize\cs}\fi 
        {\sanitize\temp}
        {\noexpand\folio}
      }%
    }%
    \@wr
  \fi
  \ignorespaces
}%
\def\writenumberedcontents@cmdname#1#2{%
  \csname open#1file\endcsname
  \edef\temp{#2}
  \expandafter\if\expandafter\isinteger\expandafter{\temp}%
    \toks0 = {\expandafter\noexpand \csname #1entry\endcsname}%
    \edef\temp{\the\toks0{\temp}}%
    \toks0 = \expandafter{\temp}%
  \else
    \toks0 = {\expandafter\noexpand \csname #1#2entry\endcsname}%
  \fi
}%
\def\readcontentsfile#1{%
   \edef\temp{%
     \noexpand\testfileexistence[\csname #1filebasename\endcsname]{#1}%
   }\temp
   \if@fileexists
      \input \csname #1filebasename\endcsname.#1\relax
   \fi
}%
\def\tocsectionentry#1#2{\line{\quad\sl #1 \dotfill\ \rm #2}}%
\def\tocsubsectionentry#1#2{\line{\qquad\rm #1 \dotfill\ #2}}%
\let\ifxrefwarning = \iftrue
\def\xrefwarningtrue{\@citewarningtrue \let\ifxrefwarning = \iftrue}%
\def\xrefwarningfalse{\@citewarningfalse \let\ifxrefwarning = \iffalse}%
\begingroup
  \catcode`\_ = 8
  \gdef\xrlabel#1{#1_x}%
\endgroup
\def\xrdef#1{%
  \begingroup
    \hldest@impl{xrdef}{#1}%
    \begingroup
      \@@hldestoff
      \definexref{#1}{\noexpand\folio}{page}%
    \endgroup
  \endgroup
  \ignorespaces
}%
\def\definexref#1#2#3{%
  \hldest@impl{definexref}{#1}%
  \edef\temp{#1}%
  \readauxfile
  \edef\@wr{\noexpand\writeaux{\string\@definelabel{\temp}{#2}{#3}}}%
  \@wr
  \ignorespaces
}%
\def\@definelabel#1{
  \begingroup 
    \expandafter\ifx\csname\xrlabel{#1}\endcsname \relax
      \expandafter\@definelabel@nocheck
    \else
      \expandafter\@definelabel@warn
    \fi
    {#1}%
}%
\def\@definelabel@nocheck#1#2#3{%
    \expandafter\gdef\csname\xrlabel{#1}\endcsname{#2}%
    \setpropertyglobal{\xrlabel{#1}}{class}{#3}%
  \endgroup 
}%
\def\@definelabel@warn#1#2#3{%
  \message{^^J\linenumber Label `#1' multiply defined,
           value `#2' of class `#3' overwriting value
           `\csname\xrlabel{#1}\endcsname' of class
           `\getproperty{\xrlabel{#1}}{class}'.}%
  \@definelabel@nocheck{#1}{#2}{#3}%
}%
\def\reftie{\penalty\@M \ }
\let\refspace\ 
\def\xrefn{\@getoptionalarg\@finxrefn}%
\def\@finxrefn#1{%
  \hlstart@impl{ref}{#1}%
  \ifx\@optionalarg\empty \else
    \let\@xrefnoptarg\@optionalarg
    \readauxfile
    {\@@hloff\@xrefnoptarg}\reftie
  \fi
  \plain@xrefn{#1}%
  \hlend@impl{ref}%
}%
\def\plain@xrefn#1{%
  \readauxfile
  \expandafter \ifx\csname\xrlabel{#1}\endcsname\relax
    \if@citewarning
       \message{\linenumber Undefined label `#1'.}%
    \fi
    \expandafter\def\csname\xrlabel{#1}\endcsname{%
      `{\tt
        \escapechar = -1
        \expandafter\string\csname#1\endcsname
      }'%
    }%
  \fi
  \csname\xrlabel{#1}\endcsname 
}%
\let\refn = \xrefn
\def\xrefpageword{p.\thinspace}%
\def\xref{\@getoptionalarg\@finxref}%
\def\@finxref#1{%
  \hlstart@impl{xref}{#1}%
  \ifx\@optionalarg\empty \else
    {\@@hloff\@optionalarg}\refspace
  \fi
  \xrefpageword\plain@xrefn{#1}%
  \hlend@impl{xref}%
}%
\def\@maybewarnref{%
  \ifundefined{amsppt.sty}%
  \else
    \message{Warning: amsppt.sty and Eplain both define \string\ref. See
             the Eplain manual.}%
    \let\amsref = \ref
  \fi
  \let\ref = \eplainref
  \ref
}
\let\ref = \@maybewarnref
\def\eplainref{\@getoptionalarg\@fineplainref}%
\def\@fineplainref{\@generalref{1}{}}%
\def\refs{\let\@optionalarg\empty \@generalref{0}s}%
\def\@generalref#1#2#3{%
  \let\@generalrefoptarg\@optionalarg
  \readauxfile
  \ifcase#1 \else \hlstart@impl{ref}{#3}\fi
  \edef\@generalref@class{\getproperty{\xrlabel{#3}}{class}}%
  \expandafter\ifx\csname \@generalref@class word\endcsname\relax
    \ifx\@generalrefoptarg\empty \else {\@@hloff\@generalrefoptarg\reftie}\fi
  \else
    \begingroup
      \@@hloff
      \ifx\@generalrefoptarg\empty \else \@generalrefoptarg \refspace \fi
      \csname \@generalref@class word\endcsname
      #2\reftie
    \endgroup
  \fi
  \ifcase#1 \hlstart@impl{ref}{#3}\fi
  \plain@xrefn{#3}%
  \hlend@impl{ref}%
}%
\newcount\eqnumber
\newcount\subeqnumber
\def\eqdefn{\@getoptionalarg\@fineqdefn}%
\def\@fineqdefn#1{%
  \ifx\@optionalarg\empty
    \global\advance\eqnumber by 1
    \def\temp{\eqconstruct{\number\eqnumber}}%
  \else
    \def\temp{\@optionalarg}%
  \fi
  \global\subeqnumber = 0
  \gdef\@currenteqlabel{#1}%
  \toks0 = \expandafter{\@currenteqlabel}%
  \begingroup
    \def\eqrefn{\noexpand\plain@xrefn}%
    \def\xrefn{\noexpand\plain@xrefn}%
    \edef\temp{\noexpand\@eqdefn{\the\toks0}{\temp}}%
    \temp
  \endgroup
}%
\def\eqsubdefn#1{%
  \global\advance\subeqnumber by 1
  \toks0 = {#1}%
  \toks2 = \expandafter{\@currenteqlabel}%
  \begingroup
    \def\eqrefn{\noexpand\plain@xrefn}%
    \def\xrefn{\noexpand\plain@xrefn}%
    \def\eqsubreftext{\noexpand\eqsubreftext}%
    \edef\temp{%
      \noexpand\@eqdefn
        {\the\toks0}%
        {\eqsubreftext{\eqrefn{\the\toks2}}{\the\subeqnumber}}%
    }%
    \temp           
  \endgroup
}%
\newcount\phantomeqnumber
\def\phantomeqlabel{PHEQ\the\phantomeqnumber}%
\def\@eqdefn#1#2{%
  \ifempty{#1}%
    \global\advance\phantomeqnumber by 1
    \edef\hl@eqlabel{\phantomeqlabel}%
    \readauxfile
  \else
    \def\hl@eqlabel{#1}%
    {\@@hldestoff \definexref{#1}{#2}{eq}}%
  \fi
  \hldest@impl{eq}{\hl@eqlabel}%
  \begingroup 
    \@definelabel@nocheck{#1}{#2}{eq}%
}%
\def\eqdef{\@getoptionalarg\@fineqdef}%
\def\@fineqdef{%
  \toks0 = \expandafter{\@optionalarg}%
  \edef\temp{\noexpand\@eqdef{\noexpand\eqdefn[\the\toks0]}}%
  \temp
}%
\def\eqsubdef{\@eqdef\eqsubdefn}%
\def\@eqdef#1#2{%
  \@maybedisableeqno
  \eqnum #1{#2}
        \let\@optionalarg\empty 
        {\@@hloff\@fineqref{#2}}
  \@mayberestoreeqno
  \ignorespaces
}%
\let\@mayberestoreeqno = \relax
\def\@maybedisableeqno{%
  \ifinner
    \global\let\eqno = \relax
    \global\let\leqno = \relax
    \global\let\@mayberestoreeqno = \@restoreeqno
  \fi
}%
\let\@primitiveeqno = \eqno
\let\@primitiveleqno = \leqno
\def\@restoreeqno{%
  \global\let\eqno = \@primitiveeqno
  \global\let\leqno = \@primitiveleqno
  \global\let\@mayberestoreeqno = \empty
}%
\def\righteqnumbers{%
  \def\eqnum{\eqno}%
  \def\eqalignnum{\eqalignno}%
}%
\def\lefteqnumbers{%
  \def\eqnum{\leqno}%
  \def\eqalignnum{\leqalignno}%
}%
\righteqnumbers
\def\eqrefn{\@getoptionalarg\@fineqrefn}%
\def\@fineqrefn#1{%
  \eqref@start{#1}%
  \plain@xrefn{#1}%
  \hlend@impl{eq}%
}%
\def\eqref{\@getoptionalarg\@fineqref}%
\def\@fineqref#1{%
  \eqref@start{#1}%
  \eqprint{\plain@xrefn{#1}}%
  \hlend@impl{eq}%
}%
\def\eqref@start#1{%
  \let\@eqrefoptarg\@optionalarg
  \ifempty{#1}%
    \hlstart@impl{eq}{\phantomeqlabel}%
  \else
    \hlstart@impl{eq}{#1}%
  \fi
  \ifx\@eqrefoptarg\empty \else
    {\@@hloff\@eqrefoptarg\reftie}%
  \fi
}%
\let\eqconstruct = \identity
\def\eqprint#1{(#1)}%
\def\eqsubreftext#1#2{#1.#2}%
\let\extraidxcmdsuffixes = \empty
\outer\def\defineindex#1{%
  \def\@idxprefix{#1}%
  \expandafter\innernewif\csname if\@idxprefix dx\endcsname
  \csname \@idxprefix dxtrue\endcsname
  \for\@idxcmd:=,marked,submarked,name%
                \extraidxcmdsuffixes\do
  {%
    \@defineindexcmd\@idxcmd
  }%
  \ece\innernewwrite{@#1indexfile}%
  \ece\innernewif{if@#1indexfileopened}%
}%
\newif\ifsilentindexentry
\def\@defineindexcmd#1{%
  \@defineoneindexcmd{s}{#1}\silentindexentrytrue
  \@defineoneindexcmd{}{#1}\silentindexentryfalse
}%
\def\@defineoneindexcmd#1#2#3{%
  \toks@ = {#3}%
  \edef\temp{%
    \def
      \expandonce\csname#1\@idxprefix dx#2\endcsname 
      {\def\noexpand\@idxprefix{\@idxprefix}
       \expandonce\csname @@#1idx#2\endcsname
      }%
    \def
      \expandonce\csname @@#1idx#2\endcsname{
        \the\toks@
        \noexpand\@idxgetrange\expandonce\csname @#1idx#2\endcsname
      }%
  }%
  \temp
}%
\let\indexfilebasename = \jobname
\def\@idxwrite#1#2{%
  \csname if\@idxprefix dx\endcsname
    \@openidxfile
    \def\temp{#1}%
    \edef\@wr{%
      \expandafter\write\csname @\@idxprefix indexfile\endcsname{%
        \string\indexentry
        {\sanitize\temp}%
        {\noexpand#2}%
      }%
    }%
    \@wr
  \else
    \write-1{}%
  \fi
  \ifindexproofing
    \def\temp{#1}%
    \edef\temp{%
      \insert\@indexproof{\noexpand\indexproofterm{\sanitize\temp}}%
    }%
    \temp
    \ifhmode\allowhyphens\fi
  \fi
  \hookrun{afterindexterm}%
  \ifsilentindexentry \expandafter\ignorespaces\fi
}%
\def\@openidxfile{%
  \csname if@\@idxprefix indexfileopened\endcsname \else
    \expandafter\immediate\openout\csname @\@idxprefix indexfile\endcsname =
      \indexfilebasename.\@idxprefix dx
    \expandafter\global\csname @\@idxprefix indexfileopenedtrue\endcsname
  \fi
}%
\newif\ifindexproofing
\newinsert\@indexproof
\dimen\@indexproof = \maxdimen                  
\count\@indexproof = 0  \skip\@indexproof = 0pt 
\font\indexprooffont = cmtt8
\def\indexproofterm#1{\hbox{\strut \indexprooffont #1}}%
\let\@plainmakeheadline = \makeheadline
\def\makeheadline{%
  \expandafter\ifx\csname\idxpageanchor{\folio}\endcsname\relax \else
    {\@@hldeston \hldest@impl{idx}{\hlidxpagelabel{\folio}}}%
  \fi
  \indexproofunbox
  \@plainmakeheadline
}%
\def\indexsetmargins{%
  \ifx\undefined\outsidemargin
    \dimen@ = 1truein
    \advance\dimen@ by \hoffset
    \edef\outsidemargin{\the\dimen@}%
    \let\insidemargin = \outsidemargin
  \fi
}%
\def\indexproofunbox{%
  \ifvoid\@indexproof\else
    \indexsetmargins
    \rlap{%
      \kern\hsize
      \ifodd\pageno \kern\outsidemargin \else \kern\insidemargin \fi
      \vbox to 0pt{\unvbox\@indexproof\vss}%
    }\nointerlineskip
  \fi
}%
\def\idxrangebeginword{begin}%
\def\idxbeginrangemark{(}
\def\idxrangeendword{end}%
\def\idxendrangemark{)}%
\def\idxseecmdword{see}%
\def\idxseealsocmdword{seealso}%
\newif\if@idxsee
\newif\if@hlidxencap
\let\@idxseenterm = \relax
\def\idxpagemarkupcmdword{pagemarkup}%
\let\@idxpagemarkup = \relax
\def\@idxgetrange#1{%
  \let\@idxrangestr = \empty
  \let\@afteridxgetrange = #1%
  \begingroup
    \catcode\idxargopen=1
    \@getoptionalarg\@finidxgetopt
}%
\def\@finidxgetopt{%
    \global\let\@idxgetrange@arg\@optionalarg
  \endgroup
  \@hlidxencaptrue
  \for\@idxarg:=\@idxgetrange@arg\do{%
    \expandafter\@idxcheckpagemarkup\@idxarg=,%
    \ifx\@idxarg\idxrangebeginword
      \def\@idxrangestr{\idxencapoperator\idxbeginrangemark}%
    \else
      \ifx\@idxarg\idxrangeendword
        \def\@idxrangestr{\idxencapoperator\idxendrangemark}%
        \@hlidxencapfalse
      \else
        \ifx\@idxarg\idxseecmdword
          \def\@idxpagemarkup{indexsee}%
          \@idxseetrue
          \@hlidxencapfalse
        \else
          \ifx\@idxarg\idxseealsocmdword
            \def\@idxpagemarkup{indexseealso}%
            \@idxseetrue
            \@hlidxencapfalse
          \else
             \ifx\@idxpagemarkup\relax
               \errmessage{Unrecognized index option `\@idxarg'}%
             \fi
          \fi
        \fi
      \fi
    \fi
  }%
  \ifnum\hldest@place@idx < 0 \else
    \if@hlidxencap
      \ifx\@idxpagemarkup\relax
        \let\@idxpagemarkup\empty
      \fi
      \ifcase\hldest@place@idx \relax
        \edef\@idxpagemarkup{hlidxpage{\@idxpagemarkup}}%
        \definepageanchor{\noexpand\folio}%
      \else
        \global\advance\hlidxlabelnumber by 1
        \edef\@idxpagemarkup{hlidx{\hlidxlabel}{\@idxpagemarkup}}%
        \hldest@impl{idx}{\hlidxlabel}%
      \fi
    \fi
  \fi
  \@afteridxgetrange
}%
\def\@idxcheckpagemarkup#1=#2,{%
  \def\temp{#1}%
  \ifx\temp\idxpagemarkupcmdword
    \if ,#2, 
      \errmessage{Missing markup command to `pagemarkup'}%
    \else
      \def\temp##1={##1}%
      \edef\@idxpagemarkup{\temp\string#2}%
    \fi
  \fi
}%
\def\hldest@place@idx{-1}%
\begingroup
  \catcode`\_ = 8
  \gdef\idxpageanchor#1{#1_p}%
\endgroup
\def\definepageanchor#1{%
  \readauxfile
  \edef\@wr{\noexpand\writeaux{\string\@definepageanchor{#1}}}%
  \@wr
  \ignorespaces
}%
\def\@definepageanchor#1{%
  \expandafter\gdef\csname\idxpageanchor{#1}\endcsname{}%
}%
\newcount\hlidxlabelnumber
\def\hlidxlabel{IDX\the\hlidxlabelnumber}%
\def\hlidxpagelabel#1{IDXPG#1}%
\def\hlidx#1#2#3{%
  \ifempty{#2}%
    \let\@idxpageencap\relax
  \else
    \expandafter\let\expandafter\@idxpageencap\csname #2\endcsname
  \fi
  \hlstart@impl{idx}{#1}%
  \@idxpageencap{#3}%
  \hlend@impl{idx}%
}%
\def\hlidxpage#1#2{%
  \@hlidxgetpages{#2}%
  \ifempty{#1}%
    \let\@idxpageencap\relax
  \else
    \expandafter\let\expandafter\@idxpageencap\csname #1\endcsname
  \fi
  \hlstart@impl{idx}{\hlidxpagelabel{\@idxpageiref}}%
  \expandafter\@idxpageencap\expandafter{\@idxpagei}%
  \hlend@impl{idx}%
  \ifx\@idxpageii\empty \else
    \@idxpagesep
    \hlstart@impl{idx}{\hlidxpagelabel{\@idxpageiiref}}%
    \expandafter\@idxpageencap\expandafter{\@idxpageii}%
    \hlend@impl{idx}%
  \fi
}%
\def\@hlidxgetpages#1{%
  \idxparselist{#1}%
  \ifx\idxpagei\empty
    \idxparserange{#1}%
    \ifx\idxpagei\empty
      \def\@idxpageiref{#1}
    \else
      \let\@idxpageiref\idxpagei 
    \fi
    \def\@idxpagei{#1}%
    \let\@idxpageii\empty
  \else
    \let\@idxpagei\idxpagei
    \let\@idxpageii\idxpageii
    \let\@idxpageiref\idxpagei 
    \let\@idxpageiiref\idxpageii 
    \let\@idxpagesep\idxpagelistdelimiter
  \fi
}%
\def\setidxpagelistdelimiter#1{%
  \gdef\idxpagelistdelimiter{#1}%
  \gdef\@removepagelistdelimiter##1#1{##1}%
  \gdef\@idxparselist##1#1##2\@@{%
    \ifempty{##2}%
      \let\idxpagei\empty
    \else
      \def\idxpagei{##1}%
      \edef\idxpageii{\@removepagelistdelimiter##2}%
    \fi
  }%
  \gdef\idxparselist##1{\@idxparselist##1#1\@@}%
}%
\def\setidxpagerangedelimiter#1{%
  \gdef\idxpagerangedelimiter{#1}%
  \gdef\@removepagerangedelimiter##1#1{##1}%
  \gdef\@idxparserange##1#1##2\@@{%
    \ifempty{##2}%
      \let\idxpagei\empty
    \else
      \def\idxpagei{##1}%
      \edef\idxpageii{\@removepagerangedelimiter##2}%
    \fi
  }%
  \gdef\idxparserange##1{\@idxparserange##1#1\@@}%
}%
\setidxpagelistdelimiter{, }%
\setidxpagerangedelimiter{--}%
\def\idxsubentryseparator{!}%
\def\idxencapoperator{|}%
\def\idxmaxpagenum{99999}%
\newtoks\@idxmaintoks
\newtoks\@idxsubtoks
\def\@idxtokscollect{%
  \edef\temp{\the\@idxsubtoks}%
  \edef\@indexentry{%
    \the\@idxmaintoks
    \ifx\temp\empty\else \idxsubentryseparator\the\@idxsubtoks \fi
    \@idxrangestr
  }%
  \if@idxsee
    \@idxseefalse 
    \edef\temp{\noexpand\idx@getverbatimarg
      {\noexpand\@finidxtokscollect{\idxmaxpagenum}}}%
  \else
    \def\temp{\@finfinidxtokscollect\folio}%
  \fi
  \temp
}%
\def\@finidxtokscollect#1#2{%
  \def\@idxseenterm{#2}%
  \@finfinidxtokscollect{#1}%
}%
\def\@finfinidxtokscollect#1{%
  \ifx\@idxpagemarkup\relax \else
    \toks@ = \expandafter{\@indexentry}%
    \edef\@indexentry{%
      \the\toks@
      \ifx\@idxrangestr\empty \idxencapoperator \fi
      \@idxpagemarkup
    }%
    \let\@idxpagemarkup = \relax
  \fi
  \ifx\@idxseenterm\relax \else
    \toks@ = \expandafter{\@indexentry}%
    \edef\@indexentry{\the\toks@{\sanitize\@idxseenterm}}%
    \let\@idxseenterm = \relax
  \fi
  \expandafter\@idxwrite\expandafter{\@indexentry}{#1}%
}%
\def\@idxcollect#1#2{%
  \@idxmaintoks = {#1}%
  \@idxsubtoks = {#2}%
  \@idxtokscollect
}%
\def\idxargopen{`\{}%
\def\idxargclose{`\}}%
\def\idx@getverbatimarg#1{%
  \def\idx@getverbatimarg@cmd{#1}
  \begingroup
    \uncatcodespecials
    \catcode\idxargopen=1
    \catcode\idxargclose=2
    \catcode`\ =10   
    \catcode`\^^I=10 
    \catcode`\^^M=5  
    \idx@fingetverbatimarg
}%
\def\idx@fingetverbatimarg#1{\endgroup\idx@getverbatimarg@cmd{#1}}%
\def\idx@getverboptarg#1{%
  \def\idx@optionaltemp{#1}
  \let\idx@optionalnext = \relax
  \begingroup
    \if@idxsee \catcode\idxargopen=1 \fi
    \@futurenonspacelet\idx@optionalnext\idx@bracketcheck
}%
\def\idx@bracketcheck{%
  \ifx [\idx@optionalnext
    \endgroup 
    \expandafter\idx@finbracketcheck
  \else
    \endgroup 
    \let\@optionalarg = \empty
    \expandafter\idx@optionaltemp
  \fi
}%
\def\idx@finbracketcheck{%
  \begingroup
    \uncatcodespecials
    \catcode`\ =10   
    \catcode`\^^I=10 
    \catcode`\^^M=5  
    \idx@@getoptionalarg
}%
\def\idx@@getoptionalarg[#1]{%
  \endgroup
  \def\@optionalarg{#1}%
  \idx@optionaltemp
}%
{\catcode`\&=0
\gdef\idx@scanterm#1{%
  \edef\idx@scanterm@nl@catcode{\the\catcode`\^^M\relax}%
  \catcode`\^^M=5
  \scantokens{#1&endinput}%
  \catcode`\^^M=\idx@scanterm@nl@catcode
}}%
\def\@idx{\idx@getverbatimarg\@finidx}%
\def\@finidx#1{%
  \idx@scanterm{#1}
  \@idxcollect{#1}{}%
}%
\def\@sidx{\idx@getverbatimarg\@finsidx}%
\def\@finsidx#1{\@idxmaintoks = {#1}\idx@getverboptarg\@finfinsidx}%
\def\@finfinsidx{%
  \@idxsubtoks = \expandafter{\@optionalarg}%
  \@idxtokscollect
}%
\def\idxsortkeysep{@}
\def\@idxconstructmarked#1#2#3{%
  \toks@ = {#2}
  \toks2 = {#3}
  \edef\temp{\the\toks2 \idxsortkeysep \the\toks@{\the\toks2}}%
  #1 = \expandafter{\temp}%
}%
\def\@idxmarked#1{\idx@getverbatimarg{\@finidxmarked{#1}}}%
\def\@finidxmarked#1#2{%
  \idx@scanterm{#1{#2}}
  \@idxconstructmarked\@idxmaintoks{#1}{#2}%
  \@idxsubtoks = {}%
  \@idxtokscollect
}%
\def\@sidxmarked#1{\idx@getverbatimarg{\@finsidxmarked{#1}}}%
\def\@finsidxmarked#1#2{%
  \@idxconstructmarked\toks@{#1}{#2}%
  \edef\temp{{\the\toks@}}%
  \expandafter\@finsidx\temp
}%
\def\@idxsubmarked{%
  \let\sidxsubmarked@print\idxsubmarked@print
  \idx@getverbatimarg\@finsidxsubmarked
}%
\def\idxsubmarked@print{\expandafter\idx@scanterm\expandafter}%
\def\@sidxsubmarked{%
  \let\sidxsubmarked@print\gobble
  \idx@getverbatimarg\@finsidxsubmarked
}%
\def\@finsidxsubmarked#1{\@idxmaintoks = {#1}\@@finsidxsubmarked}
\def\@@finsidxsubmarked#1{\idx@getverbatimarg{\@@@finsidxsubmarked{#1}}}
\def\@@@finsidxsubmarked#1#2{
  \sidxsubmarked@print 
    {\the\@idxmaintoks\space #1{#2}}
  \@idxconstructmarked\@idxsubtoks{#1}{#2}%
  \@idxtokscollect
}%
\def\idxnameseparator{, }
\def\@idxcollectname#1#2{%
  \def\temp{#1}%
  \ifx\temp\empty
    \toks@ = {}%
  \else
    \toks@ = \expandafter{\idxnameseparator #1}%
  \fi
  \toks2 = {#2}%
  \edef\temp{\the\toks2 \the\toks@}%
}%
\def\@idxname{\idx@getverbatimarg\@finidxname}%
\def\@finidxname#1{\idx@getverbatimarg{\@finfinidxname{#1}}}%
\def\@finfinidxname#1#2{%
  \idx@scanterm{#1 #2}
  \@idxcollectname{#1}{#2}%
  \expandafter\@idxcollect\expandafter{\temp}{}%
}%
\def\@sidxname{\idx@getverbatimarg\@finsidxname}%
\def\@finsidxname#1{\idx@getverbatimarg{\@finfinsidxname{#1}}}%
\def\@finfinsidxname#1#2{%
  \@idxcollectname{#1}{#2}%
  \expandafter\@finsidx\expandafter{\temp}%
}%
\let\indexfonts = \relax
\def\readindexfile#1{%
  \edef\@idxprefix{#1}%
  \testfileexistence[\indexfilebasename]{\@idxprefix nd}%
  \iffileexists \begingroup
    \ifx\begin\undefined
      \def\begin##1{\@beginindex}%
      \let\end = \@gobble
    \fi
    \input \indexfilebasename.\@idxprefix nd
    \singlecolumn
  \endgroup
  \else
    \message{No index file \indexfilebasename.\@idxprefix nd.}%
  \fi
}%
\def\@beginindex{%
  \let\item = \@indexitem
  \let\subitem = \@indexsubitem
  \let\subsubitem = \@indexsubsubitem
  \indexfonts
  \doublecolumns
  \parindent = 0pt
  \hookrun{beginindex}%
}%

\newskip\aboveindexitemskipamount  \aboveindexitemskipamount = 0pt plus2pt
\def\aboveindexitemskip{\vskip\aboveindexitemskipamount}%
\def\@indexitem{\begingroup
  \@indexitemsetup
  \leftskip = 0pt
  \aboveindexitemskip
  \penalty-100 
  \def\par{\endgraf\endgroup\nobreak}%
}%
\def\@indexsubitem{%
  \@indexitemsetup
  \leftskip = 1em
}%
\def\@indexsubsubitem{%
  \@indexitemsetup
  \leftskip = 2em
}%
\def\@indexitemsetup{%
  \par
  \hangindent = 1em
  \raggedright
  \hyphenpenalty = 10000
  \hookrun{indexitem}%
}%
\def\seevariant{\it}%
\def\indexseeword{see}%
\def\indexsee{\idx@getverbatimarg\@finindexsee}%
\def\@finindexsee#1#2{{\seevariant \indexseeword\/ }\idx@scanterm{#1}}%
\def\indexseealsowords{see also}%
\def\indexseealso{\idx@getverbatimarg\@finindexseealso}%
\def\@finindexseealso#1#2{{\seevariant \indexseealsowords\/ }\idx@scanterm{#1}}%
\defineindex{i}%
\begingroup
  \catcode `\^^M = \active %
  \gdef\flushleft{%
    \def\@endjustifycmd{\@endflushleft}%
    \def\@eoljustifyaction{\null\hfil\break}%
    \let\@firstlinejustifyaction = \relax
    \@startjustify %
  }%
  \gdef\flushright{%
    \def\@endjustifycmd{\@endflushright}%
    \def\@eoljustifyaction{\break\null\hfil}%
    \def\@firstlinejustifyaction{\hfil\null}%
    \@startjustify %
  }%
  \gdef\center{%
    \def\@endjustifycmd{\@endcenter}%
    \def\@eoljustifyaction{\hfil\break\null\hfil}%
    \def\@firstlinejustifyaction{\hfil\null}%
    \@startjustify %
  }%
  \gdef\@startjustify{%
    \parskip = 0pt
    \catcode`\^^M = \active %
    \def^^M{\futurelet\next\@finjustifyreturn}%
    \def\@eateol##1^^M{%
      \def\temp{##1}%
      \@firstlinejustifyaction %
      \ifx\temp\empty\else \temp^^M\fi %
    }%
    \expandafter\aftergroup\@endjustifycmd %
    \checkenv \environmenttrue %
    \par\noindent %
    \@eateol %
  }%
  \gdef\@finjustifyreturn{%
    \@eoljustifyaction %
    \ifx\next^^M%
      \def\par{\endgraf\vskip\blanklineskipamount \global\let\par = \endgraf}%
      \@endjustifycmd %
      \noindent %
      \@firstlinejustifyaction %
    \fi %
  }%
\endgroup
\def\@endflushleft{\unpenalty{\parfillskip = 0pt plus1fil\par}\ignorespaces}%
\def\@endflushright{
   \unskip \setbox0=\lastbox \unpenalty
   {\parfillskip = 0pt \par}\ignorespaces
}%
\def\@endcenter{
   \unskip \setbox0=\lastbox \unpenalty
   {\parfillskip = 0pt plus1fil \par}\ignorespaces
}%
\newcount\abovecolumnspenalty   \abovecolumnspenalty = 10000
\newcount\@linestogo         
\newcount\@linestogoincolumn 
\newcount\@columndepth       
\newdimen\@columnwidth       
\newtoks\crtok  \crtok = {\cr}%
\newcount\currentcolumn
\def\makecolumns#1/#2 {\par \begingroup 
   \@columndepth = #1
   \advance\@columndepth by -1
   \divide \@columndepth by #2
   \advance\@columndepth by 1
   \@linestogoincolumn = \@columndepth
   \@linestogo = #1
   \currentcolumn = 1
   \def\@endcolumnactions{%
      \ifnum \@linestogo<2 
         \the\crtok \egroup \endgroup \par 
      \else
         \global\advance\@linestogo by -1
         \ifnum\@linestogoincolumn<2
            \global\advance\currentcolumn by 1
            \global\@linestogoincolumn = \@columndepth
            \the\crtok
         \else
            &\global\advance\@linestogoincolumn by -1
         \fi
      \fi
   }%
   \makeactive\^^M
   \letreturn \@endcolumnactions
   \@columnwidth = \hsize
     \advance\@columnwidth by -\parindent
     \divide\@columnwidth by #2
   \penalty\abovecolumnspenalty
   \noindent 
   \valign\bgroup
     &\hbox to \@columnwidth{\strut \hsize = \@columnwidth ##\hfil}\cr
}%
\newcount\footnotenumber
\newcount\hlfootlabelnumber
\newdimen\footnotemarkseparation \footnotemarkseparation = .5em
\newskip\interfootnoteskip \interfootnoteskip = 0pt
\newtoks\everyfootnote
\newdimen\footnoterulewidth \footnoterulewidth = 2in
\newdimen\footnoteruleheight \footnoteruleheight = 0.4pt
\newdimen\belowfootnoterulespace \belowfootnoterulespace = 2.6pt
\let\@plainfootnote = \footnote
\let\@plainvfootnote = \vfootnote
\def\hlfootlabel{FOOT\the\hlfootlabelnumber}%
\def\hlfootbacklabel{FOOTB\the\hlfootlabelnumber}%
\def\@eplainfootnote#1{\let\@sf\empty 
  \ifhmode\edef\@sf{\spacefactor\the\spacefactor}\/\fi
  \global\advance\hlfootlabelnumber by 1
  \hlstart@impl{foot}{\hlfootlabel}%
  \hldest@impl{footback}{\hlfootbacklabel}%
  #1%
  \hlend@impl{foot}%
  \@sf\vfootnote{#1}%
}%
\let\footnote\@eplainfootnote
\def\vfootnote#1{\insert\footins\bgroup
  \interlinepenalty\interfootnotelinepenalty
  \splittopskip\ht\strutbox 
  \advance\splittopskip by \interfootnoteskip
  \splitmaxdepth\dp\strutbox
  \floatingpenalty\@MM
  \leftskip\z@skip \rightskip\z@skip \spaceskip\z@skip \xspaceskip\z@skip
  \everypar = {}%
  \parskip = 0pt 
  \ifnum\@numcolumns > 1 \hsize = \@normalhsize \fi
  \the\everyfootnote
  \vskip\interfootnoteskip
  \indent\llap{%
    \hlstart@impl{footback}{\hlfootbacklabel}%
    \hldest@impl{foot}{\hlfootlabel}%
    #1%
    \hlend@impl{footback}%
    \kern\footnotemarkseparation
  }%
  \footstrut\futurelet\next\fo@t
}%
\def\footnoterule{\dimen@ = \footnoteruleheight
  \advance\dimen@ by \belowfootnoterulespace
  \kern-\dimen@
  \hrule width\footnoterulewidth height\footnoteruleheight depth0pt
  \kern\belowfootnoterulespace
  \vskip-\interfootnoteskip
}%
\def\numberedfootnote{%
  \global\advance\footnotenumber by 1
  \@eplainfootnote{{\number\footnotenumber}}
}%
\newdimen\paperheight 
\ifnum\mag=1000
  \paperheight = 11in 
\else
  \paperheight = 11truein 
\fi
\def\topmargin{\afterassignment\@finishtopmargin \dimen@}%
\def\@finishtopmargin{%
  \dimen2 = \voffset    
  \voffset = \dimen@ \advance\voffset by -1truein
  \advance\dimen2 by -\voffset  
  \advance\vsize by \dimen2 
}%
\def\advancetopmargin{%
  \dimen@ = 0pt \afterassignment\@finishadvancetopmargin \advance\dimen@
}%
\def\@finishadvancetopmargin{%
  \advance\voffset by \dimen@
  \advance\vsize by -\dimen@
}%
\def\bottommargin{\afterassignment\@finishbottommargin \dimen@}%
\def\@finishbottommargin{%
  \@computebottommargin   
  \advance\dimen2 by -\dimen@ 
  \advance\vsize by \dimen2 
}%
\def\advancebottommargin{%
  \dimen@ = 0pt \afterassignment\@finishadvancebottommargin \advance\dimen@
}%
\def\@finishadvancebottommargin{%
  \advance\vsize by -\dimen@
}%
\def\@computebottommargin{%
  \dimen2 = \paperheight  
  \advance\dimen2 by -\vsize  
  \advance\dimen2 by -\voffset  
  \advance\dimen2 by -1truein 
}%
\newdimen\paperwidth
\ifnum\mag=1000
  \paperwidth = 8.5in 
\else
  \paperwidth = 8.5truein 
\fi
\def\leftmargin{\afterassignment\@finishleftmargin \dimen@}%
\def\@finishleftmargin{%
  \dimen2 = \hoffset    
  \hoffset = \dimen@ \advance\hoffset by -1truein
  \advance\dimen2 by -\hoffset  
  \advance\hsize by \dimen2 
}%
\def\advanceleftmargin{%
  \dimen@ = 0pt \afterassignment\@finishadvanceleftmargin \advance\dimen@
}%
\def\@finishadvanceleftmargin{%
  \advance\hoffset by \dimen@
  \advance\hsize by -\dimen@
}%
\def\rightmargin{\afterassignment\@finishrightmargin \dimen@}%
\def\@finishrightmargin{%
  \@computerightmargin    
  \advance\dimen2 by -\dimen@ 
  \advance\hsize by \dimen2 
}%
\def\advancerightmargin{%
  \dimen@ = 0pt \afterassignment\@finishadvancerightmargin \advance\dimen@
}%
\def\@finishadvancerightmargin{%
  \advance\hsize by -\dimen@
}%
\def\@computerightmargin{%
  \dimen2 = \paperwidth   
  \advance\dimen2 by -\hsize  
  \advance\dimen2 by -\hoffset  
  \advance\dimen2 by -1truein 
}%
\let\@plainm@g = \m@g
\def\m@g{\@plainm@g \paperwidth = 8.5 true in \paperheight = 11 true in}%
\newskip\abovecolumnskip \abovecolumnskip = \bigskipamount
\newskip\belowcolumnskip \belowcolumnskip = \bigskipamount
\newdimen\gutter \gutter = 2pc
\newbox\@partialpage
\newdimen\@normalhsize
\newdimen\@normalvsize  
\newtoks\previousoutput
\def\quadcolumns{\@columns4}%
\def\triplecolumns{\@columns3}%
\def\doublecolumns{\@columns2}%
\def\begincolumns#1{\ifcase#1\relax \or \singlecolumn \or \@columns2 \or
                            \@columns3 \or \@columns4 \else \relax \fi}%
\let\@ndcolumns = \relax
\chardef\@numcolumns = 1
\mathchardef\@ejectpartialpenalty = 10141
\chardef\@col@minlines = 3
\chardef\@col@extralines = 3
\newdimen\@col@extraheight
\def\@columns#1{%
  \@ndcolumns
  \global\let\@ndcolumns = \@endcolumns
  \global\chardef\@numcolumns = #1
  \global\previousoutput = \expandafter{\the\output}%
  \global\output = {%
    \ifnum\outputpenalty = -\@ejectpartialpenalty
      \dimen@ = \vsize
      \advance\dimen@ by \@col@minlines\baselineskip
      \global\setbox\@partialpage =
        \vbox  \ifdim \pagetotal > \vsize  to \dimen@  \fi  {%
    \unvbox255 \unskip
  }%
    \else
      \the\previousoutput
    \fi
  }%
  \vskip \abovecolumnskip
  \vskip \@col@minlines\baselineskip
  \penalty -\@ejectpartialpenalty
  \global\output = {\@columnoutput}%
  \global\@normalhsize = \hsize
  \global\@normalvsize = \vsize
  \count@ = \@numcolumns
  \advance\count@ by -1
  \global\advance\hsize by -\count@\gutter
  \global\divide\hsize by \@numcolumns
  \advance\vsize by -\ht\@partialpage
  \advance\vsize by -\ht\footins
  \ifvoid\footins\else \advance\vsize by -\skip\footins \fi
  \multiply\count\footins by \@numcolumns
  \advance\vsize by -\ht\topins
  \ifvoid\topins\else \advance\vsize by -\skip\topins \fi
  \multiply\count\topins by \@numcolumns
  \global\vsize = \@numcolumns\vsize
  \@col@extraheight=\@col@extralines\baselineskip
  \multiply\@col@extraheight by \@numcolumns
  \global\advance\vsize by \@col@extraheight
}%
\def\gutterbox{\vbox to \dimen0{\vfil\hbox{\hfil}\vfil}}%
\def\@columnsplit{%
  \splittopskip = \topskip
  \splitmaxdepth = \baselineskip
  \begingroup
    \vbadness = 10000
    \global\setbox1 = \vsplit255 to \dimen@  \global\wd1 = \hsize
    \global\setbox3 = \vsplit255 to \dimen@  \global\wd3 = \hsize
    \ifnum\@numcolumns > 2
      \global\setbox5 = \vsplit255 to \dimen@ \global\wd5 = \hsize
    \fi
    \ifnum\@numcolumns > 3
      \global\setbox7 = \vsplit255 to \dimen@ \global\wd7 = \hsize
    \fi
  \endgroup
  \setbox0 = \box255
  \global\setbox255 = \vbox{%
    \unvbox\@partialpage
    \ifcase\@numcolumns \relax\or\relax
      \or \hbox to \@normalhsize{\box1\hfil\gutterbox\hfil\box3}%
      \or \hbox to \@normalhsize{\box1\hfil\gutterbox\hfil\box3%
                                      \hfil\gutterbox\hfil\box5}%
      \or \hbox to \@normalhsize{\box1\hfil\gutterbox\hfil\box3%
                                      \hfil\gutterbox\hfil\box5%
                                      \hfil\gutterbox\hfil\box7}%
    \fi
  }%
  \setbox\@partialpage = \box0
}%
\def\@columnoutput{%
  \dimen@ = \ht255
    \advance\dimen@ by -\@col@extraheight
    \divide\dimen@ by \@numcolumns
  \@columnsplit
  \@recoverclubpenalty 
  \hsize = \@normalhsize 
  \vsize = \@normalvsize
  \the\previousoutput
  \unvbox\@partialpage
  \penalty\outputpenalty
  \global\vsize = \@numcolumns\@normalvsize
  \global\advance\vsize by \@col@extraheight
}%
\def\singlecolumn{%
  \@ndcolumns
  \chardef\@numcolumns = 1
  \vskip\belowcolumnskip
  \nointerlineskip
}%
\newbox\@singlecolumnbox 
\newdimen\column@pagegoal
\newdimen\column@vsize
\def\@endcolumns{%
  \global\let\@ndcolumns = \relax
  \par 
  \column@pagegoal = \pagegoal
  \advance\column@pagegoal by-\@col@extraheight
  \ifdim \pagetotal > \column@pagegoal
    \column@vsize = \column@pagegoal
  \else
    \column@vsize = \pagetotal
  \fi
  \global\output = {\global\setbox1 = \box255}%
  \pagegoal = \pagetotal
  \break                     
  \setbox2 = \box1           
  \global\output = \expandafter{\the\previousoutput}%
  \setbox\@singlecolumnbox = \box\@partialpage
  \@balancecolumns
}%
\def\@balancecolumns{%
  \global\setbox255 = \copy2  
  \dimen@ = \column@vsize
    \divide\dimen@ by \@numcolumns
  \@columnsplit
  \ifvoid\@partialpage
    \global\vsize = \@normalvsize
    \global\hsize = \@normalhsize
    \dump@balanced@columns
    \let\next\relax
  \else
    \advance \column@vsize by \@numcolumns pt
    \test@spill@columns
    \ifspill@columns
      \begingroup
        \vsize = \@normalvsize
        \hsize = \@normalhsize
        \dump@balanced@columns
        \break
        \@recoverclubpenalty
      \endgroup
      \unvbox\@partialpage
      \let\next\@endcolumns
    \else
      \setbox0=\box\@partialpage 
      \let\next\@balancecolumns
    \fi
  \fi
  \next
}%
\def\dump@balanced@columns{%
  \ifvoid\topins\else\topinsert\unvbox\topins\endinsert\fi
  \unvbox\@singlecolumnbox
  \nointerlineskip
  \box255
}%
\newif\ifspill@columns
\def\test@spill@columns{%
  \spill@columnsfalse
  \ifdim \column@vsize > \column@pagegoal
    \ifvoid\footins
      \ifvoid\topins
        \spill@columnstrue
      \fi
    \fi
  \fi
}%
\def\@saveclubpenalty{
  \edef\@recoverclubpenalty{%
     \global\clubpenalty=\the\clubpenalty\relax%
     \global\let\noexpand\@recoverclubpenalty\relax
  }
}%
\let\@recoverclubpenalty\relax
\newdimen\temp@dimen
\def\columnfill{%
  \par
  \dimen@=\pagetotal   
  \temp@dimen = \vsize 
  \divide\temp@dimen by \@numcolumns 
  \loop
    \ifdim \dimen@ > \temp@dimen
      \advance \dimen@ by -\temp@dimen
      \advance \dimen@ by \topskip 
  \repeat
  \advance \temp@dimen by -\dimen@
  \advance \temp@dimen by -\prevdepth
  \@saveclubpenalty 
  \clubpenalty=10000\relax
  \hrule height\temp@dimen width0pt depth0pt\relax
  \nointerlineskip
  \par
  \nointerlineskip
  \penalty0\vfil 
  \relax
}%
\def\@hldest{%
  \def\hl@prefix{hldest}%
  \let\after@hl@getparam\hldest@aftergetparam
  \begingroup
    \hl@getparam
}%
\def\hldest@aftergetparam{%
  \ifvmode
    \hldest@driver
  \else
    \allowhyphens
    \smash{\ifx\hldest@opt@raise\empty \else \raise\hldest@opt@raise\fi
             \hbox{\hldest@driver}}%
    \allowhyphens
  \fi
  \endgroup
}%
\def\@hlstart{%
  \leavevmode
  \def\hl@prefix{hl}%
  \let\after@hl@getparam\hlstart@aftergetparam
  \begingroup
    \hl@getparam
}%
\def\hlstart@aftergetparam{%
  \ifx\color\undefined \else
    \ifx\hl@opt@color\empty \else
      \ifx\hl@opt@colormodel\empty
        \edef\temp{\noexpand\color{\hl@opt@color}}%
      \else
        \edef\temp{\noexpand\color[\hl@opt@colormodel]{\hl@opt@color}}%
      \fi
      \temp
    \fi
  \fi
  \hl@driver
}%
\def\hl@getparam#1#2{
  \edef\@hltype{#1}%
  \ifx\@hltype\empty
    \expandafter\let\expandafter\@hltype
      \csname \hl@prefix @type\endcsname
  \fi
  \expandafter\ifx\csname \hl@prefix @typeh@\@hltype\endcsname \relax
    \errmessage{Invalid hyperlink type `\@hltype'}%
  \fi
  \For\hl@arg:=#2\do{%
    \ifx\hl@arg\empty \else
      \expandafter\hl@set@opt\hl@arg=,%
    \fi
  }%
  \bgroup
    \uncatcodespecials
    \catcode`\{=1 \catcode`\}=2
    \@hl@getparam
}%
\def\@hl@getparam#1{%
  \egroup
  \edef\@hllabel{#1}%
  \after@hl@getparam
  \ignorespaces
}%
\def\hl@set@opt#1=#2,{%
  \expandafter\ifx\csname \hl@prefix @opt@#1\endcsname \relax
    \errmessage{Invalid hyperlink option `#1'}%
  \fi
  \if ,#2, 
    \errmessage{Missing value for option `#1'}%
  \else
    \def\temp##1={##1}%
    \expandafter\edef\csname \hl@prefix @opt@#1\endcsname{\temp#2}%
  \fi
}%
\def\hldest@impl#1{%
  \expandafter\ifcase\csname hldest@on@#1\endcsname
    \relax\expandafter\gobble
  \else
    \toks@=\expandafter{\csname hldest@type@#1\endcsname}%
    \toks@ii=\expandafter{\csname hldest@opts@#1\endcsname}%
    \edef\temp{\noexpand\hldest{\the\toks@}{\the\toks@ii}}%
    \expandafter\temp
  \fi
}%
\def\hlstart@impl#1{%
  \expandafter\ifcase\csname hl@on@#1\endcsname
    \leavevmode\expandafter\gobble
  \else
    \toks@=\expandafter{\csname hl@type@#1\endcsname}%
    \toks@ii=\expandafter{\csname hl@opts@#1\endcsname}%
    \edef\temp{\noexpand\hlstart{\the\toks@}{\the\toks@ii}}%
    \expandafter\temp
  \fi
}%
\def\hlend@impl#1{%
  \expandafter\ifcase\csname hl@on@#1\endcsname
  \else
    \hlend
  \fi
}%
\def\hl@asterisk@word{*}%
\def\hl@opts@word{opts}%
\newif\if@params@override
\def\hldest@groups{definexref,xrdef,li,eq,bib,foot,footback,idx}%
\def\hl@groups{ref,xref,eq,cite,foot,footback,idx,url,hrefint,hrefext}%
\def\hldesttype{%
  \def\hl@prefix{hldest}%
  \def\hl@param{type}%
  \let\hl@all@groups\hldest@groups
  \futurelet\hl@excl\hl@param@read@excl
}%
\def\hldestopts{%
  \def\hl@prefix{hldest}%
  \def\hl@param{opts}%
  \let\hl@all@groups\hldest@groups
  \futurelet\hl@excl\hl@param@read@excl
}%
\def\hltype{%
  \def\hl@prefix{hl}%
  \def\hl@param{type}%
  \let\hl@all@groups\hl@groups
  \futurelet\hl@excl\hl@param@read@excl
}%
\def\hlopts{%
  \def\hl@prefix{hl}%
  \def\hl@param{opts}%
  \let\hl@all@groups\hl@groups
  \futurelet\hl@excl\hl@param@read@excl
}%
\def\hl@param@read@excl{%
  \ifx!\hl@excl
    \let\next\hl@param@read@opt@arg
    \@params@overridetrue
  \else
    \def\next{\hl@param@read@opt@arg{!}}%
    \@params@overridefalse
  \fi
  \next
}%
\def\hl@param@read@opt@arg#1{
  \@getoptionalarg\hl@setparam
}%
\def\@hl@setparam#1{%
  \ifx\@optionalarg\empty
    \hl@setparam@default{#1}
  \else
    \let\hl@do@all@groups\gobble
    \For\hl@group:=\@optionalarg\do{%
      \ifx\hl@group\hl@asterisk@word
        \def\hl@do@all@groups{\let\@optionalarg\hl@all@groups \hl@setparam}%
      \else
        \hl@setparam@group{#1}%
      \fi
    }%
    \hl@do@all@groups{#1}%
  \fi
}%
\def\hl@setparam@group#1{%
  \ifx\hl@group\empty
    \hl@setparam@default{#1}%
  \else
    \expandafter\ifx\csname\hl@prefix @\hl@param @\hl@group\endcsname\relax
      \errmessage{Hyperlink group `\hl@prefix:\hl@param:\hl@group' is not defined}%
    \fi
    \ifx\hl@param\hl@opts@word
      \if@params@override
        \expandafter\let\csname\hl@prefix @\hl@param @\hl@group\endcsname\empty
      \fi
      \hl@update@opts@with@list{#1}
    \else
      \ece\def{\hl@prefix @\hl@param @\hl@group}{#1}%
    \fi
  \fi
}%
\def\hl@setparam@default#1{%
  \ifx\hl@param\hl@opts@word
    \For\hl@opt:=#1\do{%
      \ifx\hl@opt\empty \else
        \expandafter\hl@set@opt\hl@opt=,%
      \fi
    }%
  \else
    \expandafter\ifx\csname\hl@prefix @\hl@param\endcsname\relax
      \message{Default hyperlink parameter `\hl@prefix:\hl@param' is not defined}%
    \fi
    \ece\def{\hl@prefix @\hl@param}{#1}%
  \fi
}%
\def\hl@update@opts@with@list#1{%
  \global\expandafter\let\expandafter\hl@update@new@list
    \csname \hl@prefix @opts@\hl@group\endcsname
  \begingroup
    \For\hl@opt:=#1\do{%
      \hl@update@opts@with@opt
    }%
  \endgroup
  \ece\let{\hl@prefix @opts@\hl@group}\hl@update@new@list
}%
\def\hl@update@opts@with@opt{%
  \global\let\hl@update@old@list\hl@update@new@list
  \global\let\hl@update@new@list\empty
  \global\let\hl@update@new@opt\hl@opt
  \expandafter\hl@parse@opt@key\hl@opt=,%
  \let\hl@update@new@key\hl@update@key
  \global\let\hl@update@comma\empty
  \begingroup
    \for\hl@opt:=\hl@update@old@list\do{%
      \ifx\hl@opt\empty \else 
        \expandafter\hl@parse@opt@key\hl@opt=,%
        \toks@=\expandafter{\hl@update@new@list}%
        \ifx\hl@update@key\hl@update@new@key
          \ifx\hl@update@new@opt\empty \else 
            \toks@ii=\expandafter{\hl@update@new@opt}%
            \xdef\hl@update@new@list{\the\toks@\hl@update@comma\the\toks@ii}%
            \global\let\hl@update@new@opt\empty
            \global\def\hl@update@comma{,}%
          \fi
        \else
          \toks@ii=\expandafter{\hl@opt}%
          \xdef\hl@update@new@list{\the\toks@\hl@update@comma\the\toks@ii}%
          \global\def\hl@update@comma{,}%
        \fi
      \fi
    }%
  \endgroup
  \ifx\hl@update@new@opt\empty \else
    \toks@=\expandafter{\hl@update@new@list}%
    \toks@ii=\expandafter{\hl@update@new@opt}%
    \xdef\hl@update@new@list{\the\toks@\hl@update@comma\the\toks@ii}%
  \fi
}%
\def\hl@parse@opt@key#1=#2,{\def\hl@update@key{#1}}%
\def\hldest@opt@raise{\normalbaselineskip}%
\def\hl@opt@colormodel{cmyk}%
\def\hl@opt@color{0.28,1,1,0.35}%
\def\hldest@on@definexref{0}%
\def\hldest@on@xrdef{0}%
\def\hldest@on@li{0}%
\def\hldest@on@eq{0}
\def\hldest@on@bib{0}
\def\hldest@on@foot{0}
\def\hldest@on@footback{0}
\def\hldest@on@idx{0}
\let\hldest@type@definexref\empty
\let\hldest@type@xrdef\empty
\let\hldest@type@li\empty
\let\hldest@type@eq\empty 
\let\hldest@type@bib\empty 
\let\hldest@type@foot\empty 
\let\hldest@type@footback\empty 
\let\hldest@type@idx\empty 
\let\hldest@opts@definexref\empty
\let\hldest@opts@xrdef\empty
\let\hldest@opts@li\empty
\def\hldest@opts@eq{raise=1.7\normalbaselineskip}
\let\hldest@opts@bib\empty 
\let\hldest@opts@foot\empty 
\let\hldest@opts@footback\empty 
\let\hldest@opts@idx\empty 
\def\hl@on@ref{0}
\def\hl@on@xref{0}%
\def\hl@on@eq{0}
\def\hl@on@cite{0}
\def\hl@on@foot{0}
\def\hl@on@footback{0}
\def\hl@on@idx{0}%
\def\hl@on@url{0}
\def\hl@on@hrefint{0}
\def\hl@on@hrefext{0}
\let\hl@type@ref\empty 
\let\hl@type@xref\empty
\let\hl@type@eq\empty 
\let\hl@type@cite\empty 
\let\hl@type@foot\empty 
\let\hl@type@footback\empty 
\let\hl@type@idx\empty
\let\hl@type@url\empty 
\let\hl@type@hrefint\empty 
\let\hl@type@hrefext\empty 
\let\hl@opts@ref\empty 
\let\hl@opts@xref\empty
\let\hl@opts@eq\empty 
\let\hl@opts@cite\empty 
\let\hl@opts@foot\empty 
\let\hl@opts@footback\empty 
\let\hl@opts@idx\empty
\let\hl@opts@url\empty 
\let\hl@opts@hrefint\empty 
\let\hl@opts@hrefext\empty 
\def\@hlon{\@hlonoff@value@stub{hl}\@@hlon1 }%
\def\@hloff{\@hlonoff@value@stub{hl}\@@hloff0 }%
\def\@hldeston{\@hlonoff@value@stub{hldest}\@@hldeston1 }%
\def\@hldestoff{\@hlonoff@value@stub{hldest}\@@hldestoff0 }%
\def\@hlonoff@value@stub#1#2#3{%
  \def\hl@prefix{#1}%
  \let\hl@on@empty#2%
  \def\hl@value{#3}%
  \expandafter\let\expandafter\hl@all@groups
    \csname \hl@prefix @groups\endcsname
  \@getoptionalarg\@finhlswitch
}%
\def\@finhlswitch{%
  \ifx\@optionalarg\empty
    \hl@on@empty
  \fi
  \let\hl@do@all@groups\relax
  \For\hl@group:=\@optionalarg\do{%
    \ifx\hl@group\hl@asterisk@word
      \let\@optionalarg\hl@all@groups
      \let\hl@do@all@groups\@finhlswitch
    \else
      \ifx\hl@group\empty
        \hl@on@empty
      \else
        \expandafter\ifx\csname\hl@prefix @on@\hl@group\endcsname \relax
          \errmessage{Hyperlink group `\hl@prefix:on:\hl@group'
                      is not defined}%
        \fi
        \ece\edef{\hl@prefix @on@\hl@group}{\hl@value}%
      \fi
    \fi
  }%
  \hl@do@all@groups
}%
\def\@@hlon{%
  \let\hlstart\@hlstart
  \let\hlend\@hlend
}%
\def\@@hloff{%
  \def\hlstart##1##2##3{\leavevmode\ignorespaces}%
  \let\hlend\relax
}%
\def\@@hldeston{%
  \let\hldest\@hldest
}%
\def\@@hldestoff{%
  \def\hldest##1##2##3{\ignorespaces}%
}%
\def\hl@idxexact@word{idxexact}%
\def\hl@idxpage@word{idxpage}%
\def\hl@idxnone@word{idxnone}%
\def\hl@raw@word{raw}%
\def\enablehyperlinks{\@getoptionalarg\@finenablehyperlinks}%
\def\@finenablehyperlinks{%
  \let\hl@selecteddriver\empty
  \def\hldest@place@idx{0}%
  \for\hl@arg:=\@optionalarg\do{%
    \ifx\hl@arg\hl@idxexact@word
      \def\hldest@place@idx{1}%
    \else
      \ifx\hl@arg\hl@idxnone@word
        \def\hldest@place@idx{-1}%
      \else
        \ifx\hl@arg\hl@idxpage@word
          \def\hldest@place@idx{0}%
        \else
          \let\hl@selecteddriver\hl@arg
        \fi
      \fi
    \fi
  }%
  \ifx\hl@selecteddriver\empty
    \ifpdf
      \def\hl@selecteddriver{pdftex}%
      \message{^^JEplain: using `pdftex' hyperlink driver.}%
    \else
      \def\hl@selecteddriver{hypertex}%
      \message{^^JEplain: using `hypertex' hyperlink driver.}%
    \fi
  \else
    \expandafter\ifx\csname hldriver@\hl@selecteddriver\endcsname \relax
      \errmessage{No hyperlink driver `\hl@selecteddriver' available}%
    \fi
  \fi
  \let\hl@setparam\@hl@setparam
  \csname hldriver@\hl@selecteddriver\endcsname
  \def\@finenablehyperlinks{\errmessage{Hyperlink driver `\hl@selecteddriver'
                                        already selected}}%
  \let\hldriver@nolinks\undefined
  \let\hldriver@hypertex\undefined
  \let\hldriver@pdftex \undefined
  \let\hldriver@dvipdfm\undefined
  \let\hloff\@hloff
  \let\hlon\@hlon
  \let\hldestoff\@hldestoff
  \let\hldeston\@hldeston
  \hlon[*,]\hloff[foot,footback]%
  \hldeston[*,]\hldestoff[foot,footback]%
}%
\def\hldriver@nolinks{%
  \def\@hldest##1##2##3{%
    \edef\temp{\write-1{hldest: ##3}}%
    \ifvmode
      \temp
    \else
      \allowhyphens
      \expandafter\smash\expandafter{\temp}%
      \allowhyphens
    \fi
    \ignorespaces
  }%
  \def\@hlstart##1##2##3{%
    \leavevmode
    \begingroup 
    \edef\temp{\write-1{hlstart: ##3}}%
    \temp
    \ignorespaces
  }%
  \def\@hlend{%
    \edef\temp{\write-1{hlend}}%
    \temp
    \endgroup 
  }%
  \let\hl@setparam\gobble
}%
{\catcode`\#=\other
\gdef\hlhash{#}}%
\def\hldriver@hypertex{%
  \def\hldest@type{xyz}%
  \let\hldest@opt@cmd \empty
  \def\hldest@driver{%
    \ifx\@hltype\hl@raw@word
      \csname \hldest@opt@cmd \endcsname
    \else
    \fi
  }%
  \let\hldest@typeh@raw \empty
  \let\hldest@typeh@xyz \empty
  \def\hl@type{name}%
  \ifx\hl@type@url\empty
    \def\hl@type@url{url}%
  \fi
  \ifx\hl@type@hrefext\empty
    \def\hl@type@hrefext{url}%
  \fi
  \let\hl@opt@cmd  \empty
  \let\hl@opt@ext  \empty
  \let\hl@opt@file \empty
  \def\hl@driver{%
    \ifx\@hltype\hl@raw@word
      \csname \hl@opt@cmd \endcsname
    \else
      \def\hlstart@preamble{html:<a href="}%
      \csname hl@typeh@\@hltype\endcsname
    \fi
  }%
  \let\hl@typeh@raw \empty
  \def\hl@typeh@name{\special{\hlstart@preamble \hlhash\@hllabel">}}%
  \def\hl@typeh@filename{%
    \special{%
      \hlstart@preamble
        file:\hl@opt@file\hl@opt@ext
        \ifempty\@hllabel \else \hlhash\@hllabel\fi
      ">%
    }%
  }%
  \def\hl@typeh@url{%
    \special{%
      \hlstart@preamble
        \@hllabel
      ">%
    }%
  }%
  \def\@hlend{\endgroup}
}%
\def\hldriver@pdftex{%
\ifpdf 
  \def\hldest@type{xyz}%
  \let\hldest@opt@width  \empty
  \let\hldest@opt@height \empty
  \let\hldest@opt@depth  \empty
  \let\hldest@opt@zoom   \empty
  \let\hldest@opt@cmd    \empty
  \def\hldest@driver{%
    \ifx\@hltype\hl@raw@word
      \csname \hldest@opt@cmd \endcsname
    \else
      \pdfdest name{\@hllabel}\@hltype
        \csname hldest@typeh@\@hltype\endcsname
    \fi
  }%
  \let\hldest@typeh@raw   \empty
  \let\hldest@typeh@fit   \empty
  \let\hldest@typeh@fith  \empty
  \let\hldest@typeh@fitv  \empty
  \let\hldest@typeh@fitb  \empty
  \let\hldest@typeh@fitbh \empty
  \let\hldest@typeh@fitbv \empty
  \def\hldest@typeh@fitr{%
    \ifx\hldest@opt@width  \empty \else width  \hldest@opt@width  \fi
    \ifx\hldest@opt@height \empty \else height \hldest@opt@height \fi
    \ifx\hldest@opt@depth  \empty \else depth  \hldest@opt@depth  \fi
  }%
  \def\hldest@typeh@xyz{%
    \ifx\hldest@opt@zoom\empty \else zoom \hldest@opt@zoom \fi
  }%
  \def\hl@type{name}%
  \ifx\hl@type@url\empty
    \def\hl@type@url{url}%
  \fi
  \ifx\hl@type@hrefext\empty
    \def\hl@type@hrefext{url}%
  \fi
  \let\hl@opt@width   \empty
  \let\hl@opt@height  \empty
  \let\hl@opt@depth   \empty
  \def\hl@opt@bstyle  {S}%
  \def\hl@opt@bwidth  {1}%
  \let\hl@opt@bcolor  \empty
  \let\hl@opt@hlight  \empty
  \let\hl@opt@bdash   \empty
  \let\hl@opt@pagefit \empty
  \let\hl@opt@cmd     \empty
  \let\hl@opt@file    \empty
  \let\hl@opt@newwin  \empty
  \def\hl@driver{%
    \ifx\@hltype\hl@raw@word
      \csname \hl@opt@cmd \endcsname
    \else
      \let\hl@BSspec\relax 
      \ifx\hl@opt@bstyle \empty
        \ifx\hl@opt@bwidth \empty
          \ifx\hl@opt@bdash \empty
            \let\hl@BSspec\empty 
          \fi
        \fi
      \fi
      \def\hlstart@preamble{%
        \pdfstartlink
          \ifx\hl@opt@width  \empty \else width  \hl@opt@width  \fi
          \ifx\hl@opt@height \empty \else height \hl@opt@height \fi
          \ifx\hl@opt@depth  \empty \else depth  \hl@opt@depth \fi
          attr{%
            \ifx\hl@opt@bcolor\empty\else /C[\hl@opt@bcolor]\fi
            \ifx\hl@opt@hlight\empty\else /H/\hl@opt@hlight\fi
            \ifx\hl@BSspec\relax
              /BS<<%
                /Type/Border%
                \ifx\hl@opt@bstyle\empty\else /S/\hl@opt@bstyle\fi
                \ifx\hl@opt@bwidth\empty\else /W \hl@opt@bwidth\fi
                \ifx\hl@opt@bdash\empty \else /D[\hl@opt@bdash]\fi
              >>%
            \fi
          }%
      }%
      \csname hl@typeh@\@hltype\endcsname
    \fi
  }%
  \let\hl@typeh@raw\empty
  \def\hl@typeh@name{\hlstart@preamble goto name{\@hllabel}}%
  \def\hl@typeh@num{\hlstart@preamble  goto num \@hllabel}%
  \def\hl@typeh@page{%
    \count@=\@hllabel
    \advance\count@ by-1
    \hlstart@preamble
    user{%
      /Subtype/Link%
      /Dest%
        [\the\count@
          \ifx\hl@opt@pagefit\empty/Fit\else\hl@opt@pagefit\fi]%
    }%
  }%
  \def\hl@typeh@filename{\hl@file{(\@hllabel)}}%
  \def\hl@typeh@filepage{%
    \count@=\@hllabel
    \advance\count@ by-1
    \hl@file{%
      [\the\count@ \ifx\hl@opt@pagefit\empty/Fit\else\hl@opt@pagefit\fi]%
    }%
  }%
  \def\hl@file##1{%
    \hlstart@preamble
    user{%
      /Subtype/Link%
      /A<<%
        /Type/Action%
        /S/GoToR%
        /D##1%
        /F(\hl@opt@file)%
        \ifx\hl@opt@newwin\empty \else
          /NewWindow \ifcase\hl@opt@newwin false\else true\fi
        \fi
      >>%
    }%
  }%
  \def\hl@typeh@url{%
    \hlstart@preamble
    user{%
      /Subtype/Link%
      /A<<%
        /Type/Action%
        /S/URI%
        /URI(\@hllabel)%
      >>%
    }%
  }%
  \def\@hlend{\pdfendlink\endgroup}
\else 
  \message{Eplain warning: `pdftex' hyperlink driver: PDF output is^^J
           \space not enabled, falling back on `nolinks' driver.}%
  \hldriver@nolinks
\fi
}%
\def\hldriver@dvipdfm{%
  \def\hldest@type{xyz}%
  \let\hldest@opt@left   \empty
  \let\hldest@opt@top    \empty
  \let\hldest@opt@right  \empty
  \let\hldest@opt@bottom \empty
  \let\hldest@opt@zoom   \empty
  \let\hldest@opt@cmd    \empty
  \def\hldest@driver{%
    \ifx\@hltype\hl@raw@word
      \csname \hldest@opt@cmd \endcsname
    \else
      \def\hldest@preamble{%
        pdf: dest (\@hllabel) [@thispage
      }%
      \csname hldest@typeh@\@hltype\endcsname
    \fi
  }%
  \let\hldest@typeh@raw\empty
  \def\hldest@typeh@fit{%
    \special{\hldest@preamble /Fit]}%
  }%
  \def\hldest@typeh@fith{%
    \special{\hldest@preamble /FitH
      \ifx\hldest@opt@top\empty @ypos \else \hldest@opt@top \fi]}%
  }%
  \def\hldest@typeh@fitv{%
    \special{\hldest@preamble /FitV
      \ifx\hldest@opt@left\empty @xpos \else \hldest@opt@left \fi]}%
  }%
  \def\hldest@typeh@fitb{%
    \special{\hldest@preamble /FitB]}%
  }%
  \def\hldest@typeh@fitbh{%
    \special{\hldest@preamble /FitBH
      \ifx\hldest@opt@top\empty @ypos \else \hldest@opt@top \fi]}%
  }%
  \def\hldest@typeh@fitbv{%
    \special{\hldest@preamble /FitBV
      \ifx\hldest@opt@left\empty @xpos \else \hldest@opt@left \fi]}%
  }%
  \def\hldest@typeh@fitr{%
    \special{\hldest@preamble /FitR
      \ifx\hldest@opt@left\empty @xpos\else\hldest@opt@left\fi\space
      \ifx\hldest@opt@bottom\empty @ypos\else\hldest@opt@bottom\fi\space
      \ifx\hldest@opt@right\empty @xpos\else\hldest@opt@right\fi\space
      \ifx\hldest@opt@top\empty @ypos\else\hldest@opt@top \fi]}%
  }%
  \def\hldest@typeh@xyz{%
    \begingroup
      \ifx\hldest@opt@zoom\empty
        \count1=\z@ \count2=\z@
      \else
        \count2=\hldest@opt@zoom
        \count1=\count2 \divide\count1 by 1000
        \count3=\count1 \multiply\count3 by 1000
        \advance\count2 by -\count3
      \fi
      \special{\hldest@preamble /XYZ
        \ifx\hldest@opt@left\empty @xpos\else\hldest@opt@left\fi\space
        \ifx\hldest@opt@top\empty @ypos\else\hldest@opt@top\fi\space
        \the\count1.\the\count2]}%
    \endgroup
  }%
  \def\hl@type{name}%
  \ifx\hl@type@url\empty
    \def\hl@type@url{url}%
  \fi
  \ifx\hl@type@hrefext\empty
    \def\hl@type@hrefext{url}%
  \fi
  \def\hl@opt@bstyle  {S}%
  \def\hl@opt@bwidth  {1}%
  \let\hl@opt@bcolor  \empty
  \let\hl@opt@hlight  \empty
  \let\hl@opt@bdash   \empty
  \let\hl@opt@pagefit \empty
  \let\hl@opt@cmd     \empty
  \let\hl@opt@file    \empty
  \let\hl@opt@newwin  \empty
  \def\hl@driver{%
    \ifx\@hltype\hl@raw@word
      \csname \hl@opt@cmd \endcsname
    \else
      \let\hl@BSspec\relax 
      \ifx\hl@opt@bstyle \empty
        \ifx\hl@opt@bwidth \empty
          \ifx\hl@opt@bdash \empty
            \let\hl@BSspec\empty 
          \fi
        \fi
      \fi
      \def\hlstart@preamble{%
        pdf: beginann
          <<%
            /Type/Annot%
            /Subtype/Link%
            \ifx\hl@opt@bcolor\empty\else /C[\hl@opt@bcolor]\fi
            \ifx\hl@opt@hlight\empty\else /H/\hl@opt@hlight\fi
            \ifx\hl@BSspec\relax
              /BS<<%
                /Type/Border%
                \ifx\hl@opt@bstyle\empty\else /S/\hl@opt@bstyle\fi
                \ifx\hl@opt@bwidth\empty\else /W \hl@opt@bwidth\fi
                \ifx\hl@opt@bdash\empty \else /D[\hl@opt@bdash]\fi
              >>%
            \fi
      }%
      \csname hl@typeh@\@hltype\endcsname
    \fi
  }%
  \let\hl@typeh@raw\empty
  \def\hl@typeh@name{\special{\hlstart@preamble /Dest(\@hllabel)>>}}%
  \def\hl@typeh@page{%
    \count@=\@hllabel
    \advance\count@ by-1
    \special{%
      \hlstart@preamble
      /Dest[\the\count@
            \ifx\hl@opt@pagefit\empty/Fit\else\hl@opt@pagefit\fi]%
     >>%
    }%
  }%
  \def\hl@typeh@filename{\hl@file{(\@hllabel)}}%
  \def\hl@typeh@filepage{%
    \count@=\@hllabel
    \advance\count@ by-1
    \hl@file{%
      [\the\count@ \ifx\hl@opt@pagefit\empty/Fit\else\hl@opt@pagefit\fi]%
    }%
  }%
  \def\hl@file##1{%
    \special{%
      \hlstart@preamble
      /A<<%
        /Type/Action%
        /S/GoToR%
        /D##1%
        /F(\hl@opt@file)%
        \ifx\hl@opt@newwin\empty \else
          /NewWindow \ifcase\hl@opt@newwin false\else true\fi
        \fi
      >>%
     >>%
    }%
  }%
  \def\hl@typeh@url{%
    \special{%
      \hlstart@preamble
      /A<<%
        /Type/Action%
        /S/URI%
        /URI(\@hllabel)%
      >>%
     >>%
    }%
  }%
  \def\@hlend{\endgroup}
}%
\def\href{%
  \bgroup
    \uncatcodespecials
    \catcode`\{=1 \catcode`\}=2
    \@href
}%
\def\@href#1{
  \egroup
  \edef\@hreftmp{\ifempty{#1}{}\fi}
  \expandafter\@@href\@hreftmp#1\@@
}%
\def\href@end@int{\hlend@impl{hrefint}}%
\def\href@end@ext{\hlend@impl{hrefext}}%
\def\@@href#1#2\@@{%
  \def\@hreftmp{#1}%
  \ifx\@hreftmp\hlhash
    \let\href@end\href@end@int
    \hlstart@impl{hrefint}{#2}%
  \else
    \let\href@end\href@end@ext
    \hlstart@impl{hrefext}{#1#2}%
  \fi
  \@@@href
}%
\def\@@@href{%
  \futurelet\@hreftmp\href@
}%
\def\href@{%
  \ifcat\bgroup\noexpand\@hreftmp
    \let\@hreftmp\href@@
  \else
    \let\@hreftmp\href@@@
  \fi
  \@hreftmp
}%
\def\href@@{\bgroup\aftergroup\href@end \let\@hreftmp}%
\def\href@@@#1{#1\href@end}%
\def\hldeston{\errmessage{Please enable hyperlinks with
  \string\enablehyperlinks\space before using hyperlink commands
  (consider selecting the `nolinks' driver to ignore all hyperlink
  commands in your document)}}%
\let\hldestoff\hldeston \let\hlon\hldeston \let\hloff\hldeston
\let\hlstart\hldeston \let\hlend\hldeston \let\hldest\hldeston
\let\hl@setparam\hldeston
\@hloff[*]\@hldestoff[*]%
\newif\ifusepkg@miniltx@loaded
\newcount\usepkg@recursion@level
\def\usepkg@rcrs{\the\usepkg@recursion@level}%
\let\usepkg@at@begin@document\empty
\let\usepkg@at@end@of@package\empty
\def\usepkg@word@autopict{autopict}%
\def\usepkg@word@psfrag{psfrag}%
\long\def\beginpackages#1\endpackages{%
  \let\usepackage\real@usepackage
  \let\DoNotLoadEpstopdf=t
  \let\eplaininput=\input
  #1%
  \usepkg@at@begin@document
  \global\let\usepkg@at@begin@document\empty
  \global\let\usepackage\fake@usepackage
  \let\packageinput=\input
  \let\input=\eplaininput
}%
\def\fake@usepackage{\errmessage{You should not use \string\usepackage\space outside of^^J
  \@spaces\string\beginpackages...\string\endpackages\space environment}%
}%
\def\eplain@RequirePackage{%
  \global\ece\let{usepkg@save@pkg\usepkg@rcrs}\usepkg@pkg
  \global\ece\let{usepkg@save@options\usepkg@rcrs}\usepkg@options
  \global\ece\let{usepkg@save@date\usepkg@rcrs}\usepkg@date
  \global\ece\let{usepkg@at@end@of@package\usepkg@rcrs}\usepkg@at@end@of@package
  \global\advance\usepkg@recursion@level by\@ne
  \real@usepackage
}%
\let\usepackage\fake@usepackage
\def\real@usepackage{\@getoptionalarg\@finusepackage}%
\def\@finusepackage#1{%
  \let\usepkg@options\@optionalarg
  \ifempty{#1}%
    \errmessage{No packages specified}%
  \fi
  \ifusepkg@miniltx@loaded \else
    \testfileexistence[miniltx]{tex}%
    \if@fileexists
      \input miniltx.tex
      \global\usepkg@miniltx@loadedtrue
      \global\let\RequirePackage\eplain@RequirePackage
      \global\let\DeclareOption\eplain@DeclareOption
      \global\let\PassOptionsToPackage\eplain@PassOptionsToPackage
      \global\let\ExecuteOptions\eplain@ExecuteOptions
      \global\let\ProcessOptions\eplain@ProcessOptions
      \global\let\AtBeginDocument\eplain@AtBeginDocument
      \global\let\AtEndOfPackage\eplain@AtEndOfPackage
      \global\let\ProvidesFile\eplain@ProvidesFile
      \global\let\ProvidesPackage\eplain@ProvidesPackage
    \else
      \errmessage{miniltx.tex not found, cannot load LaTeX packages}%
    \fi
  \fi
  \@ifnextchar[
    {\@finfinusepackage{#1}}%
    {\@finfinusepackage{#1}[]}%
}%
\def\@finfinusepackage#1[#2]{%
  \edef\usepkg@date{#2}%
  \let\usepkg@load@list\empty
  \for\usepkg@pkg:=#1\do{%
    \toks@=\expandafter{\usepkg@load@list}%
    \edef\usepkg@load@list{%
      \the\toks@
      \noexpand\usepkg@load@pkg{\usepkg@pkg}%
    }%
  }%
  \usepkg@load@list
  \ifnum\usepkg@recursion@level>0
    \global\advance\usepkg@recursion@level by\m@ne
    \expandafter\let\expandafter\usepkg@pkg\csname usepkg@save@pkg\usepkg@rcrs\endcsname
    \expandafter\let\expandafter\usepkg@options\csname usepkg@save@options\usepkg@rcrs\endcsname
    \expandafter\let\expandafter\usepkg@date\csname usepkg@save@date\usepkg@rcrs\endcsname
    \expandafter\let\expandafter\usepkg@at@end@of@package\csname usepkg@at@end@of@package\usepkg@rcrs\endcsname
    \global\ece\let{usepkg@save@pkg\usepkg@rcrs}\undefined
    \global\ece\let{usepkg@save@options\usepkg@rcrs}\undefined
    \global\ece\let{usepkg@save@date\usepkg@rcrs}\undefined
    \global\ece\let{usepkg@at@end@of@package\usepkg@rcrs}\undefined
  \fi
}%
\def\usepkg@load@pkg#1{%
  \def\usepkg@pkg{#1}%
  \ifx\usepkg@pkg\usepkg@word@autopict
    \testfileexistence[picture]{tex}%
    \if@fileexists \else
      \errmessage{Loader `picture.tex' for package `\usepkg@pkg' not found}%
    \fi
  \else
    \ifx\usepkg@pkg\usepkg@word@psfrag
      \testfileexistence[psfrag]{tex}%
      \if@fileexists \else
        \errmessage{Loader `psfrag.tex' for package `\usepkg@pkg' not found}%
      \fi
    \fi
  \fi
  \ifundefined{ver@\usepkg@pkg.sty}%
    \expandafter\@finusepkg@load@pkg
  \else
    \immediate\write-1{^^J\linenumber Eplain: package `\usepkg@pkg' already
             loaded, skipping reloading}%
  \fi
}%
\def\@finusepkg@load@pkg{%
  \testfileexistence[\usepkg@pkg]{sty}%
  \if@fileexists \else
    \errmessage{Package `\usepkg@pkg' not found}%
  \fi
  \expandafter\let\expandafter\temp\csname usepkg@options@\usepkg@pkg\endcsname
  \ifx\temp\relax
    \let\temp\empty
  \fi
  \ifx\temp\empty
    \let\temp\usepkg@options
  \else
    \ifx\usepkg@options\empty \else
      \edef\temp{\temp,\usepkg@options}%
    \fi
  \fi
  \global\ece\let{usepkg@options@\usepkg@pkg}\temp
  \let\usepackage\eplain@RequirePackage
  \global\let\usepkg@at@end@of@package\empty
  \ifx\usepkg@pkg\usepkg@word@autopict
    \input picture.tex
  \else
    \ifx\usepkg@pkg\usepkg@word@psfrag
      \input \usepkg@pkg.tex
    \else
      \input \usepkg@pkg.sty
    \fi
  \fi
  \usepkg@at@end@of@package
  \global\let\usepkg@at@end@of@package\empty
  \let\usepackage\real@usepackage
  \global\ece\let{usepkg@options@\usepkg@pkg}\undefined
  \def\Url@HyperHook##1{\hlstart@impl{url}{\Url@String}##1\hlend@impl{url}}%
}%
\def\eplain@DeclareOption#1#2{%
  \toks@{#2}%
  \expandafter\xdef\csname usepkg@option@\usepkg@pkg @#1\endcsname{\the\toks@}%
}%
\def\eplain@PassOptionsToPackage#1#2{%
  \ifempty{#1}\else
    \for\usepkg@temp:=#2\do{%
      \expandafter\let\expandafter\temp\csname usepkg@options@\usepkg@temp\endcsname
      \ifx\temp\relax
        \let\temp\empty
      \fi
      \ifx\temp\empty
        \edef\temp{#1}%
      \else
        \edef\temp{\temp,#1}%
      \fi
      \global\ece\let{usepkg@options@\usepkg@temp}\temp
    }%
  \fi
}%
\def\usepkg@exec@options#1{%
  \for\CurrentOption:=#1\do{%
    \expandafter\let\expandafter\usepkg@temp
      \csname usepkg@option@\usepkg@pkg @\CurrentOption\endcsname
    \ifx\usepkg@temp\relax
      \expandafter\let\expandafter\temp\csname usepkg@option@\usepkg@pkg @*\endcsname
      \ifx\temp\relax
        \errmessage{Unknown option `\CurrentOption' to package `\usepkg@pkg'}%
      \else
        \temp
      \fi
    \else
      \usepkg@temp
    \fi
  }%
}%
\let\eplain@ExecuteOptions\usepkg@exec@options
\def\eplain@ProcessOptions{%
  \expandafter\usepkg@exec@options\csname usepkg@options@\usepkg@pkg\endcsname
}%
\def\usepkg@accumulate@cmds#1#2{%
  \toks@=\expandafter{#1}%
  \toks@ii={#2}%
  \xdef#1{\the\toks@\the\toks@ii}%
}%
\def\eplain@AtBeginDocument{\usepkg@accumulate@cmds\usepkg@at@begin@document}%
\def\eplain@AtEndOfPackage{\usepkg@accumulate@cmds\usepkg@at@end@of@package}%
\def\eplain@ProvidesPackage#1{%
  \@ifnextchar[
    {\eplain@pr@videpackage{#1.sty}}{\eplain@pr@videpackage#1[]}%
}%
\def\eplain@pr@videpackage#1[#2]{%
  \wlog{#1: #2}%
  \expandafter\xdef\csname ver@#1\endcsname{#2}%
  \@ifl@t@r{#2}\usepkg@date{}%
    {\message{Warning: you have requested package `\usepkg@pkg', version \usepkg@date,^^J
       \@spaces but only version `\csname ver@#1\endcsname' is available.}}%
}%
\def\eplain@ProvidesFile#1{%
  \@ifnextchar[
    {\eplain@pr@videfile{#1}}{\eplain@pr@videfile#1[]}%
}%
\def\eplain@pr@videfile#1[#2]{%
  \wlog{#1: #2}%
  \expandafter\xdef\csname ver@#1\endcsname{#2}%
}%
\def\@ifl@ter#1#2{%
  \expandafter\@ifl@t@r
    \csname ver@#2.#1\endcsname
}%
\def\@ifl@t@r#1#2{%
  \ifnum\expandafter\@parse@version#1//00\@nil<%
        \expandafter\@parse@version#2//00\@nil
    \expandafter\@secondoftwo
  \else
    \expandafter\@firstoftwo
  \fi
}%
\def\@parse@version#1/#2/#3#4#5\@nil{#1#2#3#4 }%

\def\strip@prefix#1>{}%
\def\@ifpackageloaded#1{%
  \expandafter\ifx\csname ver@#1.sty\endcsname\relax
    \expandafter\@secondoftwo
  \else
    \expandafter\@firstoftwo
  \fi
}%
\long\def\g@addto@macro#1#2{%
  \begingroup
    \toks@\expandafter{#1#2}%
    \xdef#1{\the\toks@}%
  \endgroup
}%
\def\PackageWarning#1#2{%
  \begingroup
    \newlinechar=10 %
    \def\MessageBreak{%
      ^^J(#1)\@spaces\@spaces\@spaces\@spaces
    }%
    \immediate\write16{^^JPackage #1 Warning: #2\on@line.^^J}%
  \endgroup
}%
\def\PackageWarningNoLine#1#2{%
  \PackageWarning{#1}{#2\@gobble}%
}%
\def\on@line{ on input line \the\inputlineno}%
\def\@spaces{\space\space\space\space}%
\def\@inmatherr#1{%
   \relax
   \ifmmode
     \errmessage{The command is invalid in math mode}%
   \fi
}%
\let\protected@edef\edef
\let\wlog = \@plainwlog
\catcode`@ = \@eplainoldatcode
\def\eplain{t}%
{\edef\plainversion{\fmtversion}%
 \xdef\fmtversion{3.3: 21 July 2009 (and plain \plainversion)}%
}%

\lefteqnumbers
   \def\testd{oui}
   \def\choixlat{\ifx\numadroite\testd\righteqnumbers
            \else  \lefteqnumbers\fi}
    \choixlat

\catcode`@=\letter
\def\@eplainfootnote#1{\let\@sf\empty 
  \ifhmode\edef\@sf{\spacefactor\the\spacefactor}\/\fi
  \global\advance\hlfootlabelnumber by 1
  \hlstart@impl{foot}{\hlfootlabel}%
  \hldest@impl{footback}{\hlfootbacklabel}%
  \hbox{$^{(#1)}$}%
  \hlend@impl{foot}%
  \@sf\vfootnote{#1.}%
}%
\catcode`@=\other

  \interfootnoteskip=0pt
  \let\note=\numberedfootnote
  \everyfootnote={\eightpoint\leftskip=5truemm\rightskip5truemm}
  
  \hsize150truemm\vsize 240truemm\hoffset=5truemm
  
  \def\dimart{\hsize126truemm\vsize186truemm\hoffset16truemm\voffset=24truemm}
  \pretolerance=500\tolerance=1000\brokenpenalty=5000
  \parindent3mm
  
  \countdef\temps=170
  \temps=\time
  \countdef\nminutes=171{\nminutes = \time}
  \countdef\nheure=172
  \def\heure{\begingroup                     
     \temps = \time \divide\temps by 60
     \nheure = \temps                        
     \nminutes = \time
     \multiply\temps by 60
     \advance\nminutes by -\temps            
     \ifnum\nminutes<10 \toks1 = {0}%
     \else\toks1 = {}%
     \fi
     \number\nheure h\the\toks1 \number\nminutes  
  \endgroup}%

  \newcount\chstart
  \chstart=\pageno
 \headline={\ifnum\pageno=\chstart {\hfill} \else {\hss \tenrm --\ \folio\ --\hss}\fi}
  \footline={\hfill}
  \normalbaselines
  \frenchspacing
    \def\dater{\vglue-10mm\rightline{(\the\day/\the\month/\the\year)}}
  \def\dateheure{\vglue-10mm\rightline{(\the\day/\the\month/\the\year,\ \heure)}}

  \newif\ifpagetitre \pagetitretrue
\newtoks\hautpagetitre \hautpagetitre={\hfill}
\newtoks\baspagetitre \baspagetitre={\hfill}
\newtoks\auteurcourant \auteurcourant={\hfill}
\newtoks\titrecourant \titrecourant={\hfill}
\newtoks\hautpagegauche
\newtoks\hautpagedroite
\newtoks\hautpagemilieu
\hautpagemilieu={\tenrm\hfil -- \folio\ -- \hfil}
\hautpagegauche={\ifx\midfolio\oui\the\hautpagemilieu\else\tenrm\folio\hfill\the\auteurcourant\hfill\fi}
\hautpagedroite={\ifx\midfolio\oui\the\hautpagemilieu\else\hfill\the\titrecourant\hfill\tenrm\folio\fi}
\newtoks\baspagegauche \baspagegauche={\hfil}
\newtoks\baspagedroite \baspagedroite={\hfil}
\headline={\ifpagetitre\the\hautpagetitre
\else\ifodd\pageno\the\hautpagedroite\else\the\hautpagegauche\fi\fi }
\footline={\ifpagetitre\the\baspagetitre
\else\ifodd\pageno\the\baspagedroite
\else\the\baspagegauche\fi\fi \global\pagetitrefalse}

\def\pageblanche{\vfill\eject\pagetitretrue
\null\vfill\eject
\pagetitretrue
}
\def\chgtpage{\ifodd\pageno \else
\pageblanche \fi \pagetitretrue\titreun=0\footnotenumber=0}

\def\majnombres{\ifodd\pageno \else
\pageblanche \fi \pagetitretrue\hautpoly\titreun=0\footnotenumber=0}

\def\hautspages#1#2{\auteurcourant={\ninepcap#1}\titrecourant={\nineit#2}}

\ifnum\chstart=\pageno \pagetitretrue\fi
  


  \def\leftnote#1{\vadjust{\setbox1=\vtop{\hsize 20mm\parindent=0pt\eightpoint
  \baselineskip=9pt\rightskip=4mm plus 4mm\vskip-4.7mm#1}\hbox{\kern-2cm\smash{\box1}}}}

  
  \def\raggedcenter{\leftskip=20pt plus 10em  
       \rightskip=\leftskip 
        \parfillskip=0pt 
         \spaceskip=.3333em \xspaceskip=.5em 
          \pretolerance=9999 \tolerance=9999
           \hyphenpenalty=9999 \exhyphenpenalty=9999 }
           
  \def\titrecentre#1{{\parindent0mm\raggedcenter
       \spaceskip=.6em plus .2em minus .2em\xspaceskip=.6em plus .2em minus .2em
        \tit#1\par}}
        

\def\boxit#1#2{\vbox{\hrule\hbox{\vrule\vbox spread#1{\vfil\hbox spread#1{\hfil#2\hfil}\vfil}%
\vrule}\hrule}}


  \def\oui{oui}
  
   \def\fontetitreun{\ifx\paradouze\oui\douzepts\gpdouze\twelvebf\textfont1=\twelveib\else
\quatorzepts\gpquatorze\fourteenbf\fi}

\def\fontetitreunl{\douzepts\textfont1=\twelveib\scriptfont1=\tenib\fourteenti}
 
 \def\fontetitredeux{\textfont1=\eleveni\ifx\paradouze\oui\onzepts\scriptfont1=\ninei\elevenit\else
                        \douzepts\twelveit\fi}
 
   \def\fontetitredeuxb{\ifx\paradouze\oui\onzepts\eleventi\gponze\textfont1=\elevenib\scriptfont1=\nineib
                         \else\douzepts\twelveti\scriptfont1=\twelveib\scriptfont1=\tenib\gpdouze\fi}
                         
\def\fontetitredeuxl{\onzepts\textfont1=\elevenbf\scriptfont1=\ninebf\twelvebf}
  
\def\fontetitretrois{\textfont0=\elevenrm\scriptfont0=\eightrm\textfont1=\eleveni
                      \scriptfont1=\eighti\scriptscriptfont1=\sixi\elevenit}
                      
\def\fontetitrequatre{\textfont0=\elevenrm\scriptfont0=\eightrm\textfont1=\eleveni
                      \scriptfont1=\eighti\scriptscriptfont1=\sixi\elevenrm}
  
  \newcount\titreun\titreun=0
  \newcount\titredeux\titredeux=0
  \newcount\titretrois\titretrois=0
  \newcount\titrequatre\titrequatre=0
  \newcount\enonce\enonce=0
  
  \def\incr#1{\global\advance#1 by 1 {\the #1}}
  \def\avance#1{\global\advance#1 by 1}
  \def\init#1{\global#1=0}
  
  \long\def\Indentation#1#2{\setbox10=\hbox{\fontetitreun#1}
                        \ifdim\wd10 < 4mm
                         \setbox10=\hbox to 4mm{\box10\hfill}
                       \else\ifdim\wd10 < 6mm
                         \setbox10=\hbox to 6mm{\box10\hfill}
                        \else\ifdim\wd10 < 8mm
                         \setbox10=\hbox to 8mm{\box10\hfill}
                       \else\ifdim\wd10 < 12mm
                         \setbox10=\hbox to 12mm{\box10\hfill}
                       \fi\fi\fi\fi
                       \dimen10=\hsize
                       \advance \dimen10 by -\wd10
                       \noindent \box10 %
                       \ignorespaces
                       \hbox{\vtop{\hsize=\dimen10\raggedright\noindent\fontetitreun#2}}}

  \long\def\paraun#1{\removelastskip\par\bigskip\goodbreak\vskip0pt plus.01\vsize\penalty-100
                \vskip0pt plus-.01\vsize
                  \init{\titredeux}\ifnum\optionparag=1{\init\eqnumber\init\enonce}\else{}\fi
                  \goodbreak{\fontetitreun
                    \Indentation{\incr{\titreun}.\ }{\fontetitreun #1\par}}\nobreak\medskip}

 %
 %
 \long\def\paraunc#1{\removelastskip\par\bigskip\goodbreak\vskip0pt plus.01\vsize\penalty-100
                \vskip0pt plus-.01\vsize
                  \init{\titredeux}
                 \ifnum\optionparag=1{\init{\eqnumber}\init\enonce}\else{}\fi
                  \goodbreak
                    {\parindent0mm\raggedcenter\fontetitreun\incr{\titreun}.\ 
                     \fontetitreun #1\par}\nobreak\medskip}
                     
\newtoks\titreunl
\titreunl={\ifnum\titreun=1{I}\fi%
\ifnum\titreun=2{II}\fi%
\ifnum\titreun=3{III}\fi%
\ifnum\titreun=4{IV}\fi%
\ifnum\titreun=5{V}\fi%
\ifnum\titreun=6{VI}\fi%
\ifnum\titreun=7{VII}\fi%
\ifnum\titreun=8{VIII}\fi%
\ifnum\titreun=9{IX}\fi%
\ifnum\titreun=10{X}\fi%
\ifnum\titreun=11{XI}\fi%
\ifnum\titreun=12{XII}\fi%
\ifnum\titreun=13{XIII}\fi%
}
\long\def\paraunl#1{\removelastskip\par\bigskip\goodbreak\vskip0pt plus.01\vsize\penalty-100
                \vskip0pt plus-.01\vsize
                  \init{\titredeux}\ifnum\optionparag=1{\init\eqnumber\init\enonce}\else{}\fi
                  \goodbreak{\fontetitreunl
                    \Indentation{\global\advance\titreun by 1{\the\titreunl}.\ }{\fontetitreunl #1\par}}\nobreak\smallskip}

  
  \long\def\paradeux#1{\init{\titretrois}\vskip0pt plus.01\vsize\penalty-10
                \vskip0pt plus-.01\vsize\ifx \elie\oui\medskip\ifnum\titredeux>0\medskip\fi\fi
                 \Indentation{\fontetitredeux\the\titreun${\cdot}$\incr{\titredeux}.\ }
                              {\fontetitredeux\textfont1=\eleveni#1}\nobreak\smallskip}
  
  \long\def\paradeuxb#1{\init{\titretrois}\vskip0pt plus.001\vsize\penalty-10
                \vskip0pt plus-.01\vsize{\ifx \elie\oui\medskip\ifnum\titredeux>0\medskip\fi\fi
                  \Indentation
  {\fontetitredeuxb\the\titreun${\cdot}$\incr{\titredeux}.}{ \fontetitredeuxb#1}}\nobreak\smallskip
  \par}

\newtoks\titredeuxl
\titredeuxl={\ifnum\titredeux=1{A}\fi%
\ifnum\titredeux=2{B}\fi%
\ifnum\titredeux=3{C}\fi%
\ifnum\titredeux=4{D}\fi%
\ifnum\titredeux=5{E}\fi%
\ifnum\titredeux=6{F}\fi%
\ifnum\titredeux=7{G}\fi%
\ifnum\titredeux=8{H}\fi%
\ifnum\titredeux=9{I}\fi%
\ifnum\titredeux=10{J}\fi%
\ifnum\titredeux=11{K}\fi%
\ifnum\titredeux=12{L}\fi%
\ifnum\titredeux=13{M}\fi%
}
 \long\def\paradeuxl#1{\init{\titretrois}\vskip0pt plus.001\vsize\penalty-10
                \vskip0pt plus-.01
                \vsize \bigskip%
                  \Indentation
     {\fontetitredeuxl\global\advance\titredeux by 1
  \quad \the\titreunl${\cdot}$\the\titredeuxl.}{ \fontetitredeuxl#1}
  \removelastskip\medskip\nobreak\par}
  

  \long\def\paratrois#1{\init{\titrequatre}\ifdim\lastskip<\smallskipamount
                \removelastskip\smallskip\fi
                 \vskip0pt plus.01\vsize\penalty-10
                  \vskip0pt
plus-.01\vsize{\ifx \elie\oui\ifnum\titretrois>0\medskip\fi\fi
\Indentation{\fontetitretrois\the\titreun${\cdot}$\the\titredeux${\cdot}$\incr{\titretrois}.\ }
  {\hskip0mm\fontetitretrois#1}\smallskip\nobreak\par}}
  
  
  \long\def\paratroisl#1{\init{\titrequatre}\ifdim\lastskip<\smallskipamount
                \removelastskip\smallskip\fi
                 \vskip0pt plus.01\vsize\penalty-10
                  \vskip0pt
plus-.01\vsize\ifx \elie\oui\bigskip\fi
\Indentation{\fontetitretrois\quad \quad \the\titreunl{${\cdot}$}\the\titredeuxl${\cdot}$\incr{\titretrois}.\ }
  {\hskip0mm\fontetitretrois#1}\smallskip\nobreak\par}


  \long\def\paraquatre#1{\ifdim\lastskip<\smallskipamount
                \removelastskip\smallskip\fi
                 \vskip0pt plus.01\vsize\penalty-10
                  \vskip0pt
                  plus-.01\vsize\par
                \medskip
\Indentation{\fontetitrequatre \the\titreun{${\cdot}$}\the\titredeux${\cdot}$\the\titretrois${\cdot}$\incr{\titrequatre}.\ }
{\hskip0mm\fontetitrequatre#1}\smallskip\nobreak\par}


\newtoks\titrequatrel
\titrequatrel={\ifnum\titrequatre=1{a}\fi%
\ifnum\titrequatre=2{b}\fi%
\ifnum\titrequatre=3{c}\fi%
\ifnum\titrequatre=4{d}\fi%
\ifnum\titrequatre=5{e}\fi%
\ifnum\titrequatre=6{f}\fi%
\ifnum\titrequatre=7{g}\fi%
\ifnum\titrequatre=8{h}\fi%
\ifnum\titrequatre=9{i}\fi%
\ifnum\titrequatre=10{j}\fi%
\ifnum\titrequatre=11{k}\fi%
\ifnum\titrequatre=12{l}\fi%
\ifnum\titrequatre=13{m}\fi%
}
\long\def\paraquatrel#1{\ifdim\lastskip<\smallskipamount
                \removelastskip\smallskip\fi
                 \vskip0pt plus.01\vsize\penalty-10
                  \vskip0pt
                  plus-.01\vsize{\bigskip
\Indentation{\global\advance\titrequatre by 1
\fontetitrequatre\quad \quad \quad \the\titreunl${\cdot}$\the\titredeuxl${\cdot}$\the\titretrois${\cdot}$\the\titrequatrel.\ }
{\hskip0mm\fontetitrequatre#1}\smallskip\nobreak\par}}

\ifx\optionkeys\oui
\def\drefun#1{\definexref{¤#1}{{\the\titreun}}{}} 
\def\drefdeux#1{\definexref{¤#1}{{\the\titreun}.{\the\titredeux}}{}}
\def\dreftrois#1{\definexref{¤#1}{{\the\titreun}.{\the\titredeux}.{\the\titretrois}}{}}
\else
\def\drefun#1{\definexref{prg#1}{{\the\titreun}}{}} 
\def\drefdeux#1{\definexref{prg#1}{{\the\titreun}.{\the\titredeux}}{}}
\def\dreftrois#1{\definexref{prg#1}{{\the\titreun}.{\the\titredeux}.{\the\titretrois}}{}}
\fi

%


  \long\def\propdeux#1#2#3#4{%
       \avance{\enonce}
       \leavevmode\edef\temp{#2}%
         \ifx\temp\empty 
          \else
           \definexref{#2}{#1~{\the\titreun.\the\enonce}}{enonces}
            \definexref{s#2}{{\the\titreun.\the\enonce}}{enonces}
             \fi
\medbreak
      \noindent{\bf#1\ {\bf\the\titreun.\the\enonce{#3}.}\enspace}{\sl#4\par}%
      \ifdim\lastskip<\medskipamount \removelastskip\penalty55\medskip\fi
   }

  \long\def\propun#1#2#3#4{%
      \avance{\enonce}
       \leavevmode\edef\temp{#2}%
        \ifx\temp\empty 
          \else
           \definexref{#2}{#1~{\the\enonce}}{enonces}
            \definexref{{s#2}}{{\the\enonce}}{enonces}
             \fi
   \medbreak
     \noindent{\bf#1\ {\bf\the\enonce{#3}.}\enspace}{\sl#4\par}%
     \ifdim\lastskip<\medskipamount \removelastskip\penalty55\medskip\fi
  }
  
  \long\def\prop#1#2#3#4{\ifnum\optionparag=1
                          \propdeux{#1}{#2}{\textfont1=\elevenib#3}{#4} \else\propun{#1}{#2}{\textfont1=\elevenib#3}{#4}\fi}

  \long\def\propt#1#2#3{\ifx\tpf\oui \prop{Th\'eo\-r\`eme}{#1}{#2}{#3}\par
                       \else\prop{Theorem}{#1}{#2}{#3}\par\fi}
  \long\def\Propt#1#2{\propt{#1}{}{#2}}
  \long\def\propl#1#2#3{\ifx\tpf\oui\prop{Lem\-me}{#1}{#2}{#3}\par
                         \else\prop{Lemma}{#1}{#2}{#3}\par\fi}
  \long\def\Propl#1#2{\propl{#1}{}{#2}}
  \long\def\propc#1#2#3{\ifx\tpf\oui\prop{Corol\-laire}{#1}{#2}{#3}\par
                         \else\prop{Corollary}{#1}{#2}{#3}\par\fi}
  \long\def\Propc#1#2{\propc{#1}{}{#2}}
  \long\def\propp#1#2#3{\prop{Pro\-po\-si\-tion}{#1}{#2}{#3}\par}
  \long\def\Propp#1#2{\propp{#1}{}{#2}} 
  \long\def\propd#1#2#3{\ifx\tpf\oui\prop{D\'efi\-nition}{#1}{#2}{#3}\par
                       \else\prop{Definition}{#1}{#2}{#3}\par\fi} 
  \long\def\Propd#1#2{\propd{#1}{}{#2}}
  \long\def\proptd#1#2#3{\ifx\tpf\oui\prop{Th\'eor\`eme et d\'efi\-nition}{#1}{#2}{#3}\par
                       \else\prop{Theorem and definition}{#1}{#2}{#3}\par\fi}


  
  \newcount\optionparag\optionparag=1
  
  \long\def\section#1#2{\ifnum\optionparag=1 \paraun{#2} 
                        \else\goodbreak{\fontetitreun
                    \Indentation{#1.\ }{#2}}\nobreak\smallskip\fi}

  \def\eqconstruct#1{\ifnum\optionparag=1{\the\titreun\hbox{$\cdot$}#1}\else{#1}\fi}
  \def\eqprint#1{{\rm (#1)}}

  
  
  \def\numref{oui}  
  
  \newcount\mesref\mesref=0 
  \def\defbib#1{\ifx\numref\oui\global\advance\mesref by 1 \definexref{#1}{{\the
                 \mesref}}{}\else\definexref{#1}{#1}{}\fi}
  \def\bibtem#1{\defbib{#1}\item{\citer{#1}}}
  \def\citer#1{[\ref{#1}]}

  
  \font\seventeenmsa=msam10 at 17pt    
  \font\fourteenmsa=msam10 at 14pt
  \font\twelvemsa=msam10 at 12pt
  \font\tenmsa=msam10                 
  \font\ninemsa=msam10 at 9pt 
  \font\eightmsa=msam10 at 8pt 
  \font\sevenmsa=msam7 
  \font\sixmsa=msam10 at 6pt
  \font\fivemsa=msam5
  \newfam\msafam\textfont\msafam=\tenmsa\scriptfont\msafam=\sevenmsa\scriptscriptfont\msafam=\fivemsa
  
  \font\seventeenbb=msbm10 at 17pt     
  \font\fourteenbb=msbm10 at 14pt
  \font\twelvebb=msbm10 at 12pt
  \font\tenbb=msbm10                   
  \font\ninebb=msbm10 at 9pt 
  \font\eightbb=msbm10 at 8pt 
  \font\sevenbb=msbm7 
  \font\sixbb=msbm10 at 6pt
  \font\fivebb=msbm5 
  \newfam\bbfam\textfont\bbfam=\tenbb\scriptfont\bbfam=\sevenbb\scriptscriptfont\bbfam=\fivebb
  \def\bb{\fam\bbfam\tenbb}%

  \font\seventeenscaln=eusm10 at 17pt   
  \font\twelvescaln=eusm10 at 12pt
  \font\tenscaln=eusm10                
  \font\ninescaln=eusm10 scaled 900
  \font\eightscaln=eusm10 scaled 800
  \font\sevenscaln=eusm10 scaled 700
  \font\sixscaln=eusm10 scaled 600
   
  \newfam\scalnfam\textfont\scalnfam=\tenscaln\scriptfont\scalnfam=\sevenscaln\scriptscriptfont\scalnfam=\sixscaln
  \def\scaln{\fam\scalnfam\tenscaln}%
  \def\scal{\scaln}
  
  \font\tenscalb=eusb10                

  \font\sevenscalb=eusb10 scaled 700

  \newfam\scalbfam\textfont\scalbfam=\tenscalb\scriptfont\scalbfam=\sevenscalb
  %
  
  %
  %
  \font\fourteenrm=cmr12 scaled 1200
  \font\elevenrm=cmr10 at 11pt
  \font\twelverm=cmr12
  \font\ninerm=cmr9
  \font\eightrm=cmr8      
  \font\sevenrm=cmr7
  \font\sixrm=cmr6

  \font\seventeenpcap=cmcsc10 at 17pt
  \font\tenpcap=cmcsc10                        
  \font\ninepcap=cmcsc9
  \font\eightpcap=cmcsc8
  \font\sevenpcap=cmcsc10 scaled 700
  
  \newfam\pcapfam\textfont\pcapfam=\tenpcap\scriptfont\pcapfam=\sevenpcap
  \def\pcap{\fam\pcapfam\tenpcap}
  
  \font\seventeenrm=cmbx12 scaled 1400

  \font\fourteenbf=cmbx10 scaled 1400
  
  \font\twelvebf=cmbx12
  \font\elevenbf=cmbx10 at 11pt
  \font\ninebf=cmbx9  
  \font\eightbf=cmbx8
  \font\sixbf=cmbx6
  
  \font\tengot=eufm10                           
   
   at 8truept 
  \font\sevengot=eufm7 
   at 6 truept 
   
  \newfam\gotfam
  \textfont\gotfam=\tengot\scriptfont\gotfam=\sevengot
  %

  
  \def\tit{%
  \textfont0=\seventeenrm\scriptfont0=\tenrm\def\rm{\fam0\seventeenrm}%
  \textfont1=\seventeenib\scriptfont1=\twelveib%
  \textfont2=\seventeensy\scriptfont2=\twelvesy\scriptscriptfont2=\ninesy
  \textfont3=\seventeenex
  \textfont\itfam=\seventeenti
  \def\it{\fam\itfam\seventeenti}%
  \textfont\bbfam=\seventeenbb \scriptfont\bbfam=\twelvebb
  \def\bb{\fam\bbfam\seventeenbb}%
  \textfont\msafam=\seventeenmsa\scriptfont\msafam=\twelvemsa
  \textfont\scalnfam=\seventeenscaln
  \def\pcap{\fam\pcapfam\seventeenpcap}
  \normalbaselineskip=25pt\normalbaselines\rm}

  \font\seventeenti=cmbxti10 scaled 1680
  
  \font\fourteenti=cmbxti10 at 14pt
  
  \font\twelveti=cmbxti10 scaled 1200
  \font\eleventi=cmbxti10 at 11pt

  %
  %
  \font\twelveit=cmti12 
  \font\elevenit=cmti10 scaled 1100
  \font\nineit=cmti9
  \font\eightit=cmti8
  \font\sevenit=cmti7

  %
  %
  \font\seventeenib=cmmib10 scaled 1680
  \font\fourteenib=cmmib10 scaled 1400
  \font\twelveib=cmmib10 scaled 1200
  \font\elevenib=cmmib10 scaled 1100
  \font\tenib=cmmib10
  \font\nineib=cmmib10 scaled 900

  %
  %
  
  \font\eleveni=cmmi10 scaled 1100
  \font\ninei=cmmi9
  \font\eighti=cmmi8 
  \font\seveni=cmmi7                            
  \font\sixi=cmmi6
  
  \font\ninesl=cmsl9                    
  \font\eightsl=cmsl8 
  \font\sevensl=cmsl10 at 7pt

  \font\ninett=cmtt9                    
  \font\eighttt=cmtt8
  \font\seventt=cmtt10 scaled 700
  
  \font\fivett=cmtt10 scaled 500
  
  \font\seventeensy=cmsy10 scaled 1680    
  \font\fourteensy=cmsy10 scaled 1400
  \font\twelvesy=cmsy10 scaled 1176
  
  \font\ninesy=cmsy9                      
  \font\eightsy=cmsy8
  \font\sixsy=cmsy6
  \font\seventeenex=cmex10 at 17pt
  \font\fourteenex=cmex10 at 14pt
  \font\twelveex=cmex10 at 12pt
  \font\nineex=cmex10 at 9pt
  \font\eightex=cmex10 at 8pt
  \font\sevenex=cmex10 at 7pt
  \font\sixex=cmex10 at 6pt
  \font\fiveex=cmex10 at 5pt
  
   
  \font\fourteengp=cmmi10 at 14pt
  
  \font\twelvegp=cmmib10 at 12pt
  \font\elevengp=cmmib10 at 11pt
  \font\tengp=cmmib10                          
  \font\ninegp=cmmib10 at 9pt 
  \font\eightgp=cmmib8 
  \font\sevengp=cmmib7 
  \font\sixgp=cmmib6

  \def\gponze{\textfont0=\elevenbf\scriptfont0=\eightbf\scriptscriptfont0=\sixbf
           \textfont1=\elevengp\scriptfont1=\eightgp\scriptscriptfont1=\sixgp}
  \def\gpdouze{\textfont0=\twelvebf\scriptfont0=\tenbf\scriptscriptfont0=\ninebf
           \textfont1=\twelvegp\scriptfont1=\tengp\scriptscriptfont1=\ninegp}        
  
 \def\gpquatorze{\textfont0=\fourteenbf\scriptfont0=\twelvebf\scriptscriptfont0=\elevenbf
           \textfont1=\fourteengp\scriptfont1=\twelvegp\scriptscriptfont1=\elevengp}

  
  \expandafter\chardef\csname pre amssym.def at\endcsname=\the\catcode`\@
  \catcode`\@=11
  \def\undefine#1{\let#1\undefined}
  \def\newsymbol#1#2#3#4#5{\let\next@\relax
   \ifnum#2=\@ne\let\next@\msafam@\else
   \ifnum#2=\tw@\let\next@\bbfam@\fi\fi
   \mathchardef#1="#3\next@#4#5}
  \def\mathhexbox@#1#2#3{\relax
   \ifmmode\mathpalette{}{\m@th\mathchar"#1#2#3}%
   \else\leavevmode\hbox{$\m@th\mathchar"#1#2#3$}\fi}
  \def\hexnumber@#1{\ifcase#1 0\or 1\or 2\or 3\or 4\or 5\or 6\or 7\or 8\or
   9\or A\or B\or C\or D\or E\or F\fi}
  
  \def\setboxz@h{\setbox\z@\hbox}
  \def\wdz@{\wd\z@}
  \def\boxz@{\box\z@}
  
  \edef\msafam@{\hexnumber@\msafam}
  \mathchardef\dabar@"0\msafam@39
  
  \edef\bbfam@{\hexnumber@\bbfam}
  \def\widehat#1{\setboxz@h{$\m@th#1$}%
   \ifdim\wdz@>\tw@ em\mathaccent"0\bbfam@5B{#1}%
   \else\mathaccent"0362{#1}\fi}
  \def\widetilde#1{\setboxz@h{$\m@th#1$}%
   \ifdim\wdz@>\tw@ em\mathaccent"0\bbfam@5D{#1}%
   \else\mathaccent"0365{#1}\fi}
  \newsymbol\leqq 1335          
  \newsymbol\leqslant 1336
  \newsymbol\lessgtr 1337       
  \newsymbol\backprime 1038     
  \newsymbol\risingdotseq 133A  
  \newsymbol\fallingdotseq 133B 
  \newsymbol\succcurlyeq 133C   
  \newsymbol\geqq 133D          
  \newsymbol\geqslant 133E
  \newsymbol\nmid 232D
  \newsymbol\nexists 2040
  \newsymbol\smallsetminus 2272
  \newsymbol\varnothing 203F
  
  \catcode`\@=\active

  
  \catcode`\@=11
  
  \newcount\typofr\typofr=1
  
  \catcode`\;=\active
  \def;{\ifnum\typofr=1\relax\ifhmode\ifdim\lastskip>\z@\unskip\fi
     \kern.2em\fi\string;\else\string;\fi}
  
  \catcode`\:=\active
  \def:{\ifnum\typofr=1\relax\ifhmode\ifdim\lastskip>\z@\unskip\fi
  \penalty\@M\ \fi\string:\else\string:\fi}
  
  \catcode`\!=\active
  \def!{\ifnum\typofr=1\relax\ifhmode\ifdim\lastskip>\z@\unskip\fi
     \kern.2em\fi\string!\else\string!\fi}
  
  \catcode`\?=\active
  \def?{\ifnum\typofr=1\relax\ifhmode\ifdim\lastskip>\z@\unskip\fi
     \kern.2em\fi\string?\else\string?\fi}

  \def\francais{\typofr=1\def\tpf{oui}}
  
  \def\oui{oui}
  \francais
  
  \catcode`\@=12


  
  \def\og{\leavevmode\raise.24ex\hbox{$\scriptscriptstyle\langle\!\langle\>$}}    
  \def\fg{\leavevmode\raise.24ex\hbox{$\scriptscriptstyle\>\rangle\!\rangle$}}    

  \def\d{\,{\rm d}}

  \def\z{{\bb Z}}

  \def\N{{\bb N}}

  \def\PP{{\bb P}}

  \def\A{{\cal A}}

  \def\D{{\cal D}}

  \def\O{{\cal O}}
  \def\P{{\scaln P}}

  \def\frac#1#2{{#1\over #2}}
  \def\di#1#2{\sct#1\atop{\sct#2}}

  \def\qedbox{$\rlap{$\sqcap$}\sqcup$}           
  \def\qed{\nobreak\hfill\penalty250 \hbox{}\nobreak\hfill\qedbox\par\medskip}

  \def\¤{\S\thinspace}

  \def\¥{$\bullet$ }
  
  
  \def\e{{\rm e}}

  \def\epsilon{\varepsilon}

  \def\phi{\varphi}
  \def\theta{\vartheta}
  \def\rho{\varrho}
  \def\dm{{\textstyle{1\over 2}}}

  \def\sct{\scriptstyle}
  \def\pf{\noi{\it Proof. }}
  \def\nid{\ifnum\typofr=1\par\noindent{\it D\'emonstration. }\else\pf\fi}
  \def\noi{\noindent}
  \def\rem{\ifnum\typofr=1\noi{\it Remarque.}\ \else\noi{\it Remark.}\ \fi}
  \def\rems{\ifnum\typofr=1\noi{\it Remarques.}\ \else\noi{\it Remarks.}\ \fi}

  \def\1{{\bf 1}}
  \def\|{\Vert}

  \def\leq{\leqslant}
  \def\geq{\geqslant}


  \def\li{\mathop{\rm li}\nolimits}

  \def\log{\mathop{\rm log}\nolimits}




  \def\pmb#1{\setbox0=\hbox{#1}%
  \kern-.025em\copy0\kern-\wd0\kern.05em\copy0\kern-\wd0\kern-.025em\raise .0433em\box0 }

  
  \skewchar\eighti='177 \skewchar\sixi='177
  \skewchar\eightsy='60 \skewchar\sixsy='60
  
  \def\eightpoint{%
  \textfont0=\eightrm\scriptfont0=\sixrm\scriptscriptfont0=\fiverm
  \def\rm{\fam0\eightrm}%
  \textfont1=\eighti\scriptfont1=\sixi
  \scriptscriptfont1=\fivei\def\oldstyle{\fam1\seveni}%
  \textfont2=\eightsy\scriptfont2=\sixsy\scriptscriptfont2=\fivesy
  \textfont3=\eightex\scriptfont3=\sixex
  \textfont\itfam=\eightit
  \def\it{\fam\itfam\eightit}%
  \textfont\slfam=\eightsl
  \def\sl{\fam\slfam\eightsl}%
  \textfont\bbfam=\eightbb \scriptfont\bbfam=\sixbb\scriptscriptfont\bbfam=\fivebb
  \def\bb{\fam\bbfam\eightbb}%
  \textfont\msafam=\eightmsa\scriptfont\msafam=\sixmsa
  \textfont\scalnfam=\eightscaln
  \def\scaln{\fam\scalnfam\eightscaln}
  \textfont\ttfam=\eighttt
  \def\tt{\fam\ttfam\eighttt}%
  \textfont\bffam=\eightbf\scriptfont\bffam=\sixbf\scriptscriptfont\bffam=\fivebf
  \def\bf{\fam\bffam\eightbf}%
  \textfont\pcapfam=\eightpcap
  \def\pcap{\fam\pcapfam\eightpcap}
  \abovedisplayskip=2pt plus2pt minus 2pt
  \belowdisplayskip=2pt plus1pt minus 2pt
  \abovedisplayshortskip= 1pt plus 2pt minus 1pt
  \belowdisplayshortskip= 1pt plus 2pt minus 1pt
  \smallskipamount=2pt plus 1pt minus 2pt
  \medskipamount=3pt plus 2pt minus 2pt
  \bigskipamount=7pt plus 3pt minus 3pt
  \setbox\strutbox=\hbox{\vrule height 5pt depth 2pt width 0pt}%
  \normalbaselineskip=9pt\normalbaselines\rm}

  \def\({\left(}
  \def\){\right)}
  
  \def\footnoterule{\kern -2pt\hrule width 7truecm\kern 2.4pt}
  
  \def\xnotedef#1{\definexref{#1}{\noexpand\number\footnotenumber}{Note}}%

  
  
  \def\ninepoint{%
  \textfont0=\ninerm\scriptfont0=\sixrm\scriptscriptfont0=\fiverm
  \def\rm{\fam0\ninerm}%
  \textfont1=\ninei\scriptfont1=\sixi
  \scriptscriptfont1=\fivei\def\oldstyle{\fam1\ninei}%
  \textfont2=\ninesy\scriptfont2=\sixsy\scriptscriptfont2=\fivesy
  \textfont3=\nineex\scriptfont3=\sixex
  \textfont\itfam=\nineit
  \def\it{\fam\itfam\nineit}%
  \textfont\slfam=\ninesl
  \def\sl{\fam\slfam\ninesl}%
  \textfont\bbfam=\ninebb\scriptfont\bbfam=\sixbb\scriptscriptfont\bbfam=\fivebb
  \def\bb{\fam\bbfam\ninebb}%
  \textfont\msafam=\ninemsa\scriptfont\msafam=\sixmsa\scriptscriptfont\msafam=\fivemsa
  \textfont\scalnfam=\ninescaln
  \def\scaln{\fam\scalnfam\ninescaln}
  \textfont\ttfam=\ninett
  \def\tt{\fam\ttfam\ninett}%
  \textfont\bffam=\ninebf\scriptfont\bffam=\sixbf\scriptscriptfont\bffam=\fivebf
  \def\bf{\fam\bffam\ninebf}%
  \abovedisplayskip=3pt plus2pt minus 2pt
  \belowdisplayskip=3pt plus1pt minus 2pt
  \abovedisplayshortskip= 2pt plus 2pt minus 1pt
  \belowdisplayshortskip= 2pt plus 2pt minus 1pt
  \smallskipamount=2pt plus 1pt minus 2pt
  \medskipamount=3pt plus 2pt minus 2pt
  \bigskipamount=7pt plus 3pt minus 3pt
  \setbox\strutbox=\hbox{\vrule height 5pt depth 2pt width 0pt}%
  \normalbaselineskip=10.5pt plus.3pt minus.3pt\normalbaselines\rm}

  \def\sevenpoint{%
  \textfont0=\sevenrm\scriptfont0=\sixrm\scriptscriptfont0=\fiverm
  \def\rm{\fam0\sevenrm}%
  \textfont1=\seveni\scriptfont1=\sixi
  \scriptscriptfont1=\fivei\def\oldstyle{\fam1\seveni}%
  \textfont2=\sevensy\scriptfont2=\sixsy\scriptscriptfont2=\fivesy
  \textfont3=\sevenex\scriptfont3=\fiveex
  \textfont\itfam=\sevenit
  \def\it{\fam\itfam\sevenit}%
  \textfont\slfam=\sevensl
  \def\sl{\fam\slfam\sevensl}%
  \textfont\bbfam=\sevenbb \scriptfont\bbfam=\sixbb\scriptscriptfont\bbfam=\fivebb
  \def\bb{\fam\bbfam\sevenbb}%
  \textfont\msafam=\sevenmsa\scriptfont\msafam=\sixmsa
  \textfont\scalnfam=\sevenscaln
  \def\scaln{\fam\scalnfam\sevenscaln}
  \textfont\bffam=\sevenbf\scriptfont\bffam=\sixbf\scriptscriptfont\bffam=\fivebf
  \def\bf{\fam\bffam\sevenbf}%
  \textfont\ttfam=\seventt
  \abovedisplayskip=2pt plus2pt minus 2pt
  \belowdisplayskip=2pt plus1pt minus 2pt
  \abovedisplayshortskip= 1pt plus 2pt minus 1pt
  \belowdisplayshortskip= 1pt plus 2pt minus 1pt
  \smallskipamount=2pt plus 1pt minus 2pt
  \medskipamount=3pt plus 2pt minus 2pt
  \bigskipamount=7pt plus 3pt minus 3pt
  \setbox\strutbox=\hbox{\vrule height 5pt depth 2pt width 0pt}%
  \normalbaselineskip=9pt\normalbaselines\rm}

 \def\onzepts{%
 \textfont0=\elevenrm\scriptfont0=\ninerm
 \textfont1=\elevenib\scriptfont1=\ninei}

\def\douzepts{%
  \textfont0=\twelverm\scriptfont0=\tenrm\def\rm{\fam0\twelverm}%
  \textfont1=\twelveib\scriptfont1=\teni%
  \textfont2=\twelvesy\scriptfont2=\tensy\scriptscriptfont2=\eightsy
  \textfont3=\twelveex
  \textfont\itfam=\twelveti
  \def\it{\fam\itfam\twelveti}%
  \textfont\bffam=\twelvebf\scriptfont\bffam=\tenbf\scriptscriptfont\bffam=\eightbf
  \def\bf{\fam\bffam\twelvebf}%
  \textfont\bbfam=\twelvebb \scriptfont\bbfam=\tenbb
  \def\bb{\fam\bbfam\twelvebb}%
  \textfont\msafam=\twelvemsa\scriptfont\msafam=\tenmsa
  \textfont\scalnfam=\twelvescaln
  \normalbaselineskip=15pt\normalbaselines\rm}

\def\quatorzepts{%
  \textfont0=\fourteenrm\scriptfont0=\twelverm\def\rm{\fam0\fourteenrm}%
  \textfont1=\fourteenib\scriptfont1=\twelveib%
  \textfont2=\fourteensy\scriptfont2=\twelvesy\scriptscriptfont2=\tensy
  \textfont3=\fourteenex
  \textfont\itfam=\fourteenti
  \def\it{\fam\itfam\fourteenti}%
  \textfont\bffam=\fourteenbf\scriptfont\bffam=\twelvebf\scriptscriptfont\bffam=\tenbf
  \def\bf{\fam\bffam\fourteenbf}%
  \textfont\bbfam=\fourteenbb \scriptfont\bbfam=\twelvebb
  \def\bb{\fam\bbfam\fourteenbb}%
  \textfont\msafam=\fourteenmsa\scriptfont\msafam=\twelvemsa
  \textfont\scalnfam=\twelvescaln
  \normalbaselineskip=18pt\normalbaselines\rm}


\def\AA{{\it Acta Arith.}}

\def\picture #1 by #2 (#3){\leavevmode\vbox to #2{
     \hrule width #1 height 0pt depth 0pt
      \vfill
       \special{picture #3}}}

\def\illustration #1 by #2 (#3) scaled #4{\dimen1=#2
  \divide\dimen1 by 1000
  \multiply\dimen1 by #4
  \vtop to \dimen1{\dimen1=#1
  \divide\dimen1 by 1000
  \multiply\dimen1 by #4
  \hsize=\dimen1\vss
  \noindent\special{illustration #3 scaled #4}}}

\ifx\optionkeymacros\undefined\else\endinput\fi

\catcode`\å=\active\defå{{\aa}}       
\catcode`\†=\active\def†{\int}        
\catcode`\ç=\active\defç{\c c}        
\catcode`\=\active\def{\partial}    
\catcode`\Ÿ=\active\defŸ{\oint}       
\catcode`\=\active\def{\triangle}   
\catcode`\¬=\active\def¬{\neg}        
\catcode`\µ=\active\defµ{\mu}         
\catcode`\ø=\active\defø{{\o}}        
\catcode`\¼=\active\def¼{\pi}         
\catcode`\¦=\active\def¦{{\oe}}       
\catcode`\ß=\active\defß{{\ss}}       
\catcode`\Ý=\active\defÝ{\dagger}     
\catcode`\ˆ=\active\defˆ{\sqrt}       
\catcode`\…=\active\def…{\Sigma}      
\catcode`\‰=\active\def‰{\approx}     
\catcode`\‡=\active\def‡{\Omega}      
\catcode`\£=\active\def£{{\it\$}}     
\catcode`\ƒ=\active\defƒ{\infty}      
\catcode`\§=\active\def§{{\S}}        
\catcode`\¶=\active\def¶{{\P}}        
\catcode`\€=\active\def€{\bullet}     
\catcode`\ª=\active\defª{\leavevmode\raise.585ex\hbox{\b a}}      
\catcode`\º=\active\defº{\leavevmode\raise.6ex\hbox{\b o}}        
\catcode`\'=\active\def'{\not=}       
\catcode`\¾=\active\def¾{\leq}        
\catcode`\"=\active\def"{\geq}        
\catcode`\÷=\active\def÷{\div}        
\catcode`\Š=\active\defŠ{{\dots}}     
\catcode`\æ=\active\defæ{{\ae}}       
\catcode`\«=\active\def«{\og}         
\catcode`\³=\active\def³{``}          
\catcode`\¡=\active\def¡{!`}          
\catcode`\¢=\active\def¢{\rlap/c}     
\catcode`\Œ=\active\defŒ{`}           
\catcode`\¹=\active\def¹{'}           


\catcode`\Å=\active\defÅ{{\AA}}       
\catcode`\Ç=\active\defÇ{\c C}        
\catcode`\Ø=\active\defØ{{\O}}        
\catcode`\½=\active\def½{\Pi}         
\catcode`\'=\active\def'{{\OE}}       
\catcode`\Æ=\active\defÆ{{\AE}}       
\catcode`\×=\active\def×{\diamond}    
\catcode`\°=\active\def°{\accent'27}  
\catcode`\²=\active\def²{''}          
\catcode`\±=\active\def±{\pm}         
\catcode`\»=\active\def»{\fg}         
\catcode`\¿=\active\def¿{?`}          
\catcode`\­=\active\def­{--}          
\catcode`\‹=\active\def‹{---}         


\catcode`\ä=\active\defä{\"a}        
\catcode`\ë=\active\defë{\"e}        
\catcode`\ï=\active\defï{\"{\i}}     
\catcode`\ö=\active\defö{\"o}        
\catcode`\ü=\active\defü{\"u}        
\catcode`\ÿ=\active\defÿ{\"y}        
\catcode`\Â=\active\defÂ{\^A}        
\catcode`\Ä=\active\defÄ{\"A}        
\catcode`\Ö=\active\defÖ{\"O}        
\catcode`\Ü=\active\defÜ{\"U}        
\catcode`\á=\active\defá{\'a}        
\catcode`\é=\active\defé{\'e}        
\catcode`\í=\active\defí{\'{\i}}     
\catcode`\ó=\active\defó{\'o}        
\catcode`\ú=\active\defú{\'u}        
\catcode`\É=\active\defÉ{\'E}        
\catcode`\Ê=\active\defÊ{\^E}        
\catcode`\à=\active\defà{\`a}        
\catcode`\è=\active\defè{\`e}        
\catcode`\ì=\active\defì{\`{\i}}     
\catcode`\ò=\active\defò{\`o}        
\catcode`\ù=\active\defù{\`u}        
\catcode`\À=\active\defÀ{\`A}        
\catcode`\ã=\active\defã{\~a}        
\catcode`\ñ=\active\defñ{\~n}        
\catcode`\õ=\active\defõ{\~o}        
\catcode`\Ã=\active\defÃ{\~A}        
\catcode`\Ñ=\active\defÑ{\~N}        
\catcode`\Õ=\active\defÕ{\~O}        
\catcode`\â=\active\defâ{\^a}        
\catcode`\ê=\active\defê{\^e}        
\catcode`\î=\active\defî{\^{\i}}     
\catcode`\ô=\active\defô{\^o}        
\catcode`\û=\active\defû{\^u}        

\let\optionkeymacros\null

\def\EE{{\bb E}}
\def\PP{{\bb P}}

\def\D{{\scal D}}

\def\paraunn#1{\paraunc{#1}\writetocentry{section}{#1}}
\def\paradeuxn#1{\paradeuxb{#1}\writetocentry{subsection}{#1}}

\newcount\constantes\constantes=0
\def\ddconst#1{\definexref{#1}{c_{\the\constantes}}{constantes}}
\def\dconst#1{\global\advance \constantes by 1\ddconst{#1}\refn{#1}}
\def\const#1{\refn{#1}}

\font\seventt=cmtt10 at 7pt
\font\fivett=cmtt10 at 5pt
\scriptfont\ttfam=\seventt
\scriptscriptfont\ttfam=\fivett
\font\tengp=cmmib10
\font\ninegp=cmmib10 at 9pt
\font\eightgp=cmmib8
\font\sevengp=cmmib7

\newfam\gpfam
\textfont\gpfam=\tengp\scriptfont\gpfam=\sevengp

\def\sump_#1^#2{\mathop{\sum_{#1}}^{#2}}
\def\ps{{\rm ps} }
\newcount\constantes\constantes=0
\def\ddconst#1{\definexref{#1}{c_{\the\constantes}}{constantes}}
\def\dconst#1{\global\advance \constantes by 1\ddconst{#1}\refn{#1}}
\def\const#1{\refn{#1}}

\def\dateheurejob{\vglue-10mm\rightline{(\jobname.tex,\ \the\day/\the\month/\the\year,\ \heure)}}

\optionparag=2

\dimart
\optionparag=1

\def\Psif#1#2{\Psi_{#1}(#2)}
\def\Psifbig#1#2{\Psi_{#1}\Big(#2\Big)}
\def\Psisf#1{\Psi^*(#1)}
\def\Psisfbig#1{\Psi^* \Big(#1\Big)}
\def\FF{{\cal F}}

\anote{}{{\it 2010 Mathematics Subject Classification.} Primary: 11N25, 11N37. Secondary: 60G42.
{\it Key words and phrases:} multiplicative functions, random multiplicative functions, friable integers, martingales, Doob's inequality.}

\def\titcour{Sommes friables de fonctions multiplicatives aléatoires}
\def\auteur{Joseph Basquin}
\bigskip\bigskip
\hautspages{\auteur}{\titcour}
\titrecentre{\titcour\anote{1}{Version au 18/09/2011.\qquad\qquad\qquad\qquad\qquad\qquad\qquad\qquad\qquad\qquad\qquad\qquad\qquad\qquad}}
\bigskip\bigskip
\centerline{\auteur}
\bigskip\bigskip
{\rightskip1cm\leftskip1cm
\ninepoint
\noi {\bf Abstract.}
We consider a sequence $\{f(p)\}_{p\ {\rm  prime}}$ of independent random variables taking values $\pm 1$ with probability $1/2$, and extend $f$ to a multiplicative arithmetic function defined on the squarefree integers.
We investigate upper bounds for $\Psi_f(x,y)$, the summatory function of $f$ on $y$-friable integers $\leq x$.

We obtain estimations of the type $\Psi_f(x,y) \ll \Psi(x,y)^{1/2+\epsilon}$, more precise formulas being given in suitable regions for $x,y$.

In the special case $y=x$, this leads to the estimate $M_f(x) = \sum_{n \leq x} f(n) \ll \sqrt{x}\, (\log \log x)^{2+\epsilon}$, which improves on previous bounds. \par}

\newcount\paras\paras=0
\newcount\sparas\sparas=0
\def\tocsectionentry#1#2{\init{\sparas}\avance\paras
      {\quad\bf\the\paras\quad }{\hskip-2mm #1\dotfill\hskip3mm\rm#2}\par}%
\def\tocsubsectionentry#1#2{\avance\sparas
      {\qquad\eightpoint\the\paras.\the\sparas}
{\hskip-2mm\eightpoint#1\dotfill\hskip3mm\rm#2}\par}%

\bigskip\bigskip
\noi{\qquad\quad \bf Sommaire}
\smallskip
{\eightpoint\leftskip.6cm\rightskip.8cm
\readtocfile}
  \bigskip\medskip

\paraunn{Introduction et description des résultats}\drefun{intro}
\bigskip
\paradeuxn{Fonctions multiplicatives aléatoires}

L'étude des sommes de fonctions multiplicatives
est un domaine central de la théorie analytique
des nombres.  Ces fonctions possèdent, dans la
plupart des cas intéressants, un comportement
statistique complexe. \par D'où l'idée, mise en
pratique par plusieurs auteurs, de considérer des
fonctions
multiplicatives {\it aléatoires}, dont les
variations, que l'on peut appréhender par les
outils de la théorie des probabilités,
fournissent un modèle pertinent de la situation
arithmétique. Selon les hypothèses effectuées sur
les variables aléatoires en cause, de tels
modèles peuvent être adaptés aux diverses
situations concrètes rencontrées en théorie des
nombres.

La voie privilégiée par Wintner dans \citer{W44}
consiste à considérer, pour chaque nombre premier
$p$, une
variable aléatoire de Bernoulli $f(p)$, prenant
les valeurs $1$ et $-1$ avec probabilité
$1/2$, et à étendre par {\it multiplicativité}
cette fonction $f$ à l'ensemble des entiers sans
facteur carré. Ceci donne lieu à la définition suivante.

\Propd{alea}{Notant $\P$ l'ensemble des nombres premiers, soient $\{f(p)\}_{p \in
\P}$ une suite de variables aléatoires de
Bernoulli indépendantes prenant
les valeurs $+1$ et $-1$ avec probabilité $1/2$,
et $f$ la fonction multiplicative définie pour $n\in \N^*$ 
par $$f(n):=\mu(n)^2 \prod_{p | n} f(p).$$
Sous ces hypothèses, nous dirons que $f$ est une {\rm fonction 
multiplicative aléatoire au sens de Wintner}.}
\goodbreak

Posons $$M_f(x):=\sum_{n \leq x} f(n) \qquad \qquad (x\geqslant 1).$$

Une mesure de l'indépendance statistique des facteurs premiers des entiers est obtenue en comparant la fonction sommatoire de telles fonctions multiplicatives aléatoires à celle d'un modèle probabiliste, comme par exemple celui d'une somme de variables aléatoires centrées indépendantes relevant du théorème central limite.
\par 
Citons un premier résultat de Wintner \citer{W44} allant dans cette direction.
Nous utilisons dorénavant la mention \ps pour qualifier une assertion 
aléatoire valide presque sûrement.

\par 
\proclaim Théorème A. Soit $f$ une fonction multiplicative aléatoire au sens de Wintner. Pour tout $\epsilon > 0$, nous
avons
$$ M_f(x) \ll x^{1/2 + \epsilon} \qquad (x \geq 1) \qquad \ps.$$

L'enjeu essentiel du problème est de comparer le cas
étudié ici au cas sans contrainte arithmétique, c'est-à-dire
$f(n)=\pm 1$ avec probabilité $1/2$, pour lequel
la loi du logarithme itéré fournit l'estimation\note{Ici et dans la suite, nous désignons par $\log_k$ la k-ième itérée de la fonction logarithme.}
$$M_f(x) \ll  \sqrt { x\log_2 x },$$ et ainsi d'élucider 
l'influence de  la condition arithmétique de
multiplicativité. \medskip

Dans un travail non publié, Erd\H os obtient l'existence d'une constante 
$\dconst{erd} > 0$ 
 pour laquelle on a la majoration
$$M_f(x) \ll \sqrt{x}\, (\log x)^{\const{erd}} \qquad \ps.$$

Hal\'asz \citer{H82} précise cette majoration et établit le résultat suivant.

\proclaim Théorème B. Soit $f$ une fonction multiplicative aléatoire au sens de Wintner.
Il existe une constante $\dconst{int1} > 0$ telle que l'on ait 
$$M_f(x) \ll \sqrt{x}\, \e^{\const{int1} \sqrt{(\log_2 x) \log_3 x}} \qquad (x \geq 16) \qquad \ps.$$

\medskip
Dans une version préliminaire \citer{LTW10} (améliorée depuis) d'un récent travail, Lau, Tenenbaum et Wu précisent encore cette majoration et obtiennent, pour une fonction multiplicative aléatoire au sens de Wintner $f$, 
$$M_f(x) \ll \sqrt{x}\, (\log_2 x)^{5/2+\epsilon} \qquad \ps.\eqdef{LTW}$$

\bigskip
\bigskip
\paradeuxn{Entiers friables}

Nous utilisons la notation $P^+(n)$ (resp. $P^-(n)$) pour le plus grand (resp. petit) facteur premier d'un entier $n$ avec la convention $P^+(1):=1$ (resp.
$P^-(1):=\infty$). Un entier dont le plus grand facteur premier ne dépasse
pas $y$  est dit $y$-friable. Nous désignons par
$S(x,y)$ l'ensemble des entiers $y$-friables
inférieurs ou égaux à $x$ et par $\Psi(x,y)$ son cardinal. Pour une fonction arithmétique $f$, nous posons de plus
    $$\Psif{f}{x,y} := \sum_{n \in S(x,y)} f(n).$$

Nous notons $$\Psisf{x,y} := \Psif{\mu^2}{x,y}$$ pour le nombre des entiers $y$-friables, sans facteur carré et n'excédant pas $x$.

L'étude de la friabilité et de son interaction avec les critères usuels
de description utilisés en théorie probabiliste des nombres
constitue une branche essentielle de l'arithmétique. Ce sujet fait l'objet d'une importante littérature depuis les
années 1950, en particulier ces vingt dernières
années.

\bigskip
\bigskip
\paradeuxn{Résultats}

Nous nous proposons ici d'étudier les sommes friables $\Psif{f}{x,y}$ d'une fonction
multiplicative aléatoire au sens de Wintner $f$.

Un raisonnement statistique laissant augurer que les termes
de la somme se comportent comme des variables aléatoires
indépendantes, il est raisonnable d'espérer que la somme
se comporte comme la racine carrée du nombre de
termes, éventuellement multipliée
par un facteur à faible croissance en $x$ et $y$.

Posons $\theta(y) := \sum_{p \leq y} \log p.$
Nous pouvons à présent énoncer notre résultat principal. 
\Propt{mainthm}{Soit $f$ une fonction multiplicative aléatoire au sens
de Wintner. Pour tout $\epsilon > 0$, il existe des constantes positives $\dconst{m17}, \dconst{expos}, \dconst{exposant}, \dconst{ctriv}$, telles que nous ayons presque sûrement,
$$\Psif{f}{x,y} 
\cases{\ll \sqrt{x}\, (\log_2 x)^{2+\epsilon} & si $x^{\const{m17}\,(\log_3 x)/ \log_2 x} \leq y \leq x$ \cr & et $x \geq 16,$\cr \cr
 \ll \sqrt{\Psisf{x,y}}\  \e^{\const{expos} u } (\log x)^{\const{exposant}}  & si $\const{ctriv} \log x \leq y < x^{\const{m17}\,(\log_3 x)/ \log_2 x}$,  \cr &et $x \geq 16$.}$$}
\bigskip                                                                                             
\rems (i) Dans le cas particulier $y=x$, il est à noter que le théorème précédent fournit une amélioration du résultat \eqref{LTW} :
$$M_f(x) \ll \sqrt{x}\, (\log_2 x)^{2+\epsilon} \qquad \ps.$$
\bigskip
(ii) Lorsque $\theta(y) \leq \log x$ et $x\geq 2$, nous avons
$$\Psif{f}{x,y} = \cases{0 & avec probabilité $1-2^{-\pi(y)}$, \cr 2^{\pi(y)} & avec probabilité $2^{-\pi(y)}$.}
$$
En effet, dans ce domaine, on a
$$\Psif{f}{x,y} = \sum_{P^+(n) \leq y} f(n) = \prod_{p \leq y} \Big( 1 + f(p) \Big) = 0 $$
si et seulement s'il existe $p \leq y$ tel que $f(p)=-1$, ce qui intervient avec probabilité $1-2^{-\pi(y)}$.

\bigskip
Le \ref{mainthm} implique le résultat suivant.
\Propc{cor1}{Soient $f$ une fonction multiplicative au sens de Wintner et $\epsilon > 0$. Il existe une constante une constante positive $C_\epsilon$ telle que, pour $x,y \geq 2$ et $y \geq C_\epsilon \log x$, nous ayons
$$\Psi_f(x,y) \ll \Psi^*(x,y)^{1/2+\epsilon}.$$}
\bigskip
\bigskip
\paraunn{Lemmes}
\bigskip
\paradeuxn{\'Evaluation d'une somme de fonction multiplicative}

Nous désignons par $\sum_n^{y,z}$ une sommation dont l'indice entier $n$ est soumis à la condition $p\, |\, n \Rightarrow y < p \leq z$.

Pour $y\geq 2$, nous posons de plus  $L(y) := \exp \{(\log y)^{3/5} / (\log_2 y)^{1/5} \}$.
\medskip
La majoration suivante, obtenue par Lau, Tenenbaum et Wu (\citer{LTW10}) pour certaines sommes de fonctions multiplicatives, prolonge le lemme 3 de \citer{H82}.
\Propl{lemsommecourte}{
Soient $\gamma, \delta > 0$, $y \geq 2$, $\kappa \geq 1$ et $a, b  \in \N^*$ tels que $b > a$.
Il existe des constantes $\dconst{csomcourte}, \dconst{calpha}, \dconst{clem36} > 0$ telles que nous ayons,
uniformément pour $\gamma \geq {\const{calpha}}/{L(y)^{\const{clem36}}}$ et $y \geq \{{b}/(b-a)\}^{1+\delta}$,

$$\sump_{a < r \leq b}^{y,y^{1+\gamma}} \mu(r)^2 \kappa^{\omega(r)} \leq \const{csomcourte} 
\frac{b-a}{\log y} \kappa\, \e^{2 \gamma\, \kappa}.$$
}

\nid En suivant la méthode employée par Hal\'asz dans \citer{H82} (lemme 3 (ii)), il vient
$$A :=\sump_{a < r \leq b}^{y,y^{1+\gamma}} \mu(r)^2 \kappa^{\omega(r)}
\leq \frac{1}{\log a} \sump_{m \leq b /y}^{y,y^{1+\gamma}} \mu(m)^2 \kappa^{\omega(m)+1} \sum_{a/m < p \leq b/m} \log p.
$$
Notons que l'on peut supposer $b > y$, car dans le cas contraire, la somme à estimer est nulle. 
Observons également que les hypothèses $y \geq \{{b}/(b-a)\}^{1+\delta}$ et $m \leq b/y$
impliquent $$\frac{b-a}m \geq \frac b {m\, y^{1/(1+\delta)}} \geq \Big(\frac b m\Big)^{1-1/(1+\delta)} > 1. \eqdef{sup1}$$
Le théorème de Brun-Titchmarsh (voir \citer{T08}, théorème I.4.16) nous permet alors, pour $m \leq b/y$, de majorer la somme en $p$ de la façon suivante
$$
\sum_{a/m < p \leq b/m} \log p \leq \log (b/m) \sum_{a/m < p \leq b/m} 1 \leq \dconst{bt} \frac{(b-a) \log (b/m)}{m \log ((b-a)/m)}\cdot
$$

La relation \eqref{sup1} implique par ailleurs que $$\eqalign{A &\leq \const{bt} 
\frac{(b-a) \kappa \log (b/m)}{(\log a) \log ((b-a)/m)}
 \sump_{m \leq b /y}^{y,y^{1+\gamma}} \frac{\mu(m)^2 \kappa^{\omega(m)}}{m}
\cr &\leq \const{bt} \frac{(b-a) \kappa}{\log a}
\prod_{y < p \leq y^{1+\gamma}} \Big( 1 + \frac \kappa p \Big)
.}$$

On conclut en faisant usage de l'estimation suivante, découlant par exemple du lemme 3.6 de \citer{BT05b},
$$\sum_{y < p \leq y^{1+\gamma}} \frac{1}{p} = \log(1+\gamma) + O\Big(\frac{1}{L(y)^{\const{clem36}}}\Big) \leq 2 \gamma,$$
dès lors que $\gamma \geq \const{calpha} / L(y)^{\const{clem36}}$.
\qed

\bigskip
\bigskip

\paradeuxn{Espérance de fonctions multiplicatives aléatoires}
Le lemme suivant est dû Bonami \citer{B70}. Hal\'asz (\citer{H82}, lemme 2) en a donné une nouvelle démonstration.  La formulation donnée ici est adaptée au cas des fonctions multiplicatives aléatoires.
\Propl{bonami2}{
Soient $f$ une fonction multiplicative aléatoire au sens de Wintner, $\ell \geq 2$ un entier pair, et $\{a_n\}_{n\geq 1}$ une suite de nombres complexes. Nous avons 
$$\Big| \EE\, \Big\{ \Big( \sum_{n \geq 1} a_n f(n) \Big)^\ell \Big\} \Big| \leq \Big(\sum_{n \geq 1} \mu(n)^2 |a_n|^2 (\ell-1)^{\Omega(n)}\Big)^{\ell /2}.$$}

Cette majoration s'avère cruciale pour évaluer les moments d'ordre pair de fonctions multiplicatives aléatoires.
\medskip
\rem Pour $\ell=2$ et $\{a_n\}_{n \geq 1}$ égale à la fonction indicatrice des entiers $y$-friables inférieurs à $x$, la majoration fournie par le lemme précédent est en fait une égalité. Nous avons en effet
$$\EE \ ( \Psif{f}{x,y}^2 ) = \Psisf{x,y} \qquad \qquad (x,y \geq 2). \eqdef{momentordre2}$$

Cela résulte du développement du carré sous la forme
$$\EE \ ( \Psif{f}{x,y}^2 ) = \EE \ \Big\{ \sum_{m,n \in S(x,y)} f(m) f(n) \Big\}
=\sum_{m,n \in S(x,y)} \EE \  f(m) f(n).$$
Il suffit ensuite d'observer que, si $m$ et $n$ sont sans facteur carré et $m \neq n$, il existe un nombre premier $p$ tel
que $p \parallel m\, n$.\note{Ici $p^\nu \parallel a$ signifie : $p^\nu |a$ et $p^{\nu+1} \nmid a$.} Il en découle $$\EE \  (f(m) f(n)) = \EE\ (f(p)) \, \EE\ (f(mn /p)) = 0.$$ 

\bigskip
\bigskip
\paradeuxn{Entiers friables sans facteur carré}

Le résultat suivant, obtenu par La Bretèche et Tenenbaum \citer{BT02}, 
fournit une description de la répartition des entiers friables sans facteur carré dans les petits intervalles.
\Propl{bretechetenen}{Soit $\kappa \geq 1$. Nous avons, uniformément sous les conditions
$x \geq y \geq 2$ et $\max(1,x y^{-\kappa})\leq~z\leq x$, 
$$\Psisf{x+z,y}-\Psisf{x,y} \ll \frac{z}{x} \Psisf{x,y}.$$}

Posons à présent $\phi(s,y) := \sum_{p\leq y} \log (1+1/p^s)$ $(s \geq 0)$ et désignons par
$\alpha = \alpha(x,y)$ l'unique solution positive de l'équation $$-\phi'(\alpha,y) =
\sum_{p\leq y} \frac{\log p}{1+ p^\alpha} = \log x.
\eqdef{def-alpha}$$ 
Par ailleurs, nous désignons par $\beta := \beta(x,y)$ l'unique solution positive de l'équation
$$\sum_{p\leq y} \frac{\log p}{p^\beta-1} = \log x.$$ 

Pour $t>0$, $t \neq 1$, nous définissons $\xi(t)$ comme l'unique solution réelle non nulle de l'équation $$\e^\xi = 1 + t \xi,$$ et posons $\xi(1) = 0$.

Rappelons que pour $\epsilon > 0$, $x \geq x_0(\epsilon)$ et $(\log x)^{1+\epsilon} < y \leq x$, nous avons les estimations (voir (7.8) de \citer{HT86} et (2.17) de \citer{BT02})
$$\eqalign{\alpha(x,y) \cr \beta(x,y)} \Big\} = 1- \frac{\xi(u)}{\log y} + O\Big(\frac{1}{L_\epsilon(y)}+\frac{1}{u (\log y)^2}\Big), \eqdef{estalphabeta}$$
où nous avons posé
$$L_\epsilon(y) := \exp\big\{(\log y)^{3/5-\epsilon}\big\}.$$
\goodbreak
Le résultat suivant est un analogue du théorème 2.4 de \citer{BT05b},
énoncé ici dans le cas des entiers friables sans facteur carré.
\Propl{semi-asymp-sf}{
Il existe des constantes absolues positives $b_1, b_2$, et une fonction $b = b(x,y;d) \in [b_1, b_2]$ telles que, sous les conditions
$$x \geq y \geq 2, \qquad \vartheta(y) > 2 \log x, \qquad 1 \leq d \leq x,$$
nous ayons uniformément
$$\Psisfbig{\frac{x}{d},y} = \frac{\Psisf{x,y}}{d^{\alpha}}
\e^{- b t^2 / u} \Big\{ 1 + O\Big(\frac{t}{\sqrt{u}}+\frac{1}{u_y}\Big) \Big\}, \eqdef{sasfc}$$
où $t := (\log d) / \log y$, $u_y := u+(\log y)/\log (u+2)$,
et $\alpha = \alpha(x,y)$ est défini par \eqref{def-alpha}.}
\goodbreak 
\nid 
Notons que nous pouvons supposer $x$ et $y$ assez grands.
En effet, lorsque $x$ est borné, le résultat est acquis sous réserve que les constantes implicites des termes d'erreurs soient choisies suffisamment grandes. Si $y$ est borné, la relation $\vartheta(y) > 2 \log x$ implique que $x$ est borné à son tour, d'où le résultat.
Nous supposons donc dans toute la suite que $x$ et $y$ sont suffisamment grands.

Plaçons-nous tout d'abord  dans le domaine
$$ y > (\log x)^3, \qquad  x \geq x_0. \eqdef{dom1}$$
Nous pouvons faire appel aux résultats de \citer{BT05b} concernant le comportement local de $\Psi(x,y)$.
Au vu de l'estimation \eqref{estalphabeta}, nous avons 
$\beta \geq 3/5$, dès que $y$ est assez grand.  
Les quantités $\zeta(2 \beta,y)$ et $\prod_{p \leq y} (1+p^{-2\beta})$ étant bornées, le corollaire 2.6 de \citer{BT05b} nous permet alors d'écrire
$$\Psisf{x, y} = \frac{\Psi (x, y )}{\zeta(2 \beta,y)}\Big\{ 1 + O\Big(\frac{1}{u_y}\Big)\Big\},$$ 
où nous avons posé, pour $y \geq 2$,
$$\zeta(s,y) := \sum_{P^+(n)\leq y} \frac{1}{n^s} = \prod_{p \leq y} \Big(1 - \frac 1 {p^s}\Big)^{-1}\cdot$$

Nous avons d'une part
$$\zeta(2 \beta, y) = \zeta(2 \beta) + O(1/u_y),$$
et d'autre part,
$$\sum_{p\leq y}\frac{(\beta-\alpha) (\log p)^2}{p^\beta -1} \ll
\sum_{p\leq y}\frac{(p^\beta-p^\alpha) \log p}{(p^\alpha+1)(p^\beta -1)} 
= \sum_{p\leq y}\frac{2 \log p}{(p^\alpha+1)(p^\beta -1)} \ll 1,$$
ce qui implique encore
$$\beta-\alpha \ll 1/u(\log y)^2.$$

Le théorème 2.4 de \citer{BT05b} implique alors
$$\Psisfbig {\frac x d, y } = \frac{\Psisf{x,y}}{d^{\alpha}} \exp\Big\{\frac{-b\, t^2} u + (\alpha - \beta) \log d \Big\} 
\Big\{ 1 + O\Big(\frac{t}{{u}} + \frac{1}{u_y} \Big) \Big\},$$
sous la condition \eqref{dom1}, ce qui fournit bien la formule \eqref{sasfc} dans ce domaine.
De plus, on peut y remplacer $t/\sqrt{u}$ par $t/u$ dans le terme d'erreur.
\medskip
Il reste à examiner le cas 
$$x \geq x_0, \qquad 2 \log x < \vartheta(y) \leq (\log x)^3,$$
Nous pouvons nous limiter à prouver le résultat lorsque $t \leq u-1$, i.e. $d \leq x /y$.  En effet, si nous supposons le résultat établi dans ce sous-domaine, nous avons, lorsque $x / y < d \leq x$,
$$\eqalign{\Psisfbig{\frac{x}{d},y } \asymp \frac{x}{d} = \frac{x}{d\, y} y &\ll\frac{x}{d\, y} \frac{\Psisf{x,y}}{(x/y)^{\alpha}} \sqrt{u}\,  \e^{-b\, (u-1)^2 / u}
\cr &\ll \frac{\Psisf{x,y}}{d^{\alpha}} \e^{-b u / 2}
\ll \frac{\Psisf{x,y}}{d^{\alpha}} \e^{-b t^2 / 2 u},
}$$
car $\alpha < 1$ d'après le lemme 2.8 de \citer{BT02} sous la forme
$$\alpha = \Big\{ 1 + O\Big(\frac{1}{\log y}\Big) \Big\} \frac{\log (1+z)}{\log y},$$
où $z:=z(x,y) > 0$ est implicitement défini par la relation $\theta(y) = (2+ z) \log y$.

Notons que dans le domaine considéré, nous avons $u \gg u_y$ et
le théorème 2.1 de \citer{BT02} fournit, pour $1\leq d \leq x/y$,
$$\Psisfbig{\frac{x}{d},y} = \e^{h(u-t)} \Big\{1 + O \Big(\frac{t+1}{u}\Big)\Big\},$$
où nous avons posé 
$$h(v) := \alpha_v v \log y + \phi(\alpha,y) + g(\alpha \sqrt{\sigma_2}), \qquad \alpha_v := \alpha(y^v,y) \qquad (v \geq 1),$$
$$g(z) := \log \Big\{ \frac {\e^{z^2/2}} { \sqrt{2 \pi}} \int_z^\infty \e^{-t^2 /2} \d t \Big\}, \qquad \sigma_2 := \frac{\d^2 \phi}{\d s^2}(\alpha,y).$$

Il est établi en (2.35) de \citer{BT02} que l'on a, $$h'(v) = \alpha_v \log y +O\Big( \frac{1}{\sqrt{v} + v \log (1+z)} + \frac{1}{v}  \Big)
\qquad (v\geq 1). \eqdef{hprimv}$$ Il vient alors, pour $t \leq u-1$,
$$h(u-t)=h(u)-\alpha_u \, t \log y - I + R,$$
avec $R \ll \int_{u-t}^u \d v / \sqrt{v} \ll t/\sqrt{u}$ et 
$$I := (\log y) \int_{u-t}^u \{ \alpha_v - \alpha_u \} \d v
\asymp \int_{u-t}^u \int_v^u \frac{\d w \d v}{w} = - (u-t) \log \Big(\frac{u}{u-t}\Big) + t,$$
grâce à l'esimation $\alpha_v' \asymp -1 / (v \log y)$ $(v\geq 1)$ obtenue dans \citer{BT02}. Il vient ainsi $$I \asymp{t^2}/{u},$$ ce qui fournit le résultat
car $-I+R \asymp -t^2 / u + O(1)$ lorsque $t>\sqrt{u}$.

Notons cependant que l'on peut estimer $R$ plus précisément lorsque $t \gg \sqrt{u} \log u$ :
$$\eqalign{R &\ll \int_1^u \frac{\log (2 u / v ) \d v}{\sqrt{v} + v \log (2 u / v -1)}  \asymp \int_1^{2u-1} \frac{u \log (w+1) \d w}{w \sqrt{w u } + u w \log w}\cr &\ll
\int_1^{1+1/\sqrt{u}} \sqrt{u} + \int_{1+1/\sqrt{u}}^{2u-1} \frac{\log (w+1) \d w}{w \log w} \ll \log u.}$$
\qed

\rem La démonstration précédente montre également que l'on peut remplacer
$t/\sqrt{u}$ par $t/u$ dans le terme d'erreur de \eqref{sasfc}, dans tout domaine du type 
$$x \geq y \geq 2, \qquad \theta(y) \geq (2+\epsilon) \log x,$$
puisque dans ce cas, \eqref{hprimv} fournit $R \ll \int_{u-t}^u \d v / v \ll \log (u / (u-t)) \ll t/u$
pour $t\leq u/2$.
\bigskip

\goodbreak
Le lemme suivant s'avère utile pour évaluer certaines sommes portant sur les nombres premiers.
\Propl{sum1surp}{
Nous avons, uniformément pour $x\geq y \geq 2$ et  $2 \log x < \vartheta(y)$,
$$\sum_{p \leq y} \frac1{p^\alpha} = \log_2 y + \frac{u w}{w-1} \Big(1 + 
O\Big(\frac{1}{\log y} + \frac{1}{\log 2u} \Big) \Big), \eqdef{eqsum1surp}$$
où $w:=\vartheta(y) / \log x$ et $\alpha = \alpha(x,y)$ est défini par \eqref{def-alpha}.
}

\nid
Le lemme 3.6 de \citer{BT05b} fournit, pour $y \geq 2$,
$$\sum_{p \leq y} \frac1{p^\alpha} = \log_2 y + \Big\{1 + O\Big(\frac{1-\alpha}{L(y)^{\dconst{lemme36}}}\Big)\Big\}
\int_1^{v_\alpha} t \xi'(t) \d t + O(1),$$
où nous avons posé $$v_\alpha := \frac{y^{1-\alpha} -1 }{(1-\alpha)\log y}\cdot$$

Il découle du lemme 2.7 de \citer{BT02} que
$$\log x = \frac{y^{1-\alpha}-1}{(1+y^{-\alpha}) (1-\alpha)} \Big\{ 1 + O\Big(\frac{1}{\log y}\Big) \Big\} + O(1)$$
et des formules (2.17) et (2.18) de \citer{BT02} que
$$1+y^{-\alpha} = \frac{w}{w-1} \Big\{ 1 + O\Big(\frac{1}{\log y}\Big) \Big\} \qquad\big(x \geq y \geq 2,\  2 \log x < \vartheta(y)\big).$$
Nous obtenons alors l'estimation
$$v_\alpha = \frac{u w}{w-1} \Big\{ 1 + O\Big(\frac{1}{\log y}\Big) \Big\}$$
et la relation
$$\xi'(t) = \frac{1}{t} \Big\{1 + O \Big( \frac{1}{\log 2t}\Big) \Big\} \qquad (t \geq 1)$$
permet de conclure.
\qed

\bigskip
\bigskip
\paraunn{Comportement local de $\Psif{f}{x,y}$}
\bigskip
\paradeuxn{Martingale $y \mapsto \Psi_f(x,y)$ et inégalité maximale}

Soient $(\Omega, \A, \PP)$ un espace probabilisé, $\{X_n\}_{n \geq 0}$ une suite de variables aléatoires et $\FF = \{\FF_n\}_{n \geq 0}$ une filtration, i.e. une famille croissante de sous-tribus de $\A$.

Nous rappelons que $\{X_n\}_{n \geq 0}$ est dite une $\FF$-martingale
(resp. une $\FF$-sous-martingale) si, pour tout entier $n \geq 0$, $X_n$ est $\FF_n$-mesurable,
$\EE\ ( |X_n|) < \infty$ et pour tout entier $n \geq 1$,
$$\EE\ ( X_n | \FF_{n-1} ) = X_{n-1} \qquad \ps \qquad ({\rm resp. } \geq X_{n-1}).$$

Rappelons l'inégalité de Doob pour les sous-martingales positives (voir par exemple \citer{W91}, §14.6, {\it Theorem}).
\Propl{doobinequ}{Soit $\{X_n\}_{n \geq 0}$ une $\FF$-sous-martingale positive. Pour tous $t>0$ et $n \geq 0$, nous avons $$\PP\Big(\max_{k \leq n} X_k \geq t\Big) \leq \frac{\EE\ (X_n)}{t} \cdot$$
}

Rappelons également l'assertion suivante (voir \citer{W91}, §14.6, {\it Lemma (b)}).
\Propl{jens}{Si $\{X_n\}_{n \geq 0}$ est une $\FF$-martingale, $\ell$ est un nombre entier $\geq 1$, et
$\EE \ | X_n^\ell | < \infty$ $(n \geq 0)$, alors $\{X_n^\ell\}_{n \geq 0}$ est une $\FF$-sous-martingale.}

Considérant une fonction aléatoire au sens de Wintner $f$, nous désignons à présent par $\FF_y := \sigma \{f(p) : p \leq y\}$ la 
tribu engendrée par les variables aléatoires $\{f(p)\}_{p \leq y}$, et par $\FF$ la filtration $\{\FF_y\}_{y \geq 2}$.

\Propp{marting}{Soit $f$ une fonction multiplicative aléatoire au sens de Wintner.
Pour tout $x \geq 1$, la suite de variables aléatoires $\{\Psi_f(x,y)\}_{y \geq 2}$ est une $\FF$-martingale.
}
\nid  Notons d'emblée que $\Psi_f(x,y)$ est $\FF_y$-mesurable, pour $y \geq 2$. Par ailleurs, nous avons pour tout $ y \geq 2$ premier,
$$\eqalign{\EE\, \big( \Psi_f(x,y) \ | \ \FF_{y-1} \big) &= 
\EE\, \big( \Psi_f(x,y-1) \ | \ \FF_{y-1} \big) \cr &\qquad \qquad + \EE\, ( f(y) \Psi_f(x/y,y-1)
\ | \ \FF_{y-1} ) \cr
&= \Psi_f(x,y-1) + \Psi_f(x/y,y-1) \EE\, ( f(y) )  \cr &=  \Psi_f(x,y-1),}$$
car $\EE\, (f(y)) = 0$. 

Pour $y$ non premier, la relation précédente
est bien encore vérifiée car $\Psi_f(x,y)=\Psi_f(x,y-1)$.

 Nous avons montré que $\{\Psi_f(x,y)\}_{y \geq 2}$ est une $\FF$-martingale.
\qed
\bigskip
Le résultat suivant est un cas particulier du théorème 2.1 de \citer{TL06}.
\Propl{lemtl}{Soient $\{u_n\}_{n \geq 0}$ une suite de réels positifs, et
$\{ X_n \}_{n \geq 0}$ une suite de variables aléatoires.
Supposons qu'il existe une constante $\dconst{tl06c1} >0$ telle que pour tout $n \geq 0$ et $t>0$, nous ayons
$$\PP\Big(\sup_{k \leq n} |X_k| \geq t\Big) \leq \frac{\const{tl06c1}}{t} \sum_{k \leq n} u_k.$$
Alors il existe une constante $\dconst{tl06c2} >0$, telle que pour toute suite croissante de réels positifs $\{v_n\}_{n \geq 0}$, nous ayons, pour tout $n \geq 0$ et $t>0$,
$$\PP\Big(\sup_{k \leq n} \frac{|X_k|}{v_k} \geq t\Big) \leq \frac{\const{tl06c2}}{t} \sum_{k \leq n} \frac{u_k}{v_k}\cdot$$}
\bigskip

Nous pouvons en déduire l'estimation suivante.
\Propl{lem_mart}{
Soit $$N_j :=\int_{1}^{X} \frac{\Psif{f}{v,z_j}^2}{v^2} \d v. \eqdef{defNj}$$
Nous avons, pour tout $t >0$, 
$$\PP\Big(\sup_{j \leq J} \frac{N_j}{\log z_j} > t\Big) \ll \frac{\log h}{t},$$
avec $X \geq 2$, $z_0 \geq 2$, $h := \log X / \log z_0$, $\gamma >0$, $z_j := {z_0}^{(1+\gamma)^j}$ ($j \geq 1$), $J := \min \{j \geq 1 : z_j \geq X \}$.}

\nid
La relation \eqref{momentordre2} et l'estimation du théorème III.5.1 de \citer{T08}
$$\Psi(x,y) \ll x\, \e^{-u/2} \qquad (x \geq y \geq 2) \eqdef{typerankin}$$
impliquent pour $v \geq 2$,
$$\EE \ \big( \Psif{f}{v,z_j}^2 \big)= \Psisf{v,z_j} \leq \Psi(v,z_j) \ll v\, \e^{-(\log v) / (2 \log z_j)}.$$

Par conséquent, pour tout $j\leq J$, nous obtenons
$$\eqalign{\EE\ N_j &= 
\int_{1}^{X} \frac{\EE \ \Psif{f}{v,z_j}^2}{v^2} 
\d v \ll \int_1^X \frac{\e^{- (\log v) / (2 \log z_j) }}{v} \d v \ll \log z_j.}$$

Le \ref{marting} nous assure que $\{\Psi_f(v,z_j)\}_{j \leq J}$ est une $\{\FF_{z_j}\}_{j\leq J}$-martingale et donc que 
la suite $\{N_j / \log z_0\}_{j \leq J}$ est une $\{\FF_{z_j}\}_{j\leq J}$-sous-martingale. Nous avons ainsi, par l'inégalité de Doob pour les sous-martingales positives, pour tout $t >0$ et tout $j \leq J$, 
$$\PP\Big(\sup_{k \leq j} \frac{N_k}{\log z_0} > t\Big) \leq \frac{\EE\, N_j}{t \log z_0} \leq \dconst{lemc1} \frac{(1+\gamma)^j}{t} \leq \const{lemc1} \frac{\sum_{k \leq j} u_k}{t},$$
où nous avons posé $u_0 := 1$ et $u_k := \gamma (1+\gamma)^{k}$ ($k\geq 1$).
En effet, nous avons
$\sum_{k \leq j} u_k = (1+\gamma)^{j+1} - \gamma \geq (1+\gamma)^j$ pour tout $j\leq J$.

Appliquons le \ref{lemtl} avec la suite croissante $v_k := (1+\gamma)^k$ ($k\geq 0$), et il vient
$$\PP\Big(\sup_{j \leq J} \frac{N_j}{\log z_j} > t\Big) 
= \PP\Big(\sup_{j \leq J} \frac{N_j}{v_j (\log z_0)} > t\Big)
\leq \dconst{lemc2} \frac{\sum_{j \leq J} u_j / v_j}{t} \ll \frac{1+\gamma J}{t},$$
d'où le résultat, puisque nous avons $ J \asymp {(\log h)}/{\gamma}$.
\qed

\medskip

\rem La majoration triviale donnant
$$\PP\Big(\sup_{j \leq J} \frac{N_j}{\log z_j} > t\Big) \leq \sum_{j \leq J}
\PP\Big(\frac{N_j}{\log z_j} > t\Big) \ll \frac{J}{t} \asymp \frac{\log h}{\gamma\, t},$$
le précédent lemme permet donc de gagner un facteur $1/\gamma$.
\bigskip
\bigskip

\paradeuxn{\'Etude des variations de $\Psif{f}{x,y}$}

Nous ferons usage des notations additionnelles suivantes.
\'Etant donné une suite strictement croissante d'entiers $\{x_k\}_{k\geq 1}$, 
nous posons, pour tout $j, k\geq 1$, 
$$y_{k,j} := {x_k}^{1/j},$$
et désignons l'ensemble de tous les intervalles dyadiques inclus
dans $]x_{k-1},x_k]$ par 
$$\D_k := \Big\{
x_{k-1} + \big] (s-1) 2^m , s 2^m \big] : s \geq 1,\ m\geq 0,\
s 2^m \leq x_k - x_{k-1} \Big\}. \eqdef{defDk}$$ \bigskip

En vue d'obtenir une majoration des quantités 
$|\Psif{f}{x,y} - \Psif{f}{x_k,y}|$, pour $x_{k-1} < x \leq x_k$ et $y \leq x_k$,
nous établissons le lemme suivant, analogue du lemme 1 de \citer{H82}. La différence réside ici dans la prise en compte d'un paramètre de friabilité.
Nous faisons ici pleinement usage de la martingale $y \mapsto \Psi_f(x,y)$, afin d'obtenir des majorations uniformes sur de larges plages du paramètre $y$.

\Propl{lem1}{Soient $f$ une fonction multiplicative aléatoire au sens de Wintner.
Il existe des constantes $\dconst{xk}\in]0,1[$ et $\dconst{expoz} > 0$ telles que, posant\note{Nous désignons par $[x]$ la partie entière d'un nombre réel $x$.} 
$$x_k := [ \exp\{k^{\const{xk}}\} ], \qquad J_k := \max \Big\{ j \geq 1 :  2 \log x_k  < \vartheta (y_{k,j}) \Big\} \qquad (k \geq 1),$$
nous ayons $$\max_{\di{]a,b] \in \D_k}{y \leq y_{k,j}}} \big|\Psif{f}{b,y} - \Psif{f}{a,y} \big| \ll \frac{\sqrt{\Psisf{x_k,y_{k,j}}}\  \e^{\const{expoz}\, j} }{\log x_k} \quad (k \geq 1,\  j \leq J_k) \quad \ps.$$}

\nid Utilisant le \ref{bonami2}, nous avons pour tous entiers $b > a \geq 2$ et $y \geq 2$,
$$\EE\  \big\{\big(\Psif{f}{b,y} - \Psif{f}{a,y} \big)^4 \big\}
\leq \Big\{ \sum_{\di{a<n\leq
b}{P^+(n) \leq y}} \mu(n)^2 3^{\omega(n)}
\Big\}^2 =: E_{a,b}^{y}.$$

\'Etant donnés $k\geq 1$, $1 \leq j \leq J_k$ et $R_{k,j} \geq 1$, un paramètre qui sera fixé ultérieurement, nous considérons l'événement
$$A_{k,j} := \Big\{
\max_{\di{]a,b] \in \D_k}{y \leq y_{k,j}}}
\big|\Psif{f}{b,y} - \Psif{f}{a,y} \big| \geq R_{k,j} \Big\}.$$

 Une application de l'inégalité de Doob à la sous-martingale
positive (cf. \ref{marting} et \ref{jens}) 
$\{(\Psif{f}{b,y}-\Psif{f}{a,y})^4\}_{y\geq 2}$
 fournit
$$\PP\Big(\max_{y \leq y_{k,j}} \big|\Psif{f}{b,y} - \Psif{f}{a,y} \big| \geq R_{k,j} \Big) \leq
\frac{1} {R_{k,j}^4} E_{a,b}^{y_{k,j}}.$$
Il vient alors, par sommation,
$$\PP(A_{k,j}) \leq \frac{1}{R_{k,j}^4
} \sum_{]a,b] \in \D_k} E_{a,b}^{y_{k,j}}.$$

Utilisons à présent le fait que $E_{a,b}^y + 
E_{b,c}^y \leq E_{a,c}^y$ pour tous entiers $c > b > a \geq 2$ et $y \geq 2$, pour obtenir la majoration
$$\eqalign{\sum_{]a,b] \in \D_k} E_{a,b}^{y_{k,j}} &= 
\sum_{2^m \leq x_k - x_{k-1}} \sum_{s \leq 
(x_k - x_{k-1}) / 2^m} E_{x_{k-1} + (s-1) 
2^m, x_{k-1} + s 2^m}^{y_{k,j}} \cr
&\leq \sum_{2^m \leq x_k - x_{k-1}} E_{x_{k-1},x_k}^{y_{k,j}}
 \ll (\log x_k) E_{x_{k-1}, x_k}^{y_{k,j}}.}$$

Ainsi
$$ \PP (A_{k,j}) \ll \frac{\log x_k}{R_{k,j}^4} E_{x_{k-1},x_k}^{y_{k,j}}. \eqdef{probexceptionnelle}$$

Nous avons, par l'inégalité de Hölder,
$$\eqalign{E_{x_{k-1},x_k}^{y_{k,j}} &= \Big\{\sum_{\di{x_{k-1} < n \leq x_k}{P^+(n) \leq y_{k,j}}} \mu(n)^2 3^{\omega(n)}\Big\}^2 \cr &\leq \Big\{\Psisf{x_k,y_{k,j}} - \Psisf{x_{k-1},y_{k,j}}\Big\}^{4/3}  
 \Big\{\!\!\!\!\sum_{n \in S(x_k,y_{k,j})} \mu(n)^2 3^{3 \omega(n)}\Big\}^{2/3}\!\!.}$$

Nous avons, pour $x,y \geq 2$ et $\theta(y) > 2 \log x$, 
$$\eqalign{\sum_{n \in S(x,y)} \mu(n)^2 3^{3 \omega(n)} &\leq\! \sum_{d \in S(x,y)} \mu(d)^2 26^{\omega(d)} \Psisfbig{\frac x d,y} \leq \Psisf{x,y} \!\! \sum_{P^+(d) \leq y} \!\! \frac{\mu(d)^2 26^{\omega(d)}}{d^\alpha} \cr &\leq \Psisf{x,y} \prod_{p \leq y} \Big(1+\frac{26}{p^\alpha}\Big) \leq \Psisf{x,y} (\log y)^{26} \exp \{\const{expoz} u \},}$$
où nous avons fait appel à \eqref{sasfc} et \eqref{eqsum1surp}. En appliquant de plus le \ref{bretechetenen}, valable sous réserve de l'existence d'une constante positive $K$ telle que 
$$\frac{x_k - x_{k-1}}{x_{k-1}} \geq \frac1{x_k^{K/j}}, \eqdef{petitint}$$ nous obtenons enfin, pour tout $j \leq J_k$,
$$E_{x_{k-1},x_k}^{y_{k,j}} \ll \Big(\frac{x_k - x_{k-1}}{x_{k-1}}\Big)^{4/3} \Psisf{x_k,y_{k,j}}^{2}(\log x_k)^{52/3} \e^{\const{expoz} j},$$

Considérons l'événement
$$A_k :=\bigcup_{j \leq J_k}
A_{k,j}$$ et posons $$R_{k,j} := \frac{\sqrt{\Psisf{x_k, y_{k,j}}} \ \e^{\const{expoz}\, j}} {\log x_k} \qquad (k \geq 1, \, 1 \leq j \leq J_k).$$
Nous avons, compte tenu de \eqref{probexceptionnelle},
$$\eqalign{\PP(A_{k}) \leq \sum_{j \leq J_k}
\PP(A_{k,j})  &\ll
\Big(\frac{x_k-x_{k-1}}{x_{k-1}}\Big)^{4/3} 
(\log x_k)^{67/3} \sum_{j \leq J_k} \e^{-3 \const{expoz} j} \cr &\ll \Big(\frac{x_k-x_{k-1}}{x_{k-1}}\Big)^{4/3} 
(\log x_k)^{67/3}.}$$

Soit $0 < \const{xk} < 1$. Posons $x_k =  [ \exp\{k^{\const{xk}}\} ]$.
L'estimation $$\frac{x_k - x_{k-1}}{x_{k-1}} \sim \frac{
\const{xk}}{k^{1-\const{xk}}} \eqdef{xkmoinsxkmoins1}$$ assure tout d'abord que la condition \eqref{petitint} est bien vérifiée. En effet, le fait que $2 \log x_k \leq \vartheta(y_{k,j})$ implique qu'il existe une constante $\dconst{kapp} > 0$  telle que $x_k^{\const{kapp}/j} \geq \log x_k$. 
D'autre part, la relation \eqref{xkmoinsxkmoins1} fournit $$\PP(A_k) \ll \frac{1}{k^{5/4}},$$
en choisissant par exemple $\const{xk}= 1/356$.

Le lemme de Borel-Cantelli permet alors d'obtenir 
$$\PP\Big(\limsup_{k \geq 1} A_k\Big) = 0$$
puisque $\sum_{k \geq 1} \PP(A_k) < \infty.$
Ainsi, presque sûrement, il existe un entier $k_0$ tel que,
pour tout $k \geq k_0$,  et pour tout $j \leq J_k$,
$$\max_{\di{]a,b] \in \D_k}{y \leq y_{k,j}}}
\big|\Psif{f}{b,y}-\Psif{f}{a,y}\big| \leq R_{k,j}.$$

\qed

\Propc{cordulem1}{Sous les mêmes hypothèses, nous avons
pour $k \geq 1$ et $j \leq J_k$,
$$\max_{\di{x_{k-1}\leq x \leq x_k}{y \leq y_{k,j}}}
\big|\Psif{f}{x,y} - \Psif{f}{x_k,y}\big|
\ll \sqrt{\Psisf{x_k,y_{k,j}}}\ \e^{\const{expoz}\, j} \quad \ps.$$}

\nid 
Soit $K \geq 0$. Tout entier positif $n \leq 2^{K+1}-1$ peut s'écrire
$$ n = \sum_{j \leq K} a_j(n) 2^j \qquad \big(a_j(n) \in \{0, 1\}, \forall j \in [0,K] \big). \eqdef{binaire}$$
Pour tout $k \geq 2$, le choix $K=\log(x_k - x_{k-1}+1)/(\log 2) - 1 $ est donc suffisant pour pouvoir écrire tout entier $n \leq x_k - x_{k-1}$ sous la forme \eqref{binaire}.
Cela implique, pour $x_{k-1}< x \leq x_k$, que tout intervalle $]x_{k-1}, x]$ peut s'écrire comme la réunion disjointe d'au plus $O(\log x_k)$
intervalles dyadiques de $\D_k$, où $\D_k$ a été défini en \eqref{defDk}.

Désignons alors par $\D(]x_{k-1}, x])$ l'ensemble des
intervalles dyadiques apparaissant dans la
décomposition de $]x_{k-1}, x]$ dont le
nombre de termes est minimal. 

Nous avons alors, pour tout $x_{k-1} < x \leq x_k$,
$$\eqalign{\big| \Psif{f}{x,y} -
\Psif{f}{x_{k-1},y} \big| &\leq \sum_{{]a,b] \in
\D(]x_{k-1},x])}} \big|\Psif{f}{b,y}-\Psif{f}{a,y}\big| \cr
&\ll (\log x_k) \max_{\di{]a,b] \in \D_k}{y \leq y_{k,j}}}
\big|\Psif{f}{b,y}-\Psif{f}{a,y}\big|.}$$

On conclut par le lemme précédent.
\qed

\bigskip
\bigskip\goodbreak
\paradeuxn{Moyenne de $\Psif{f}{x,y}$ autour de points-tests}
Nous pouvons à présent démontrer le résultat suivant.
\Propl{lem}{
Soit $f$ une fonction multiplicative aléatoire au sens de Wintner. Soit
$\{x_k\}_{k \geq 1}$ définie comme au \ref{lem1}. Pour tout $\epsilon >0$, il existe une constante $\dconst{m18}>0$, telle que, uniformément pour $$\exp\Big\{ \const{m18} \frac{\log x_k}{\log_2 x_k} \Big\}\leq y \leq x_k, \qquad k \geq k_0,$$
nous ayons
$$ \frac{1}{x_k-x_{k-1}}\int_{x_{k-1}}^{x_k}\, \Psif{f}{x,y} \d x  \ll \sqrt{x_k}\, (\log_2 x_k)^{2+\epsilon} \qquad \ps.$$}

\medskip
\nid 
Il est nécessaire d'introduire tout d'abord un certain nombre de paramètres qui seront fixés ultérieurement. Soient donc $\gamma > 0$, $X \geq 2$ et $z_0 \geq 2$. 
Posons de plus $$k_X := \Big\{ k \geq 1 : \sqrt{X} < x_k \leq X \Big\}, \qquad h := \frac{\log X}{\log z_0}, \eqdef{defdeh}$$
en notant d'emblée que $|k_X| \leq (\log X)^{1/\const{xk}}.$  Enfin, soit 
$\{z_j\}_{j\geq 0}$ une suite de nombres réels positifs définie de telle sorte que $z_j := z_{j-1}^{1+\gamma} = {z_0}^{(1+\gamma)^j}$ $(j \geq 1)$.
\medskip
Nous pouvons, lorsque $2 \leq x \leq X$ et $y > z_0$, décomposer la somme $\Psif{f}{x,y}$ selon la taille de $P^+(n)$ de la façon suivante 
$$\Psif{f}{x,y} = \Psif{f}{x,z_0} + \sum_{j \leq J_0(y)-1} \sum_{\di{n \leq x}{P^+(n) \in ]z_j,z_{j+1}] }} f(n) + \sum_{\di{n \leq x}{P^+(n) \in ]z_{J_0},y] }} f(n),$$
où l'on a posé $J_0 = J_0(y) := \max \{j : z_j \leq y\}$ en observant que $J_0(y) \leq J := \min \{j : z_j \geq X\}$. Cela implique
$$\eqalign{\Psif{f}{x,y} &= \Psif{f}{x,z_0} + \sum_{j \leq J_0(y)-1} \sum_{1 < r\leq x}^{z_j,z_{j+1}} f(r) \Psifbig{f}{\frac{x}{r},z_j} \cr &\qquad \qquad \qquad \qquad+ \sum_{1 < r\leq x}^{z_{J_0},y} f(r) \Psifbig{f}{\frac{x}{r},z_{J_0}}.}$$

\medskip
Nous pouvons à présent écrire, pour la suite $\{x_k\}_{k\geq 1}$ définie au \ref{lem1},
 $$\frac{1}{x_k-x_{k-1}} \Big| \int_{x_{k-1}}^{x_k} \Psif{f}{x,y} \Big| \d x \leq V_k + \sum_{j\leq J} W_{k,j}^*,
 \eqdef{decoup1}$$
avec
$$\eqalign{V_k &:= \frac{1}{x_k-x_{k-1}} \Big| \int_{x_{k-1}}^{x_k} \Psif{f}{x,z_0} \d x\Big|, \cr W_{k,j}(z) &:= \frac{1}{x_k-x_{k-1}} \sum_{1 < r\leq x_k}^{z_j,z} f(r) \int_{x_{k-1}}^{x_k} \Psifbig{f}{\frac{x}{r},z_j} \d x,}$$
et $$W_{k,j}^* := \max_{z_j < z \leq z_{j+1}} |W_{k,j}(z)|, \qquad Z_{k,j} := W_{k,j}(z_{j+1}).$$

Soit $R$ un paramètre réel positif qui sera fixé ultérieurement (voir \eqref{choixdeR}), et considérons les événements
$$\eqalign{A &:= \bigcup_{k \in k_X} \Big\{  \frac{1}{x_k-x_{k-1}} \max_{z_0 < y \leq x_k} \Big|  \int_{x_{k-1}}^{x_k} \Psif{f}{x,y} \d x\Big| \geq 3 \sqrt{x_k} R  \Big\}, \cr
B &:= \bigcup_{k \in k_X} \Big\{ V_k \geq \sqrt{x_k} R \Big\}, \cr
C_j &:= \bigcup_{k \in k_X} \Big\{ W_{k,j}^* \geq \sqrt{x_k} R / J \Big\}, \qquad
C := \bigcup_{j\leq J} C_j.
}$$

La relation \eqref{decoup1} implique directement $$A \subset B \cup C. \eqdef{decoupsubset}$$

Voyons comment majorer $\PP(B)$. Il vient, par l'inégalité de Cauchy-Schwarz,
$$\EE\, ({V_k} ^2) \leq 
\frac{1}{x_k - x_{k-1}} \int_{x_{k-1}}^{x_k} \EE \, \Psif{f}{x,z_0}^2 \d x \leq \Psi(x_k,z_0).$$
car, d'après la relation \eqref{momentordre2},
 $$\EE\, \Psif{f}{x,z_0}^2 = \Psi^*(x_k,z_0) \leq \Psi(x_k,z_0).$$
L'inégalité de Markov et la majoration \eqref{typerankin} pour $\Psi(x,y)$ impliquent
 $$\PP\big(V_k \geq \sqrt{x_k} R\big) \leq \frac{\Psisf{x_k,z_0}}{x_k R^2} \ll 
\frac{\e^{- \log x_k / (2 \log z_0)}}{R^2},$$
et
 $$\PP(B) \ll (\log X)^{1/\const{xk}} \frac{\e^{- h/4}}{R^2} \cdot \eqdef{by}$$

Il nous reste donc à majorer $\PP(C)$. Donnons-nous un entier pair $\ell\geq 2$ 
et introduisons la quantité
$$D_{k,j} := \frac{1}{x_k} \sum_{1 < r\leq x_k}^{z_j,z_{j+1}} \mu(r)^2 (\ell-1)^{\Omega(r)} \Big\{\frac{1}{x_k - x_{k-1}}\int_{x_{k-1}}^{x_k} \Psifbig{f}{\frac{x}{r},z_j} \d x\Big\}^2.$$
Nous avons, par l'inégalité de Cauchy-Schwarz, pour tout $k\in k_X$,
$$\eqalign{D_{k,j} &\leq \frac{1}{x_k} \sum_{1 < r\leq x_k}^{z_j,z_{j+1}} 
\mu(r)^2 (\ell-1)^{\Omega(r)}
\Big\{\frac{1}{x_k - x_{k-1}} \int_{x_{k-1}}^{x_k} \Psifbig{f}{\frac{x}{r},z_j}^2 \d x\Big\}
\cr&\leq
\sum_{1 < r\leq x_k}^{z_j,z_{j+1}}
\mu(r)^2 (\ell-1)^{\Omega(r)}
\Big\{\frac{r/x_k}{x_k - x_{k-1}} \int_{x_{k-1}/r}^{x_k/r} \Psif{f}{v,z_j}^2 \d v\Big\}
\cr&\leq
\int_{1}^{x_k/z_j} \frac{\Psif{f}{v,z_j}^2}{v} 
\Big\{\frac{1}{x_k - x_{k-1}}
\sum_{x_{k-1}/v < r\leq x_k/v}^{z_j,z_{j+1}} 
\mu(r)^2 (\ell-1)^{\Omega(r)} \Big\} \d v \cr
&\leq 
\frac{\dconst{c10} \ell \e^{2\gamma \ell}}{\log z_j}
\int_{1}^{x_k/z_j} \frac{\Psif{f}{v,z_j}^2}{v^2} 
\d v \leq \frac{\const{c10} \ell \e^{2\gamma \ell}}{\log z_j} N_j.
}$$
où l'on utilise la notation $N_j$ définie en \eqref{defNj}, et où la dernière somme en $r$ a été estimée par le \ref{lemsommecourte}, 
pourvu que l'on ait, pour $\delta > 0$ convenablement fixé,
 $$z_0 > \Big(\frac{x_k}{x_k - x_{k-1}}\Big)^{1+\delta} \eqdef{condz}$$
et
$$ \gamma \geq \frac{\dconst{calpha1}}{L(z_0)^{\const{clem36}}} \cdot \eqdef{condgamma}$$ 

Notons que la condition \eqref{condz} est impliquée par l'inégalité
$$h \leq \dconst{zzero} \frac{\log X}{\log_2 X}, \eqdef{condh}$$
où $h$ a été défini en \eqref{defdeh}.
Les assertions \eqref{condh} et \eqref{condgamma} seront respectivement vérifiées en
\eqref{choixdeh}  et \eqref{verifgamma} lors du choix des paramètres $h$ et $\gamma$.

\medskip
Soit à présent $R'$ un paramètre réel positif fixé ultérieurement. Posons $$F_j :=  \{ N_j \leq R' \cdot \log z_j \}.$$ 

Nous sommes maintenant en mesure d'évaluer $\PP(C)$. Nous avons
$$C_j = \big(C_j \cap F_j\big) \cup \big(C_j \cap \overline{F_j}\big)$$
d'où 
$$\eqalign{C =\bigcup_{j\leq J} C_j &= \bigcup_{j\leq J} \big(C_j \cap F_j\big) \cup \bigcup_{j\leq J} \big(C_j \cap \overline{F_j}\big) \cr
&\subset 
\bigcup_{j\leq J} \big(C_j \cap F_j\big) \cup \bigcup_{j\leq J} \overline{F_j}
} \eqdef{csubset} $$

Grâce à l'inégalité de Doob appliquée à la sous-martingale
$\{W_{k,j}(z)^\ell\}_{z > z_j}$
  nous pouvons écrire, puisque $F_j \in \FF_{z_j}$,
$$\eqalign{\PP (C_j \cap F_j) &\leq \sum_{k \in k_X} \PP \Big(\big\{W_{k,j}^* \geq \sqrt{x_k} R / J \big\} \cap F_j\Big) \cr 
&\leq \sum_{k \in k_X} \frac{\EE\ \Big( \1_{F_j} {Z_{k,j}}^\ell \Big)}{{x_k}^{\ell /2} (R/J)^{\ell}} \leq \sum_{k \in k_X} \frac{\EE\ \Big( \1_{F_j} \EE\ \Big( {Z_{k,j}}^\ell \ | \ \FF_{z_j} \Big) \Big)}{{x_k}^{\ell /2} (R/J)^{\ell}} \cdot}$$

Une simple application du \ref{bonami2} fournit l'évaluation de l'espérance conditionnelle précédente sous la forme
$$\EE\, \big( {Z_{k,j}}^\ell \, | \,  \FF_{z_j} \big) \leq
(x_k \cdot D_{k,j})^{\ell /2}.$$

Il en résulte alors
$$
\eqalign{\PP (C_j \cap F_j) &\leq \sum_{k \in k_X} \EE\, \Big(\1_{F_j} \cdot \frac{{D_{k,j}}^{\ell /2}}{(R/J)^\ell}\Big) \cr
 &\leq (\log X)^{1/\const{xk}} \e^{\gamma \ell^2} \Big(\frac{\const{c10} \cdot \ell \cdot R' \cdot J^2}{R^2}\Big)^{\ell/2}   \cdot
}$$

D'autre part, le \ref{lem_mart} fournit
$$ \PP\Big(\bigcup_{j \leq J} \overline{F_j} \Big) = \PP \Big( \max_{j \leq J} \frac{N_j}{\log z_j} > R' \Big)\ll \frac{\log  h}{R'}\cdot$$

Il découle de \eqref{decoupsubset}, \eqref{by}, \eqref{csubset} et du fait que
$ J \asymp {\log  h}/{\gamma}$, que
$$\eqalign{\PP(A) &\leq \PP(B) + \PP(C) \cr &\ll 
(\log X)^{1/\const{xk}} \Big\{J^{\ell+1} \Big(\frac{\const{c10} \cdot \ell \cdot R'}{R^2}\Big)^{\ell/2} \e^{\gamma \ell^2}
+ \frac{\e^{- h / 4}}{R^2}
\Big\} + \frac{\log  h}{R'} \cr
&\leq 
(\log X)^{1/\const{xk}} \Big\{\frac{(\log  h)^{\ell+1}}{\gamma^{\ell+1}} \Big(\frac{\const{c10} \cdot \ell \cdot R'}{R^2}\Big)^{\ell/2} \e^{\gamma \ell^2}
\!+\!\frac{\e^{- h / 4}}{R^2}
\Big\}\!+\! \frac{\log  h}{R'} \cdot
}$$

Prenons $\gamma := 1 / \ell$. Il vient
$$\PP(A) \leq  
(\log X)^{\dconst{xk2}} \Big\{\Big(\frac{\sqrt{R' \log h}}{R}\Big)^{\ell} 
(\dconst{c215} \ell)^{3 \ell/2+1} 
+ \frac{\e^{-h/4}}{R^2}
\Big\} + \frac{\log  h}{R'} \cdot$$

Le choix $$\ell = 2 \Big[ \dm \Big(\frac{R}{(\sqrt{R'} \log  h)}\Big)^{2/3}\Big]$$ 
conduit à l'estimation
$$\PP(A) \ll 
(\log X)^{\const{xk2}} \Big\{\exp\Big\{-\frac{\dconst{c13} R^{2/3}}{R'^{1/3} (\log  h)^{2/3}} \Big\} +
\frac{\e^{- h / 4}}{R^2}
\Big\} + \frac{\log  h}{R'} \cdot$$

Remarquons qu'en fixant $$h=\dconst{macst118} \log_2 X, \eqdef{choixdeh}$$ nous avons 
$$\PP(A) \ll 
(\log X)^{\const{xk2}} \exp\Big\{- \frac{\const{c13} R^{2/3}}{(R'^{1/3} (\log_3 X)^{2/3})} \Big\} +
\frac{1}{R^2}
+ \frac{\log_3 X}{R'},$$
et que la condition \eqref{condh} est bien remplie, comme annoncé précédemment. Le choix $R' = \dconst{c15}  (\log_2 X)^{1+\epsilon}$ fournit enfin
$$\PP(A) \ll 
(\log X)^{\const{xk2}} \exp\Big\{-\frac{\const{c13} R^{2/3}}{(\log_2 X)^{1/3 + \epsilon}}\Big\} +
\frac{1}{R^2}
+ \frac{1}{(\log_2 X)^{1+\epsilon/2}} \cdot$$

Nous pouvons ainsi fixer $$R = (\log_2 X)^{2 + \epsilon}, \eqdef{choixdeR}$$ en notant que dans ce cas, la condition \eqref{condgamma} est bien vérifiée puisque
$$\gamma \gg \Big(\frac{\sqrt{R'} \log h}{R}\Big)^{2/3} \gg \frac{1}{\sqrt{\log X}} \gg \frac{1}{L(z_0)^{\const{clem36}}}\cdot \eqdef{verifgamma}$$
 Posant $X = X_s := \exp \{ 2^s \}$, nous avons clairement
$$\PP(A_{X_s}) \leq \frac{1}{s^{1+\epsilon}}.$$

Le lemme de Borel-Cantelli permet alors d'établir, puisque la série $\sum_{s \geq 1} \PP(A_{X_s})$ est convergente, que
$$\PP \Big(\limsup_{s \geq 1} A_{X_s} \Big) = 0.$$

Il existe ainsi, presque sûrement, un entier $k_1$ tel que pour $k \geq k_1$, 
$$
\max_{z_0 < y \leq x_k}
 \frac{1}{x_k-x_{k-1}} \int_{x_{k-1}}^{x_k} \Psif{f}{x,y} \d x \leq 3 \sqrt{x_k} R.$$
\qed

\bigskip
\bigskip
\paraunn{Preuve du \ref{mainthm}}
\bigskip

Pour $x$ fixé, nous définissons $k_1$ comme l'unique entier tel que $$x_{k_1-1} < x \leq x_{k_1}, \eqdef{encadrement}$$ 
Nous scindons la démonstration du \ref{mainthm} en plusieurs cas, selon la taille relative de $y$ par rapport à $x$. 

\bigskip\goodbreak
{\it Premier cas.} 

Supposons $x \geq y \geq \exp\{ \const{m18}\,{\log x}/{\log_2 x} \}$. 
Nous avons, dans ce domaine,
$$\eqalign{|\Psif{f}{x,y}| &\leq \!\!\min_{x_{k_1-1} < x' \leq x_{k_1}} | \Psif{f}{x',y} |+\!\!\! \max_{x_{k_1-1} < x' \leq x_{k_1}} | \Psif{f}{x_{k_1},y} - \Psif{f}{x',y}| \cr
& \leq \frac{1}{x_{k_1}-x_{k_1-1}} \int_{x_{k_1-1}}^{x_{k_1}} \Psif{f}{x',y} \d x' \cr &\qquad \qquad +  \max_{x_{k_1-1} < x' \leq x_{k_1}} | \Psif{f}{x_{k_1},y} - \Psif{f}{x',y}|.}$$
Il vient alors 
$${ \Psif{f}{x,y} \ll \sqrt{x_{k_1}}\, (\log_2 x_{k_1})^{2+\epsilon} + x_k \ll \sqrt{x}\, (\log_2 x)^{2+\epsilon},} \eqdef{premcas}$$
en appliquant respectivement le \ref{cordulem1} et le \ref{lem}.

\bigskip\goodbreak
{\it Deuxième cas.} 

Nous supposons ici que l'on a $\theta(y) > 2 \log x$.
Posons à présent
$$j_k := \max \Big\{1, \Big[ \frac{\log x_k}{\log y} \Big] \Big\}, \qquad y_k := x_k^{1/j_k} \qquad (k \leq k_1).$$
 Le \ref{cordulem1}, fournit directement
$$|\Psif{f}{x_k,y} - \Psif{f}{x_{k-1},y}| \leq \sqrt{\Psisf{x_k,y_k}} \ \e^{\const{expoz}\, j_k} \qquad  (k  \leq k_1)$$
et
$$|\Psif{f}{x_{k_1},y} - \Psif{f}{x,y}| \leq \sqrt{\Psisf{x_{k_1},y_{k_1}}}\  \e^{\const{expoz}\, j_{k_1}}.$$

Notons que, pour tout $k \leq k_1$, nous avons
$$\eqalign{j_k &\leq \max \Big\{1, \frac{\log x_k}{\log y}\Big\}  \leq \frac{\log x_{k_1}}{\log y} \leq 
\frac{\log x_{k_1-1}}{\log y} + \frac{\log (x_{k_1}/x_{k_1-1}) }{\log y} \cr &
\leq u + O\Big(\frac 1 {\log y}\Big),}$$
en vertu de la relation \eqref{encadrement}. 
Supposant provisoirement l'estimation
$$M_1:=\sup_{k \leq k_1} \frac{\Psisf{x_k,y_k}}{\Psisf{x,y}} = \exp \Big\{
O \Big( \frac{u}{L_\epsilon(y)} + \log(u+1) \Big) \Big\} \eqdef{m11}$$
acquise, nous pouvons en déduire la majoration
$$\eqalign{|\Psif{f}{x,y}| &\leq \sum_{k \leq k_1} | \Psif{f}{x_k,y} - \Psif{f}{x_{k-1},y} |
+ |\Psif{f}{x,y} - \Psif{f}{x_{k_1},y}| \cr &\ll
k_1 \sqrt{M_1} \sqrt{\Psisf{x,y}} \e^{\const{expoz} u} \ll
\sqrt{\Psisf{x,y}}\, \e^{\dconst{aa1277} u}\, (\log x)^{\dconst{aa1276}}
,} \eqdef{seccas}$$
où l'on a utilisé le fait que $k_1 \ll (\log x)^{1/\const{xk}}$.

\'Etablissons à présent la relation \eqref{m11}.
Lorsque $x$ est suffisamment grand et $(\log x)^3 < y \leq x$, nous avons $\beta(x,y) \geq 3/5$. En appliquant le corollaire 2.6 de \citer{BT05b}
et le théorème III.5.21 de \citer{T08}, il vient
 $$\eqalign{M_1 &\asymp \sup_{k \leq k_1} \frac{\Psi(x_k,y_k)}{\Psi(x,y)} \cr &\leq \sup_{k \leq k_1} \frac{x_{k}}{x} \frac{\rho(\log x_k/ \log y -1)}{\rho(u)} \exp \Big\{
O \Big(\frac{\log (u+1)}{\log y}  + \frac{u}{L_\epsilon(y)} + \frac 1 u \Big) \Big\}.}$$
Faisant usage de la croissance de la fonction $x \mapsto x\, \rho(\log x / \log y - 1)$ pour $x$ suffisamment grand et $y \leq x \leq \exp\{y^{1/3}\}$, nous obtenons
$$\eqalign{M_1 &\leq
\frac{x_{k_1}}{x} \frac{\rho(\log x_{k_1}/ \log y -1)}{\rho(u)} \exp \Big\{
O \Big(\frac{\log (u+1)}{\log y}  + \frac{u}{L_\epsilon(y)} + \frac 1 u \Big) \Big\} \cr
&= \exp\Big\{O\Big(\frac{u}{L_\epsilon(y)} + \log (u+1) \Big)\Big\}.}$$ 

Plaçons-nous à présent dans le domaine $$ 2 \log x < \vartheta(y) \leq (\log x)^{3}.$$
Puisque $j_k \geq {\log x_k}/{\log y} -1$, nous avons
  $y \leq y_k \leq y + y^{3/4}$ 
si l'on suppose de plus que
$$ x_k \geq \exp \{y^{1/4} (\log y)^2 \}. \eqdef{condk}$$

En remarquant que nous avons trivialement
$$\Psisf{x_k,y_k} \leq x_k \leq
\exp\Big\{O\Big(\frac{u}{L_\epsilon(y)}\Big)\Big\},$$
pour tous les entiers $k$ ne vérifiant pas \eqref{condk},  et que $x_k \leq \e\, x$ $(k \leq k_1)$,
il vient
$$\eqalign{M_1 &\ll \frac{\Psisf{x, y + y^{3/4}}}{\Psisf{x,y}} = 1+
\sum_{\di{d>1}{P^-(d)>y, P^+(d) \leq y+ y^{3/4}}} \mu(d)^2 \frac{\Psisf{x/d,y}}{\Psisf{x,y}} \cr &\ll 1 + 
\sum_{\di{d>1}{P^-(d)>y, P^+(d) \leq y+y^{3/4}}} \frac{\mu(d)^2}{d^\alpha}
 \ll \exp \Big\{ \sum_{y < p \leq y+ y^{3/4}} \frac 1{p^{\alpha}}\Big\},}$$
où $\alpha$ a été défini en \eqref{def-alpha}. Posant $$w(t) = w(t;x,y) := \frac{t^{1 - \alpha} - 1}{(1-\alpha) \log t} \qquad (t \geq 2),$$ le lemme 3.6 de \citer{BT05b} fournit, pour $y \geq 2$,
$$\eqalign{\sum_{y < p \leq y + y^{3/4}}\frac{1}{p^{\alpha}} &=
\log \Big(\frac{\log (y+y^{3/4})}{\log y}\Big) + \int_{w(y)}^{w(y+y^{3/4})}\!\!\!\!\!\!\!\!\!\! t \xi'(t) \ d t + 
O\Big(\frac{y^{1-\alpha}}{L(y)^{\dconst{cstLy}}}\Big) \cr
&\ll w(y+y^{3/4}) - w(y) + \frac{y^{1-\alpha}}{L(y)^{\const{cstLy}}}
\ll \frac{y^{3/4}}{y^{\alpha} \log y} + \frac{y^{1-\alpha}}{L(y)^{\const{cstLy}}} \cr
&\ll \frac{u}{L_\epsilon(y)},}$$
car on a $y^{1-\alpha} \asymp \log x$ dans le domaine considéré.
Cela achève la démonstration de l'estimation \eqref{m11} et donc de \eqref{seccas}.

Un calcul simple permet enfin d'observer que l'estimation \eqref{seccas} est meilleure que celle résultant de \eqref{premcas} si
$$x \geq 16, \qquad y \leq \exp\Big\{\dconst {meil1} \frac{\log x \, \log_3 x}{\log_2 x}\Big\}.$$
où $\const {meil1}$ est une constante positive convenable.

\bigskip
{\it Troisième cas.}

Il s'agit à présent d'évaluer le domaine dans lequel la majoration triviale
$$|\Psif{f}{x,y}| \leq \Psisf{x,y} \eqdef{majtriv}$$ est meilleure que l'estimation \eqref{seccas}. Nous supposons ici que $\theta(y) > 2 \log x$.

D'une part, pour tout $\epsilon > 0$, il existe $\delta > 0$ tel que
l'on ait $$v := \log (\theta(y) / \log x -1 ) \leq \epsilon \log y, \eqdef{def_v}$$ lorsque
$y \leq (\log x)^{1+\delta}$.
 Aussi, sous cette dernière hypothèse, le corollaire 2.3 de \citer{BT02} fournit la minoration 
$$ \Psisf{x,y} \geq \exp \{ (1-\epsilon) \pi(y) v / (\e^v +1 ) \} 
\geq \exp \{ (1-\epsilon') v u \}, \eqdef{minorPsisf}$$
d'où
$$\Psisf{x,y} \gg \sqrt{\Psisf{x,y}} \, \e^{\const{expos} u } (\log x)^{\const{exposant}}$$
lorsque $v$ est suffisamment grand.

D'autre part, lorsque $y \geq (\log x)^{1+\delta}$, il découle
de l'estimation de $\Psi(x,y)$ par Hildebrand et Tenenbaum (théorème III.5.21 de \citer{T08}) que 
$$\sqrt{\Psisf{x,y}} \asymp \sqrt{\Psi(x,y)} \gg \e^{\const{expos} u } (\log x)^{\const{exposant}}$$
et donc que l'estimation \eqref{seccas} est toujours de qualité supérieure à la majoration triviale \eqref{majtriv} lorsque $y \geq (\log x)^{1+\delta}$ et $x$ est suffisamment grand.

Au vu des observations précédentes, la majoration fournie par \eqref{seccas} est meilleure que la majoration triviale \eqref{majtriv} lorsque $\theta(y) \geq \dconst{aa991} \log x$, où $\const{aa991}$ est une constante positive.

\qed
\bigskip
\bigskip
\paraunn{Preuve du \ref{cor1}}
\bigskip

Déterminons dans quel domaine les majorations apportées par le \ref{mainthm} fournissent une estimation du type $\Psif{f}{x,y} \ll_\epsilon \Psi^*(x,y)^{1/2+\epsilon}$ pour tout $\epsilon > 0$.
Tout d'abord, il apparaît clairement (en utilisant par exemple la formule asymptotique issue du corollaire III.5.19 de \citer{T08}) que l'on a
$$\sqrt{x}\, (\log_2 x)^{2+\epsilon'} \ll \Psi(x,y)^{1/2+\epsilon} \asymp \Psi^*(x,y)^{1/2+\epsilon}$$
lorsque
 $x \geq y \geq x^{\const{m17}\,(\log_3 x)/ \log_2 x}$ et $x$ suffisamment grand.
 
Par ailleurs, pour 
$(\log x)^{1+\delta} \leq y \leq x^{\const{m17}\,(\log_3 x)/ \log_2 x}$ et $x$ suffisamment grand, le théorème III.5.21 de \citer{T08} nous assure que

$$\Psisf{x,y}\  \e^{\const{expos} u } (\log x)^{\const{exposant}} \ll \Psi^*(x,y)^{1/2+\epsilon}$$
 puisque l'on a alors
 $$\Psi^*(x,y)^\epsilon \asymp \Psi(x,y)^\epsilon \gg \e^{\const{expos} u } (\log x)^{\const{exposant}}.$$

Enfin, la minoration \eqref{minorPsisf},
valable dans le domaine $y \leq (\log x)^{1+\delta}$, $\theta(y) > 2 \log x$
nous assure que
 $$\Psi^*(x,y)^\epsilon \gg \e^{\const{expos} u } (\log x)^{\const{exposant}}$$
lorsque $\epsilon\, v$ est suffisamment grand (rappelons que $v$ est défini en \eqref{def_v}), c'est-à-dire lorsque $y \geq C_\epsilon \log x$.
Ceci implique directement la majoration 
$$\Psisf{x,y}\  \e^{\const{expos} u } (\log x)^{\const{exposant}} \ll \Psi^*(x,y)^{1/2+\epsilon}.$$\qed
\bigskip
{\bf Remerciements.}
L'auteur remercie Gérald Tenenbaum pour l'aide qu'il lui a apportée au cours de l'élaboration de ce travail, et tient également à exprimer sa gratitude à l'arbitre, pour ses remarques éclairantes.

\bigskip\bigskip\bigskip\goodbreak

\centerline{\twelvebf Bibliographie}
\bigskip
\eightpoint\leftskip9mm
\bibtem{B70} A. Bonami, Étude des coefficients de Fourier des fonctions de $L^p(G)$, {\it Ann. Inst. Fourier} {\bf 20} (1970), 335--402.
\bibtem{BT02} R. de la Bretèche, G. Tenenbaum, Sur les lois locales de la répartition du k-ième diviseur d'un entier, {\it Proc. London Math. Soc.} (3) {\bf 84} (2002).
\bibtem{BT05b} R. de la Bretèche, G. Tenenbaum, Propriétés statistiques des entiers friables, {\it Ramanujan J.} {\bf 9} (2005), 139--202.
\bibtem{E61} P. Erd\H os, Some unsolved problems, {\it Magyar Tud. Akad. Kutat\'o Int. Közl.} {\bf 6} (1961), 221--254.
\bibtem{H82} G. Hal\'asz, On random multiplicative functions, Hubert Delange colloquium (Orsay, 1982), {\it Publ. Math. Orsay} {\bf 83} (4) (1983), 74--96. 
\bibtem{H10} A. J. Harper, On the limit distributions of some sums of a random multiplicative function, à paraître.
\bibtem{HT86} A. Hildebrand, G. Tenenbaum, On integers free of large prime factors, {\it Trans. Amer. Math. Soc.} {\bf 296} (1986), 265--290.
\bibtem{HT93} A. Hildebrand, G. Tenenbaum, Integers without large prime factors. {\it J. Théor. Nombres Bordeaux}  {\bf 5}  (1993),  no. 2, 411--484.
\bibtem{LTW10} Y.-K. Lau, G. Tenenbaum, J. Wu, Mean values of random multiplicative functions, prépublication, sept. 2010 ; version finale :  {\it Proc. Amer. Math. Soc.}, à paraître.
\bibtem{T08} G. Tenenbaum, {\it Introduction à la théorie analytique et probabiliste des nombres}, Troisième Edition, Belin, 2008.
\bibtem{TW03} G. Tenenbaum, J. Wu, Moyennes de certaines fonctions multiplicatives sur les entiers friables, {\it J. Reine Angew. Math.} {\bf 564} (2003), 119--166.
\bibtem{TL06} T. T\'om\'acs, Z. L\'\i bor, A H\'ajek-Renyi type inequality and its applications, {\it Ann. Math. et Inform.} {\bf 33} (2006), 141--149.
\bibtem{W91} D. Williams, {\it Probability with martingales}, Cambridge University Press, 1991.
\bibtem{W44} A. Wintner, Random factorizations and Riemann's hypothesis, {\it Duke Math. J. } {\bf 11} (1944), 267--275.
\par

\vskip 5mm
{\sevenrm\baselineskip9pt\obeylines
Joseph Basquin
Institut \'Elie Cartan
Universit\'e Henri Poincar\'e--Nancy 1
   BP 239
54506 Vand\oe uvre Cedex
   France}
\smallskip
internet: \seventt joseph.basquin@iecn.u-nancy.fr\par
\goodbreak
\bye